\documentclass[11pt]{article}
\usepackage{times,amsmath,amssymb,graphicx,theorem}%,wrapfig,graphicx,epic,eepic}%multicol
\usepackage{pstricks,pst-node,pst-text,pst-tree,pstricks-add}% pstricks-add MUST come last
\makeatletter
\def\section{\@startsection{section}{1}{0pt}{-3.25ex plus -1ex minus 
-.2ex}{1.5ex plus .2ex minus .3ex}{\normalfont\large\bf}}
\renewcommand\subsection{\@startsection{subsection}{2}{\z@}%
                                     {-3.25ex\@plus -1ex \@minus -.2ex}%
                                     {1.5ex \@plus .2ex}%
                                     {\normalfont\normalsize\bfseries}}
\makeatother
\renewenvironment{abstract}
{\vspace{-1.8ex}\begin{quotation}\small%\textbf{Abstract. }
}{\end{quotation}}
\newcommand{\parag}[1]{\paragraph{#1.}\!\!\!}
\newcommand{\defn}[1]{{\textit{\textbf{#1}}}}
\newcommand{\myitem}[1]{\item[\textnormal{(#1)}]}

\theoremheaderfont{\scshape}
\theorembodyfont{\normalfont\slshape}
\theoremstyle{plain}

\newtheorem{lemma}{Lemma}
\newtheorem{proposition}{Proposition}
\newtheorem{corollary}{Corollary}
\newtheorem{theorem}{Theorem}
\newenvironment{proof}{\begin{trivlist}\item{}\normalfont\textit{Proof.}}{\hfill$\square$\end{trivlist}}
\newenvironment{proofof}[1]{\begin{trivlist}\item{}\normalfont\textit{Proof of #1.}}{\hfill$\square$\end{trivlist}}

\newcommand{\nth}{^{\text{th}}}

\newcommand{\ie}{\emph{i.e.}}
\newcommand{\eg}{\emph{e.g.}}
\newcommand{\cf}{\emph{cf.}}
\newcommand{\etc}{\emph{etc.}}

\newdimen\arrayruleHwidth
\makeatletter
\def\Hline{\noalign{\ifnum0=`}\fi\hrule \@height \arrayruleHwidth
  \futurelet \@tempa\@xhline}
\makeatother

\newdimen\proofrulebreadth \proofrulebreadth=.05em
\newdimen\proofdotseparation \proofdotseparation=1.25ex
\newdimen\proofrulebaseline \proofrulebaseline=2ex
\newcount\proofdotnumber \proofdotnumber=3
\let\then\relax
\def\hfi{\hskip0pt plus.0001fil}
\mathchardef\squigto="3A3B
\newif\ifinsideprooftree\insideprooftreefalse
\newif\ifonleftofproofrule\onleftofproofrulefalse
\newif\ifproofdots\proofdotsfalse
\newif\ifdoubleproof\doubleprooffalse
\let\wereinproofbit\relax
\newdimen\shortenproofleft
\newdimen\shortenproofright
\newdimen\proofbelowshift
\newbox\proofabove
\newbox\proofbelow
\newbox\proofrulename
\def\shiftproofbelow{\let\next\relax\afterassignment\setshiftproofbelow\dimen0 }
\def\shiftproofbelowneg{\def\next{\multiply\dimen0 by-1 }%
\afterassignment\setshiftproofbelow\dimen0 }
\def\setshiftproofbelow{\next\proofbelowshift=\dimen0 }
\def\setproofrulebreadth{\proofrulebreadth}

\def\prooftree{% NESTED ZERO (\ifonleftofproofrule)
\ifnum  \lastpenalty=1
\then   \unpenalty
\else   \onleftofproofrulefalse
\fi
\ifonleftofproofrule
\else   \ifinsideprooftree
        \then   \hskip.5em plus1fil
        \fi
\fi
\bgroup% NESTED ONE (\proofbelow, \proofrulename, \proofabove,
\setbox\proofbelow=\hbox{}\setbox\proofrulename=\hbox{}%
\let\justifies\proofover\let\leadsto\proofoverdots\let\Justifies\proofoverdbl
\let\using\proofusing\let\[\prooftree
\ifinsideprooftree\let\]\endprooftree\fi
\proofdotsfalse\doubleprooffalse
\let\thickness\setproofrulebreadth
\let\shiftright\shiftproofbelow \let\shift\shiftproofbelow
\let\shiftleft\shiftproofbelowneg
\let\ifwasinsideprooftree\ifinsideprooftree
\insideprooftreetrue
\setbox\proofabove=\hbox\bgroup$\displaystyle % NESTED TWO
\let\wereinproofbit\prooftree
\shortenproofleft=0pt \shortenproofright=0pt \proofbelowshift=0pt
\onleftofproofruletrue\penalty1
}

\def\eproofbit{% NESTED TWO
\ifx    \wereinproofbit\prooftree
\then   \ifcase \lastpenalty
        \then   \shortenproofright=0pt  % 0: some other object, no indentation
        \or     \unpenalty\hfil         % 1: empty hypotheses, just glue
        \or     \unpenalty\unskip       % 2: just had a tree, remove glue
        \else   \shortenproofright=0pt  % eh?
        \fi
\fi
\global\dimen0=\shortenproofleft
\global\dimen1=\shortenproofright
\global\dimen2=\proofrulebreadth
\global\dimen3=\proofbelowshift
\global\dimen4=\proofdotseparation
\global\count255=\proofdotnumber
$\egroup  % NESTED ONE
\shortenproofleft=\dimen0
\shortenproofright=\dimen1
\proofrulebreadth=\dimen2
\proofbelowshift=\dimen3
\proofdotseparation=\dimen4
\proofdotnumber=\count255
}

\def\proofover{% NESTED TWO
\eproofbit % NESTED ONE
\setbox\proofbelow=\hbox\bgroup % NESTED TWO
\let\wereinproofbit\proofover
$\displaystyle
}%
\def\proofoverdbl{% NESTED TWO
\eproofbit % NESTED ONE
\doubleprooftrue
\setbox\proofbelow=\hbox\bgroup % NESTED TWO
\let\wereinproofbit\proofoverdbl
$\displaystyle
}%
\def\proofoverdots{% NESTED TWO
\eproofbit % NESTED ONE
\proofdotstrue
\setbox\proofbelow=\hbox\bgroup % NESTED TWO
\let\wereinproofbit\proofoverdots
$\displaystyle
}%
\def\proofusing{% NESTED TWO
\eproofbit % NESTED ONE
\setbox\proofrulename=\hbox\bgroup % NESTED TWO
\let\wereinproofbit\proofusing
\kern0.3em$
}

\def\endprooftree{% NESTED TWO
\eproofbit % NESTED ONE
  \dimen5 =0pt% spread of hypotheses
\dimen0=\wd\proofabove \advance\dimen0-\shortenproofleft
\advance\dimen0-\shortenproofright
\dimen1=.5\dimen0 \advance\dimen1-.5\wd\proofbelow
\dimen4=\dimen1
\advance\dimen1\proofbelowshift \advance\dimen4-\proofbelowshift
\ifdim  \dimen1<0pt
\then   \advance\shortenproofleft\dimen1
        \advance\dimen0-\dimen1
        \dimen1=0pt
        \ifdim  \shortenproofleft<0pt
        \then   \setbox\proofabove=\hbox{%
                        \kern-\shortenproofleft\unhbox\proofabove}%
                \shortenproofleft=0pt
        \fi
\fi
\ifdim  \dimen4<0pt
\then   \advance\shortenproofright\dimen4
        \advance\dimen0-\dimen4
        \dimen4=0pt
\fi
\ifdim  \shortenproofright<\wd\proofrulename
\then   \shortenproofright=\wd\proofrulename
\fi
\dimen2=\shortenproofleft \advance\dimen2 by\dimen1
\dimen3=\shortenproofright\advance\dimen3 by\dimen4
\ifproofdots
\then
        \dimen6=\shortenproofleft \advance\dimen6 .5\dimen0
        \setbox1=\vbox to\proofdotseparation{\vss\hbox{$\cdot$}\vss}%
        \setbox0=\hbox{%
                \advance\dimen6-.5\wd1
                \kern\dimen6
                $\vcenter to\proofdotnumber\proofdotseparation
                        {\leaders\box1\vfill}$%
                \unhbox\proofrulename}%
\else   \dimen6=\fontdimen22\the\textfont2 % height of maths axis
        \dimen7=\dimen6
        \advance\dimen6by.5\proofrulebreadth
        \advance\dimen7by-.5\proofrulebreadth
        \setbox0=\hbox{%
                \kern\shortenproofleft
                \ifdoubleproof
                \then   \hbox to\dimen0{%
                        $\mathsurround0pt\mathord=\mkern-6mu%
                        \cleaders\hbox{$\mkern-2mu=\mkern-2mu$}\hfill
                        \mkern-6mu\mathord=$}%
                \else   \vrule height\dimen6 depth-\dimen7 width\dimen0
                \fi
                \unhbox\proofrulename}%
        \ht0=\dimen6 \dp0=-\dimen7
\fi
\let\doll\relax
\ifwasinsideprooftree
\then   \let\VBOX\vbox
\else   \ifmmode\else$\let\doll=$\fi
        \let\VBOX\vcenter
\fi
\VBOX   {\baselineskip\proofrulebaseline \lineskip.2ex
        \expandafter\lineskiplimit\ifproofdots0ex\else-0.6ex\fi
        \hbox   spread\dimen5   {\hfi\unhbox\proofabove\hfi}%
        \hbox{\box0}%
        \hbox   {\kern\dimen2 \box\proofbelow}}\doll%
\global\dimen2=\dimen2
\global\dimen3=\dimen3
\egroup % NESTED ZERO
\ifonleftofproofrule
\then   \shortenproofleft=\dimen2
\fi
\shortenproofright=\dimen3
\onleftofproofrulefalse
\ifinsideprooftree
\then   \hskip.5em plus 1fil \penalty2
\fi
}

\newcount\PLv\newcount\PLw\newcount\PLx\newcount\PLy\newdimen\PLyy\newdimen\PLX
\newbox\PLdot \setbox\PLdot\hbox{{\tiny.}} \def\scl{.07} % resettable scale
\def\PLot#1{\PLx`#1\advance\PLx-42\PLy\PLx\PLv\PLx\divide\PLy9\PLw\PLy\multiply
\PLw9\advance\PLx-\PLw\advance\PLx-4\PLy-\PLy\advance\PLy4\PLX=\the\PLx pt
\advance\PLyy\the\PLy pt\wd\PLdot=\scl\PLX\raise\scl\PLyy\copy\PLdot}
\def\draw#1{\ifx#1\end\let\next=\relax\else\PLot#1\let\next=\draw\fi\next}

\newcount\smallPLv\newcount\smallPLw\newcount\smallPLx\newcount\smallPLy\newdimen\smallPLyy\newdimen\smallPLX\newbox\smallPLdot \setbox\smallPLdot\hbox{{\tiny.}} \def\smallscl{.062} \def\smallPLot#1{\smallPLx`#1\advance\smallPLx-42\smallPLy\smallPLx\smallPLv\smallPLx\divide\smallPLy9\smallPLw\smallPLy\multiply \smallPLw9\advance\smallPLx-\smallPLw\advance\smallPLx-4\smallPLy-\smallPLy\advance\smallPLy4\smallPLX=\the\smallPLx pt \advance\smallPLyy\the\smallPLy pt\wd\smallPLdot=\smallscl\smallPLX\raise\smallscl\smallPLyy\copy\smallPLdot} \def\smalldraw#1{\ifx#1\end\let\next=\relax\else\smallPLot#1\let\next=\smalldraw\fi\next}

\newlength{\parrdp}\newlength{\parrht}
\newcommand{\parr}{\raisebox{-\parrdp}{\raisebox{\parrht}{\rotatebox{180}{$\&$}}}}
\newlength{\smallparrdp}\newlength{\smallparrht}
\newcommand{\smallparr}
{\mkern1mu\raisebox{-.1\smallparrdp}{\raisebox{\smallparrht}{\rotatebox{180}{\small$\&$}}}\mkern1mu}
\newlength{\footnoteparrdp}\newlength{\footnoteparrht}
\newcommand{\footnoteparr}
{\mkern1mu\raisebox{-\footnoteparrdp}{\raisebox{\footnoteparrht}{\rotatebox{180}{\footnotesize$\&$}}}\mkern1mu}
\newlength{\scriptparrdp}\newlength{\scriptparrht}
\newcommand{\scriptparr}
{\raisebox{-\scriptparrdp}{\raisebox{\scriptparrht}{\rotatebox{180}{\scriptsize$\&$}}}}

\newlength{\tinyparrdp}\newlength{\tinyparrht}
\newcommand{\A}{\mathbb{A}}
\newcommand{\C}{\mathbb{C}}

\newcommand{\op}{^{\mathsf{op}}}

\newcommand{\tensor}{\otimes}

\newcommand{\unit}{I}
\newcommand{\unitr}{\rho}
\newcommand{\unitl}{\lambda}
\newcommand{\unitrr}{\overline\rho}
\newcommand{\unitll}{\overline\lambda}

\newcommand{\id}{\mathsf{id}}
\newcommand{\tw}{\mathsf{tw}}
\newcommand{\assoc}{\alpha}

\newcommand{\sym}{\sigma}

\renewcommand{\perp}{^*}
\newcommand{\perpp}{{}^*{}^*}

\newcommand{\cl}[1]{\mathsf{L}#1}
\newcommand{\linka}{\cl{\A}}

\psset{radius=2pt,xunit=.5cm,yunit=1,arrowscale=1.7,arrowlength=.7,labelsep=3pt,nodesep=.5ex}%-.2mm

\newcommand{\lineanglesheight}[5]{\nccurve[angleA=#3,angleB=#4,ncurv=#5]{#1}{#2}}
\newcommand{\lineangles}[4]{\nccurve[angleA=#3,angleB=#4]{#1}{#2}}

\newcommand{\linevecpos}[3]{\ncline[ArrowInside=->,ArrowInsidePos=#3]{#1}{#2}}
\newcommand{\linevec}[2]{\linevecpos{#1}{#2}{.54}}
\newcommand{\vecanglesposheight}[6]{\nccurve[ncurv=#6,ArrowInside=->,ArrowInsidePos=#5,angleA=#3,angleB=#4]{#1}{#2}}
\newcommand{\vecanglespos}[5]{\vecanglesposheight{#1}{#2}{#3}{#4}{#5}{1}}
\newcommand{\vecanglesheight}[5]{\vecanglesposheight{#1}{#2}{#3}{#4}{.54}{#5}}
\newcommand{\vecangles}[4]{\vecanglespos{#1}{#2}{#3}{#4}{.54}}
\newcommand{\uvec}[2]{\vecangles{#1}{#2}{90}{-90}}
\newcommand{\dvec}[2]{\vecangles{#1}{#2}{-90}{90}}
\newcommand{\uvecpos}[3]{\vecanglespos{#1}{#2}{90}{-90}{#3}}
\newcommand{\dvecpos}[3]{\vecanglespos{#1}{#2}{-90}{90}{#3}}
\newcommand{\dloopvecleft}[2]{\nccurve[ncurv=1.3,angleA=-105,angleB=-75]{#1}{#2}
\lput{:0}{\psline{->}(.05,.005)(.1,.005)}}
\newcommand{\dloopvecanglesheight}[5]{\nccurve[ncurv=#5,angleA=#3,angleB=#4]{#1}{#2}
\lput{:0}{\psline{->}(.05,.005)(.1,.005)}}
\newcommand{\dloopvecheight}[3]{\dloopvecanglesheight{#1}{#2}{-75}{-105}{#3}}
\newcommand{\dloopvec}[2]{\dloopvecheight{#1}{#2}{1.3}}
\newcommand{\uloopvec}[2]{\nccurve[ncurv=1.3,angleA=75,angleB=105]{#1}{#2}
\lput{:0}{\psline{->}(.05,-.005)(.1,-.005)}}
\newcommand{\uloopvecanglesheight}[5]{\nccurve[ncurv=#5,angleA=#3,angleB=#4]{#1}{#2}
\lput{:0}{\psline{->}(.05,-.005)(.1,-.005)}}
\newcommand{\uloopvecheight}[3]{\uloopvecanglesheight{#1}{#2}{75}{105}{#3}}
\newcommand{\uloopvecleft}[2]{\nccurve[ncurv=1.3,angleA=105,angleB=75]{#1}{#2}
\lput{:0}{\psline{->}(.05,-.005)(.1,-.005)}}
\newcommand{\uloopvecleftheight}[3]{\nccurve[ncurv=#3,angleA=105,angleB=75]{#1}{#2}
\lput{:0}{\psline{->}(.05,-.005)(.1,-.005)}}

\newcommand{\one}[1]{\rnode{#1}{\unit}}
\newcommand{\nbot}[1]{\rnode{#1}{\bot}}
\newcommand{\nets}[1]{\mathsf{N}#1}
\newcommand{\netsa}{\nets{\A}}

\newcommand{\bgap}{\rule{9mm}{0mm}}% gap between GoI blobs
\newcommand{\negb}[1]{\cnode{2pt}{#1}}
\newcommand{\posb}[1]{\cnode*{2pt}{#1}}
\newcommand{\gnegb}[1]{\bgap\negb{#1}}
\newcommand{\gposb}[1]{\bgap\posb{#1}}
\newcommand{\posrestr}{^+}%{^\bullet}
\newcommand{\negrestr}{^-}%{^\circ}

\newcommand{\cutpair}[1]{{\psset{nodesep=1pt}\rnode{a}{#1}\;\,\rnode{b}{#1}\perp%
\nccurve[ncurv=1,angleA=-60,angleB=-120]{a}{b}}}
\newcommand{\dualcutpair}[1]{{\psset{nodesep=1pt}\rnode{a}{#1}\perp\;\,\rnode{b}{#1}\perpp%
\nccurve[ncurv=1,angleA=-60,angleB=-120]{a}{b}}}
\newcommand{\revcutpair}[1]{{\psset{nodesep=1pt}\rnode{a}{#1}\perp\,\rnode{b}{#1}%
\nccurve[ncurv=1,angleA=-60,angleB=-120]{a}{b}}}

\newcommand{\reduc}[1]{\underline{#1}}
\newcommand{\Matching}{\textsc{Matching}}
\newcommand{\Switching}{\textsc{Switching}}
\newcommand{\Labelling}{\textsc{Labelling}}
\newcommand{\Bijection}{\textsc{Bijection}}
\newcommand{\alabelpos}[5]{\ncline[linestyle=none]{#1}{#2}\naput[labelsep=#4,npos=#5]{#3}}
\newcommand{\blabelpos}[5]{\ncline[linestyle=none]{#1}{#2}\nbput[labelsep=#4,npos=#5]{#3}}
\newcommand{\clabelpos}[4]{\ncline[linestyle=none]{#1}{#2}\ncput[npos=#4]{#3}}
\newcommand{\alabel}[4]{\alabelpos{#1}{#2}{#3}{#4}{.5}}
\newcommand{\blabel}[4]{\blabelpos{#1}{#2}{#3}{#4}{.5}}
\newcommand{\clabel}[3]{\clabelpos{#1}{#2}{#3}{.5}}

\newcommand{\setp}{\mathsf{Setp}}

\renewcommand{\int}{\textbf{Int}}
\newcommand{\intsetp}{\int(\setp)}

\newcommand{\tighttensor}{\mkern-3mu\tensor\mkern-3mu}
\newcommand{\lf}{\mathcal{L}}
\newcommand{\leaves}[1]{|#1|}

\newcommand{\polynetstar}[2]{\text{Net}^*_{#1}(#2)}
\newcommand{\polynetstara}{\polynetstar{\mathcal{C}_\A}{E_\A}}

\newcommand{\netstara}{\mathsf{CircNet}{\A}}%{\netstar{\A}}
\newcommand{\mll}[1]{\widehat{#1}}
\newcommand{\LaxlinkingToCanonicalCircuit}[1]{#1^{\textsf{c}}}
\newcommand{\netToCircuitNet}[1]{#1^{\textsf{n}}}
\newcommand{\seqcomp}[2]{{#1};{#2}}
\newcommand{\iso}{\cong}

\title{\vspace{-5.5ex}\Large\bf 
Simple free star-autonomous categories and full coherence }%
\author{\\[-4ex]\normalsize\sc 
Dominic J.~D.~Hughes
\\[-.2ex]
\normalsize Stanford University\\
\small February 5, 2012}\date{}
\begin{document}
\maketitle\vspace{-5ex}
\begin{abstract}\noindent
This paper gives a simple presentation of the free star-autonomous
category over a category, based on Eilenberg-Kelly-MacLane graphs and
Trimble rewiring, yielding a full coherence theorem: the commutativity
of diagrams of canonical maps is decidable.
\vspace{-1.5ex}
\end{abstract}

\section{Introduction}

Eilenberg-Kelly-MacLane graphs \cite{EK66,KM71} elegantly describe
certain morphisms of closed categories.  This paper shows that little
more is needed to present the free star-autonomous category
\cite{Bar79} generated by a category, for a full coherence theorem:
the commutativity of diagrams of canonical maps is decidable.

Given a set $\A=\{a,b,\ldots\}$ of generators, we define the category
of \emph{$\A$-linkings}: objects are star-autonomous
\emph{shapes} 
(expressions)
over $\A$, 
such as $S\:=\:\big(a\tensor (b\perp\tensor
a\perpp)\big)\perpp{}\perp\tensor \unit\perp$ (with $\unit$ the unit),
and a
morphism $S\to T$
is a \emph{linking}, a function from negative leaves to positive
leaves, \eg
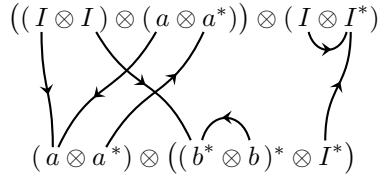
\begin{center}
\small
\label{firstexample}
\begin{psmatrix}[rowsep=3.3\baselineskip]
$\big((\,\one{u}\tensor\one{uu}\,)\tensor(\,\rnode{a}{a}\tensor\rnode{aa}{a}\perp)\big)
\tensor (\,\one{uuu}\tensor \one{uuuu}\perp)$ \\
$(\,\rnode{a1}{a}\tensor\rnode{aa1}{a}\,\perp) \tensor \big((\,\rnode{b1}{b}\perp\tensor\rnode{bb1}{b}\,)\perp
\tensor \one{u1}\perp\big)$
\dvec{u}{a1}
\vecanglespos{uu}{b1}{-60}{115}{.5}
\vecanglespos{a}{a1}{-120}{65}{.59}
\vecanglespos{aa1}{aa}{60}{-115}{.67}
\uloopvecleft{bb1}{b1}
\dloopvec{uuu}{uuuu}
\uvec{u1}{uuuu}
\end{psmatrix}
\end{center}
is a morphism from the upper shape to the lower shape.
Composition is simply path composition:\rput[b]{R}(2.5,-.5){\sf To appear in \textsl{Annals of Pure and Applied Algebra}, 2012.}
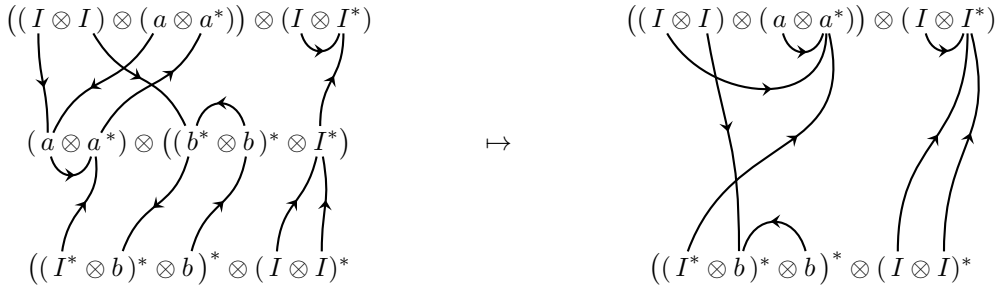
\begin{center}
\small
\begin{psmatrix}[rowsep=2.8\baselineskip]
$\big((\,\one{u}\tensor\one{uu}\,)\tensor(\,\rnode{a}{a}\tensor\rnode{aa}{a}\perp)\big)
\tensor (\one{uuu}\tensor \one{uuuu}\perp)$
&&
$\big((\,\one{u'}\tensor\one{uu'}\,)\tensor(\,\rnode{a'}{a}\tensor\rnode{aa'}{a}\perp)\big)
\tensor (\,\one{uuu'}\tensor \one{uuuu'}\perp)$
\\
$(\,\rnode{a1}{a}\tensor\rnode{aa1}{a}\,\perp) \tensor \big((\,\rnode{b1}{b}\perp\tensor\rnode{bb1}{b}\,)\perp
\tensor \one{u1}\perp\big)$
&
$\mapsto$ \\
$\big((\,\one{u2}\perp\tensor \rnode{b2}{b}\,)\perp\tensor \rnode{bb2}{b}\,\big)\perp
\tensor (\,\one{uu2} \tensor \one{uuu2})\perp$
&&
$\big((\,\one{u2'}\perp\tensor \rnode{b2'}{b}\,)\perp\tensor \rnode{bb2'}{b}\,\big)\perp
\tensor (\,\one{uu2'} \tensor \one{uuu2'})\perp$
\dvec{u}{a1}
\vecanglespos{uu}{b1}{-60}{115}{.5}
\vecanglespos{a}{a1}{-120}{65}{.59}
\vecanglespos{aa1}{aa}{60}{-115}{.67}
\uloopvecleft{bb1}{b1}
\dloopvec{uuu}{uuuu}
\uvec{u1}{uuuu}
\dloopvec{a1}{aa1}
\vecangles{u2}{aa1}{82}{-83}
\vecangles{b1}{b2}{-105}{70}
\vecangles{bb2}{bb1}{75}{-100}
\vecangles{uu2}{u1}{85}{-103}
\vecangles{uuu2}{u1}{90}{-80}
\vecanglesposheight{u'}{aa'}{-60}{-90}{.58}{1}
\vecanglespos{uu'}{b2'}{-80}{90}{.45}
\dloopvec{a'}{aa'}
\dloopvec{uuu'}{uuuu'}
\vecanglespos{u2'}{aa'}{70}{-77}{.58}
\uloopvecleft{bb2'}{b2'}
\vecangles{uu2'}{uuuu'}{90}{-90}
\vecangles{uuu2'}{uuuu'}{90}{-76}
\end{psmatrix}
\end{center}
A leaf function qualifies as a linking only if it satisfies the
standard criterion for multiplicative proof
nets\footnote{\label{assume}This paper will not assume any familiarity
with proof nets or linear logic.} \cite{DR89,Gir96}, so simple as to
be checkable in linear time \cite{Gue99,Hug12}.  Employing Trimble
rewiring
\cite{Tri94,BCST96}, we define two linkings as \emph{similar} if they differ
by an edge from an $\unit$, \eg
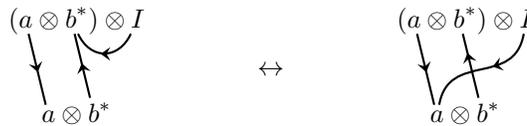
\begin{center}
\small
\vspace{-2.7ex}
\begin{psmatrix}[rowsep=.5\baselineskip]
$(\rnode{a}{a}\tensor\rnode{b}{b}\perp)\tensor\one{u}$ &&
$(\rnode{a'}{a}\tensor\rnode{b'}{b}\perp)\tensor\one{u'}$ \\ &
$\leftrightarrow$ \\ $\rnode{a1}{a}\tensor\rnode{b1}{b}\perp$ &&
$\rnode{a1'}{a}\tensor\rnode{b1'}{b}\perp$
\linevec{a}{a1}
\linevec{b1}{b}
\nccurve[ncurv=1.1,angleA=-115,angleB=-65]{u}{b}
\lput{:0}{\psline{->}(.05,.005)(.1,.005)}
\linevec{a'}{a1'}
\linevecpos{b1'}{b'}{.7}
\vecanglespos{u'}{a1'}{-100}{80}{.4}
\end{psmatrix}
\end{center}
and define an \emph{$\A$-net} as an $\A$-linking modulo similarity.
The category of $\A$-nets is the free star-autonomous category
generated by $\A$.  To emphasise the simplicity:\vspace{-1ex}
\begin{itemize}
\myitem{1} a morphism is a leaf function satisfying a standard criterion (checkable in linear time);
\vspace{-1ex}\myitem{2} composition is standard path composition;
\vspace{-1ex}\myitem{3} modulo a standard equivalence (Trimble rewiring).
\end{itemize}

The key novelty is (2), the fact that composition is simply path
composition.  This preserves an elegant feature of
Eilenberg-Kelly-MacLane graphs.  In contrast, the composition in
previous presentations of free star-autonomous categories
\cite{BCST96,KO99,LS04} is more complex.  (We return to this related
work later in the Introduction.)

Abstractly, the underlying path composition can be understood as a
forgetful functor from the category of $\A$-linkings (sketched above)
to $\intsetp$, the compact closed category obtained by applying the
$\int$ or geometry-of-interaction construction
\cite{Gir89,JSV96,Abr96} to the traced monoidal category $\setp$ of
sets and partial functions (with coproduct as tensor).  This ties in
nicely with Eilenberg-Kelly-MacLane graphs and Kelly-Laplaza graphs
\cite{KL80} for compact closed categories, since each has a forgetful
functor to $\intsetp$.

\parag{Arbitrary base category}
When $\A$ is not discrete, we simply label each edge between generators with
a morphism of $\A$.  For example, if $x:a\to b$ and $y:c\to c$ in
$\A$, then
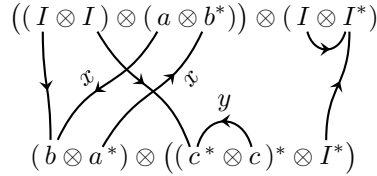
\begin{center}\small
\begin{psmatrix}[rowsep=3.3\baselineskip]
$\big((\,\one{u}\tensor\one{uu}\,)\tensor(\,\rnode{a}{a}\tensor\rnode{b}{b}\perp)\big)
\tensor (\,\one{uuu}\tensor \one{uuuu}\perp)$ \\
$(\,\rnode{b1}{b}\tensor\rnode{a1}{a}\,\perp) \tensor \big((\,\rnode{c1}{c}\,\perp\tensor\rnode{cc1}{c}\,)\perp
\tensor \one{u1}\perp\big)$
\dvec{u}{b1}
\vecanglespos{uu}{c1}{-60}{110}{.47}
\vecanglespos{a}{b1}{-120}{65}{.59}\nbput[npos=.59,nrot=:D]{$x$}
\vecanglespos{a1}{b}{60}{-115}{.67}\nbput[npos=.76,nrot=:U]{$x$}
\uloopvecleft{cc1}{c1}\nbput[npos=.5,nrot=:D]{$y$}
\dloopvec{uuu}{uuuu}
\uvec{u1}{uuuu}
\end{psmatrix}
\end{center}
is an $\A$-linking.  Composition collects labels along a path, and
composes\footnote{Througout this paper we employ sequential notation
  $\seqcomp f g$ for the composite of $f:S\to T$ and $g:T\to U$ rather
  than functional notation $gf$, since it is more natural in a
  diagrammatic, path-following setting.} them in $\A$, \eg
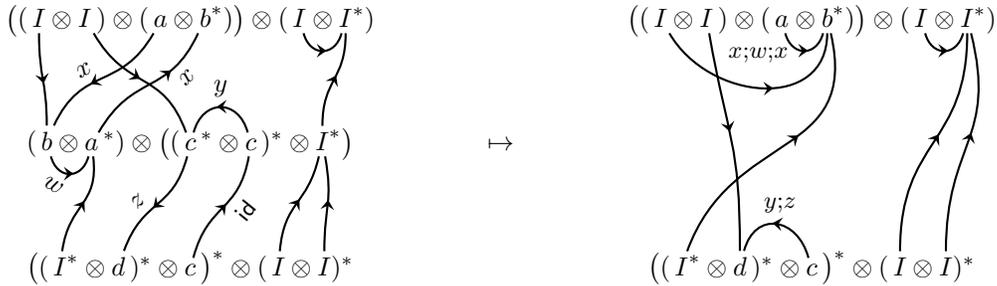
\begin{center}\small\label{label-comp-eg}
\begin{psmatrix}[rowsep=2.8\baselineskip]
$\big((\,\one{u}\tensor\one{uu}\,)\tensor(\,\rnode{a}{a}\tensor\rnode{b}{b}\perp)\big)
\tensor (\,\one{uuu}\tensor \one{uuuu}\perp)$
&&
$\big((\,\one{u'}\tensor\one{uu'}\,)\tensor(\,\rnode{a'}{a}\tensor\rnode{b'}{b}\perp)\big)
\tensor (\,\one{uuu'}\tensor \one{uuuu'}\perp)$
\\
$(\,\rnode{b1}{b}\tensor\rnode{a1}{a}\,\perp) \tensor \big((\,\rnode{c1}{c}\,\perp\tensor\rnode{cc1}{c}\,)\perp
\tensor \one{u1}\perp\big)$
&
$\mapsto$ \\
$\big((\,\one{u2}\perp\tensor \rnode{d2}{d}\,)\perp\tensor \rnode{c2}{c}\,\big)\perp
\tensor (\,\one{uu2} \tensor \one{uuu2})\perp$
&&
$\big((\,\one{u2'}\perp\tensor \rnode{d2'}{d}\,)\perp\tensor \rnode{c2'}{c}\,\big)\perp
\tensor (\,\one{uu2'} \tensor \one{uuu2'})\perp$
\dvec{u}{b1}
\vecanglespos{uu}{c1}{-60}{110}{.47}
\vecanglespos{a}{b1}{-120}{65}{.59}\nbput[npos=.6,nrot=:D]{$x$}
\vecanglespos{a1}{b}{60}{-115}{.67}\nbput[npos=.76,nrot=:U]{$x$}
\uloopvecleft{cc1}{c1}\nbput[npos=.5,nrot=:D]{$y$}
\dloopvec{uuu}{uuuu}
\uvec{u1}{uuuu}
\dloopvec{b1}{a1}\nbput[npos=.32,nrot=:U]{$w$}
\vecangles{u2}{a1}{82}{-83}
\vecangles{c1}{d2}{-105}{70}\nbput[npos=.62,nrot=:U]{$z$}
\vecangles{c2}{cc1}{75}{-100}\nbput[npos=.63,nrot=:U]{\footnotesize$\id$}
\vecangles{uu2}{u1}{85}{-103}
\vecangles{uuu2}{u1}{90}{-80}
\vecanglesposheight{u'}{b'}{-60}{-90}{.58}{1}
\vecanglespos{uu'}{d2'}{-80}{90}{.45}
\dloopvec{a'}{b'}\nbput[npos=.2,labelsep=1pt]{\footnotesize$x;\!w;\!x\mkern-1mu$}
\dloopvec{uuu'}{uuuu'}
\vecanglespos{u2'}{b'}{70}{-77}{.58}
\uloopvecleft{c2'}{d2'}\nbput[npos=.5]{\footnotesize$y;\!z$}
\vecangles{uu2'}{uuuu'}{90}{-90}
\vecangles{uuu2'}{uuuu'}{90}{-76}
\end{psmatrix}
\end{center}
The category of $\A$-nets ($\A$-linkings modulo Trimble rewiring, as
before) is the free star-autonomous category generated by the category $\A$.

\paragraph*{Full coherence.}
Equivalence modulo rewiring is decidable, by finiteness.  Thus we have
a full coherence theorem: we can decide the commutativity of diagrams
of canonical maps in star-autonomous categories.  Here are two short
illustrative examples.

\emph{Example 1: identity $\neq$ twist on $\bot=\unit\perp$.}\label{intro-rewiring-examples}\label{intro-rewiring-example1}
Let $\tw_{A\tensor B}:A\tensor B\to B\tensor A$ be the canonical twist
(symmetry) isomorphism.  
The identity and twist $\id_{\bot\tensor\bot},\tw_{\bot\tensor\bot}:\bot\tensor
\bot\to \bot\tensor \bot$ determine respective linkings $i$ and $t$:
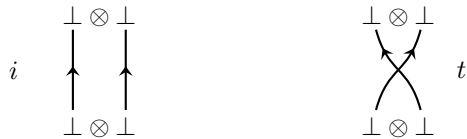
\begin{center}\small\psset{rowsep=2.5\baselineskip}$\begin{psmatrix}
\nbot{b}\tensor\nbot{bb} \\ 
\nbot{b1}\tensor\nbot{bb1}
\uvec{b1}{b}\uvec{bb1}{bb}\blabel{b}{b1}{i}{.7cm} 
\end{psmatrix}
\rule{3cm}{0cm}
\begin{psmatrix} 
\nbot{b}\tensor\nbot{bb} \\
\nbot{b1}\tensor\nbot{bb1}
\vecanglespos{b1}{bb}{75}{-105}{.78}\vecanglespos{bb1}{b}{105}{-75}{.78}
\alabel{bb}{bb1}{t}{.4cm} 
\end{psmatrix}$\end{center}
They differ on two inputs,
and there are no other%
\footnote{\label{no-others}The
other two functions from negative leaves to positive leaves fail the
proof net criterion.}
linkings $\bot\tensor\bot\to\bot\tensor\bot$, so $i$ and $t$ cannot be rewired into each other.
Thus, in general,
$\id_{\bot\tensor\bot}\neq\tw_{\bot\tensor\bot}:\bot\tensor\bot\to\bot\tensor\bot$.

\emph{Example 2: Triple-dual problem.}\label{intro-rewiring-example2}
We show that the following diagram commutes (a triple-dual problem
\cite[\S1, diagram\ (1.4)]{KM71})
\begin{equation}\label{triangle}
\psset{nodesep=1ex}
\begin{psmatrix}[ref=c]
\big((a\multimap\bot)\multimap\bot\big)\multimap\bot &
a\multimap\bot\rule{0pt}{3ex} \\
& 
\big((a\multimap\bot)\multimap\bot\big)\multimap\bot
\ncline{->}{1,1}{1,2}\naput{k_a\multimap \id}
\ncline{->}{1,2}{2,2}\naput{k_{(a\multimap\mkern2mu\bot)}}
\ncline{->}{1,1}{2,2}\nbput{\id}
\end{psmatrix}
\end{equation}
where 
$A\multimap B=(A\tensor B\perp)\perp$ and
$k_A:A\;\to\;(A\mkern-2mu\multimap\mkern-2mu\bot)\mkern-2mu\multimap\mkern-2mu\bot$ is the canonical map of
its type.  Each path in the triangle determines a
corresponding linking
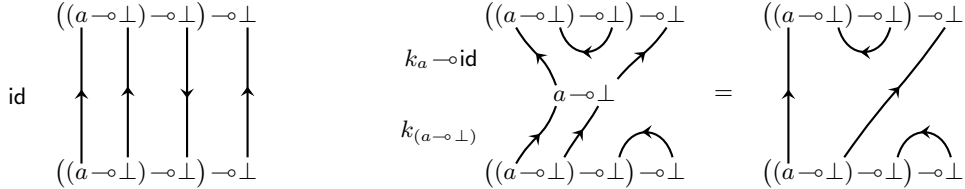
\begin{center}\footnotesize\begin{math}
\newcommand{\tightmap}{\mkern-3mu\multimap\mkern-4mu}
\psset{rowsep=1.7\baselineskip,nodesep=2pt}
\rule{.3cm}{0pt}\begin{psmatrix}
\big((\rnode{a}{a}\tightmap\nbot{b})\tightmap\nbot{bb}\big)\tightmap\nbot{bbb}\\\\
\big((\rnode{a1}{a}\tightmap\nbot{b1})\tightmap\nbot{bb1}\big)\tightmap\nbot{bbb1}
\uvec{a1}{a}
\uvec{b1}{b}
\dvec{bb}{bb1}
\uvec{bbb1}{bbb}
\blabel{a}{a1}{\id}{.7cm} 
\end{psmatrix}
\rule{3cm}{0cm}
\rnode{A}{\begin{psmatrix} 
\big((\rnode{a}{a}\tightmap\nbot{b})\tightmap\nbot{bb}\big)\tightmap\nbot{bbb}\\
\rnode{am}{a}\tightmap\nbot{bm} \\
\big((\rnode{a1}{a}\tightmap\nbot{b1})\tightmap\nbot{bb1}\big)\tightmap\nbot{bbb1}
\blabelpos{a}{a1}{k_a\tightmap \id}{.4cm}{.23} 
\blabelpos{a}{a1}{k_{(a\multimap\mkern2mu\bot)}}{.4cm}{.77} 
\vecangles{am}{a}{110}{-60}
\vecangles{a1}{am}{60}{-110}
\dloopvecleft{bb}{b}
\vecangles{bm}{bbb}{50}{-125}
\vecangles{b1}{bm}{60}{-120}
\uloopvecleft{bbb1}{bb1}
\end{psmatrix}}
\rule{1cm}{0cm}
\rnode{B}{\begin{psmatrix} 
\big((\rnode{a}{a}\tightmap\nbot{b})\tightmap\nbot{bb}\big)\tightmap\nbot{bbb}\\
\\
\big((\rnode{a1}{a}\tightmap\nbot{b1})\tightmap\nbot{bb1}\big)\tightmap\nbot{bbb1}
\uvec{a1}{a}
\dloopvecleft{bb}{b}
\vecangles{b1}{bbb}{55}{-125}
\uloopvecleft{bbb1}{bb1}
\end{psmatrix}}
\clabel{A}{B}{=}
\end{math}\end{center}
equivalent via three rewirings:
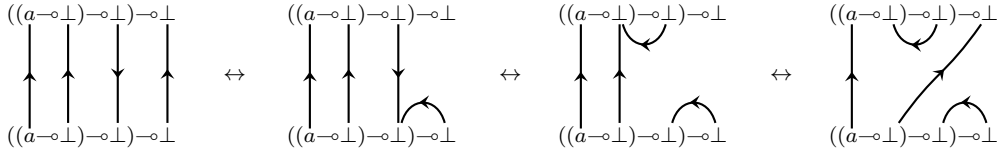
\begin{center}\scriptsize
\psset{nodesep=.5ex,rowsep=1.3cm}
\newcommand{\tightmap}{\mkern-5mu\multimap\mkern-5mu}
\newcommand{\core}[2]{%
\rnode{#1}{\psframebox[linestyle=none]{%
\psset{nodesep=.5ex,rowsep=1.3cm}\begin{psmatrix}
$((\rnode{a'}{a}\tightmap\nbot{b'})\tightmap\nbot{bb'})\tightmap\nbot{bbb'}\;$
\\
$((\rnode{a1'}{a}\tightmap\nbot{b1'})\tightmap\nbot{bb1'})\tightmap\nbot{bbb1'}\;$
\uvec{a1'}{a'}
#2
\end{psmatrix}}}}
\begin{psmatrix}
  \core{D}{%
    \dvec{bb'}{bb1'}
    \uvec{b1'}{b'}
    \uvec{bbb1'}{bbb'}
  }
  & \hspace{-3ex} &
  \core{C}{%
    \dvec{bb'}{bb1'}
    \uvec{b1'}{b'}
    \uloopvecleft{bbb1'}{bb1'}
  }
  & \hspace{-4ex} &
  \core{B}{%
    \dloopvecleft{bb'}{b'}
    \uvec{b1'}{b'}
    \uloopvecleft{bbb1'}{bb1'}
  }
  & \hspace{-4ex} &
  \core{A}{%
    \dloopvecleft{bb'}{b'}
    \vecangles{b1'}{bbb'}{55}{-125}
    \uloopvecleft{bbb1'}{bb1'}
  }
\end{psmatrix}
  \psset{linestyle=none}
  \ncline{A}{B}\ncput{$\leftrightarrow$}
  \ncline{B}{C}\ncput{$\leftrightarrow$}
  \ncline{C}{D}\ncput{$\leftrightarrow$}
\end{center}
Thus we conclude that triangle (\ref{triangle}) commutes in every
star-autonomous category.\footnote{\label{path}Compare with
  similar arguments in \cite[\S4.2]{BCST96}, \cite[\S2]{KO99}
  and \cite[\S10]{MO03}.  A key advantage here is that, because of
  the simple path composition, the composite
  $\:\seqcomp{k_a\mkern-3mu\multimap\mkern-3mu\id\:}{\:k_{(a\multimap\bot)}}\:$
  is immediate on inspection.}

\parag{Related work}\label{related-work}

This paper follows an approach which can be traced back to Todd
Trimble's Ph.D.{} thesis \cite{Tri94}\footnote{Copies of Trimble's
thesis \cite{Tri94} are not particularly easy to come by.  See
\cite[\S1,\,\S3.2]{BCST96} and \cite[\S3]{CS97} for overviews of some
of the content.}.  We call this the \emph{rewiring approach}:
\begin{itemize}
\myitem{a} represent a morphism by a structure involving attachments (`wiring') of
negative units\footnote{For history and development of the attachment
of negative units in linear logic, see
\cite{Dan90,Reg92,GSS92,Gir93,Gir96}.};
\myitem{b} quotient by rewiring: identify correct structures which differ by just 
one such attachment.
\end{itemize}
This is the fourth paper to use the rewiring approach to construct free
star-autonomous categories.  The chronological sequence is detailed 
below, and is summarized in Table~\ref{approach-table}.
\begin{table}\footnotesize\newcommand{\cell}[2]{\parbox{#2in}{\center #1}}
\newcommand{\cellone}[1]{\cell{\bf#1}{.45}}
\newcommand{\celltwo}[1]{\cell{#1}{1.33}}
\newcommand{\cellthree}[1]{\cell{#1}{.6}}
\newcommand{\cellfour}[1]{\cell{#1}{.9}}
\newcommand{\cellfive}[1]{\cell{#1}{.8}}
\newcommand{\cellsix}[1]{\cell{#1}{.8}}
\setlength{\doublerulesep}{\arrayrulewidth}
\begin{tabular}{|c|||c|c|c|c|c|}
\hline
& \celltwo{\bf\textsl{Structure}} &
\cellthree{\bf\textsl{\hspace*{-2.12315pt}Attachment of negative
unit}} & 
\cellfour{\bf\textsl{Correctness / allowability}} & 
\cellfive{\bf\textsl{Rewiring of attachments, between correct structures}}
&\cellsix{\bf\textsl{Problem with composition / normalisation}}
\\&&&&&\\\hline\hline\hline
\cellone{\makebox{\hspace*{-1ex}\cite{BCST96}\hspace*{-1ex}}} & \celltwo{Circuit diagram\\[1ex] (proof net in tensor calculus style)} &
\cellthree{Thinning link (dotted edge)} &
\cellfour{Standard proof\\ net criterion\\[1ex]\scriptsize (\mbox{sequentialisability} / contractibility for planar case)} & 
\cellfive{Surgery rules / re-target thinning link (by \cite[Prop.\,3.3]{BCST96})} 
& \cellsix{Attachments can block cut redexes}
\\&&&&&\\\hline
\cellone{\cite{KO99}} & 
\celltwo{$\lambda\mu$-style term, \eg{} \footnotesize$(\mathbf{\lambda}\beta^{A\multimap\bot}.\langle\beta\rangle z^A)\;\;\;\;\;\;$ $\;\;\;\;\;\;\;\{\mu
\alpha^{A\multimap\bot}.\langle\gamma\rangle\alpha/z\}$\!\!\!\!\!\!\!\!} &
\cellthree{Unit let term constructor $\langle x/\ast\rangle(-)$} &
\cellfour{Typability\\ (\ie, sequential\-isability)} &
\cellfive{$\pi$-congruence rule, $\:\Gamma\vdash C[\langle x/\ast\rangle t]\sim\langle x/\ast\rangle C[t] : A$}
& \cellsix{normalisation confluent only modulo rewiring / congruence}
\\&&&&&\\\hline
\cellone{\cite{LS04}} & 
\celltwo{MLL formula sharing sequent leaves, with permutation, \eg\\
\footnotesize\mbox{\hspace{-1ex}$\bot_1\tighttensor((\assoc_2\tighttensor a\perp_4)\footnoteparr(\bot_3\tighttensor \unit_5))$\hspace{-.8ex}} \mbox{\normalsize$\triangleright$} $\;\;\;\bot\footnoteparr a,(\bot\tighttensor a\perp)\footnoteparr \unit$} &
\cellthree{Formula constructor $(-)\tensor\bot$} &
\cellfour{Standard proof\\ net criterion} &
\cellfive{Formula rewrite
$Q\footnoteparr (R\tensor\bot)\leftrightarrow (Q\footnoteparr
R)\tensor\bot$}
& \cellsix{Attachments can block cut redexes}
\\&&&&&\\\hline
\cellone{This paper} & \celltwo{Leaf function} & \cellthree{Edge\\ from $\unit$} &
\cellfour{Standard proof\\ net criterion} & 
\cellfive{Re-target edge from $\unit$} &
\cellsix{}
\\&&&&&\\\hline
\end{tabular}
\caption{\label{approach-table}This paper is the fourth to construct
  free star-autonomous categories using the \emph{`rewiring
    approach}', which can be traced back to Trimble's thesis
  \cite{Tri94}: (a) represent a morphism by a structure involving
  attachments (`wiring') of negative units, (b) quotient by rewiring,
  that is, identify correct structures which differ by just one such
  attachment.  In each case a morphism of the free category is a
  finite equivalence class, hence equality of morphisms is
  decidable.}\vspace{1ex}
\end{table}

In \cite{BCST96} structures are circuit diagrams (in tensor calculus
style \cite{JS91}), attachments are (dotted) edges from negative
units, called thinning links, correctness is the standard
multiplicative proof net
criterion\footnote{Sequentialisability/contractibility \cite{Dan90} is
  used to deal with the planar case; see Section~2.7 of
  \cite{BCST96}.}, and rewiring between correct structures is
expressed in rules of surgery on circuits, which (by the empire
rewiring Proposition~3.3 of \cite{BCST96}) permit an arbitrary
re-targeting of a thinning link between correct
circuits\footnote{Decomposing rewiring into shorter steps aided the
  freeness proofs in \cite{BCST96}.}.  Equivalence classes yield the free
linearly distributive category and free star-autonomous category
(linearly distributive category with negation) generated by a
polygraph (\eg, by a category), for full coherence.

In \cite{KO99} structures are $\lambda\mu$-style terms \cite{Par92}
with explicit substitution $\{-/-\}$, for example
$(\mathbf{\lambda}\beta^{A\multimap\bot}.\langle\beta\rangle z^A)\{\mu
\alpha^{A\multimap\bot}.\langle\gamma\rangle\alpha/z\}$\:, attachments are unit let
constructs $\langle x/\ast\rangle(-)$, correctness is inductive
(typability, \ie, sequentialisability) and rewiring is by an instance
of the $\pi$-congruence rule, $\:\Gamma\vdash C[\langle x/\ast\rangle
t]\sim\langle x/\ast\rangle C[t] : A$\:.
Equivalence (congruence) classes yield an internal language
for autonomous and star-autonomous categories, for full coherence.

In \cite{LS04} the structure is a syntactic\footnote{Axiom links
$a\tensor a^\bot$ and attachments of negative units $(-)\tensor\bot$
are syntactic, enveloped in a formula sharing the leaves of the
sequent. In conventional proof nets \cite{Gir96}, axiom links and unit
attachments ($\bot$-jumps) are edges.}
proof net, a formula of multiplicative linear logic equipped with a
leaf permutation, \eg\
\mbox{\:$\bot_1\tighttensor((\assoc_2\tighttensor
a\perp_4)\parr(\bot_3\tighttensor \unit_5))$}\:
\mbox{\large$\triangleright$}\:$\bot\parr
a,(\bot\tighttensor a\perp)\parr \unit\;$ is a structure representing a
morphism $\unit\tensor a\perp\to (\bot\tensor a\perp)\parr I$, the
formula constructor $(-)\tensor\bot$ attaches negative units,
correctness is again the standard multiplicative proof net criterion,
and rewiring is by the invertible linear distributivity rewrite
$Q\parr (R\tensor\bot)\leftrightarrow (Q\parr R)\tensor\bot$\:.
Equivalence classes yield the free star-autonomous category with strict
double involution\footnote{The canonical map $A\perpp\to A$ is the
identity.  Up to equivalence this is a free star-autonomous category,
in a strict sense \cite{CHS05}.} generated by a set, for full
coherence.

At first sight, it may seem repetitive and uninteresting to employ the
rewiring approach for star-autonomous categories a fourth time.
However, the simplicity of the end product relative to the previous
approaches seems to justify the repetition.
As we remarked earlier, we preserve an elegant feature of
Eilenberg-Kelly-MacLane graphs:

\begin{itemize}
\vspace{1ex}
\item \sl Composition is simply path composition.
\vspace{1ex}
\end{itemize}
Composition is more complex in the three previous approaches to free
star-autonomous categories \cite{BCST96,KO99,LS04}. 
In each case, given normal forms $s$ and $t$ representing equivalence
classes, one first forms a `concatenation' $s;t$ (in \cite{BCST96},
pasting the circuits at a cut wire, in \cite{KO99} forming an explicit
substitution, and in \cite{LS04} forming a proof net with cuts), then
normalises $s;t$ in a rewrite system.
In \cite{BCST96} and \cite{LS04} normalisation is defined only modulo
equivalence, since unit attachments (thinning links in
\cite{BCST96} and $(-)\tensor\bot$ in \cite{LS04}) can block cut
redexes\footnote{A problem also discussed by Girard, in the context of
proof nets \cite[\S{}A.2]{Gir96}.}, and in
\cite{KO99} confluence is only modulo equivalence (congruence).
In contrast, path composition, as in this paper and
Eilenberg-Kelly-MacLane graphs, is simple and direct.

This paper is the sequel to \cite{Hug12}\footnote{It was tempting to merge the two papers.  However,
some proof theorists and linear logicians may not be interested in
star-autonomous categories and the emphasis on coherence over
correctness criteria (combinatorial characterisations of
\emph{allowability}, in the parlance of \cite{KM71}), and conversely,
some categorists may not be interested in linear logic and its sequent
calculus, and the emphasis on correctness criteria over coherence. The
present paper targets categorists, and assumes no familiarity with
linear logic.} on multiplicative proof nets with units, which relate
closely to $\A$-linkings, for $\A$ discrete.
For comprehensive background and history on free star-autonomous
categories and coherence, see the introductions of
\cite{BCST96,KO99,LS04}.

\parag{Potential future work}  

Perhaps the most direct redeployment of Trimble rewiring is
\cite{MO03}, since it is for SMCCs (symmetric monoidal closed
categories), the original case treated in \cite{Tri94}.\footnote{The
  definition of $\sim$ in Section~9 of \cite{MO03} is precisely
  Trimble rewiring: two (correct) structures (strategies / linkages /
  generating functions) are identified if they differ by the
  attachment of just one negative unit (there called a joker move).
  Syntactically (see the discussion before Proposition~45 of
  \cite{MO03}), Trimble rewiring here corresponds to that in
  \cite{KO99}: $\pi$-congruence with the unit let construct,
  $\:\Gamma\vdash C[\langle x/\ast\rangle t]\sim\langle x/\ast\rangle
  C[t] : A$\:.}  Hybridising the present paper on star-autonomous
categories with the extension of Lamarche's essential nets
\cite{Lam94} in \cite{MO03} might yield a simple presentation of free
SMCCs: objects are shapes generated by $\tensor$, $\multimap$ and
$\unit$ (\eg\ $(a\tensor \unit)\multimap\bot$), a morphism is a leaf
function satisfying Lamarche's criterion, modulo Trimble rewiring, and
composition is simply path composition.

Path composition would constitute a \emph{direct} composition of
generating functions, bypassing the complexities of strategies
(O-orientation, shortsightedness, non-determinism, conditional
exhaustion, \etc) for a more economical description of the free SMCC.
Furthermore, the approach would extend immediately to the free SMCC
generated by an arbitrary category $\A$, not just a set $\A$ (as in
\cite{MO03}), by labelling edges with morphisms of $\A$ (as in the
present paper for star-autonomous categories).  In summary, this path
composition approach would abstract away from the intricacies of the
strategies, extracting the essence: a geometry of interaction of
generating functions in $\intsetp$.

\parag{Acknowledgements} Many thanks to Robin Cockett and Robert Seely
for helping me understand the construction of free star-autonomous
categories in \cite{BCST96}, an important precursor to this
paper.  I am extremely grateful to Robin Houston for insightful
feedback, in particular for improvements to the definition of
\emph{split star-autonomous category}.  Thanks to Peter Selinger for
corrections.  I am indebted to an anonymous referee for excellent
suggestions.

\section{Split star-autonomous categories}\label{split-star}

The category of $\A$-linkings (sketched above, and defined in the next
section) is almost, but not quite, star-autonomous.  It becomes
star-autonomous upon quotienting by Trimble rewiring.  Below we
axiomatise its raw structure, prior to quotienting, as a
\emph{split star-autonomous category}, defined by relaxing the
unit isomorphisms $A\to \unit\tensor A$ and $A\to A\tensor\unit$ of a
star-autonomous category to be split monomorphisms (sections).

A \defn{star-autonomous category} \cite{Bar79} is a category $\C$
equipped with the following structure:
\begin{enumerate}
\myitem{1} \emph{Tensor}. A functor $-\tensor-:\C\times \C\to\C$.
\myitem{2} \emph{Associativity}. 
A natural isomorphism $\assoc_{A,B,C}:(A\tensor B)\tensor C\to
A\tensor(B\tensor C)$, natural in objects $A,B,C\in\C$, such that the
following pentagon commutes:
\begin{displaymath}
\begin{psmatrix}[ref=c,nodesep=1ex,colsep=-4ex]
\big((A\tensor B)\tensor C\big)\tensor D \rule{0pt}{3ex}
&& (A\tensor B)\tensor (C\tensor D) && A\tensor \big(B\tensor (C\tensor
D)\big)
\\
& \big(A\tensor(B\tensor C)\big)\tensor D
&&
A\tensor\big((B\tensor C)\tensor D\big)
\ncline{->}{1,1}{2,2}\nbput[npos=.34]{\assoc\tensor \id}
\ncline{->}{1,1}{1,3}\naput{\assoc}
\ncline{->}{1,3}{1,5}\naput{\assoc}
\ncline{->}{2,2}{2,4}\nbput{\assoc}
\ncline{->}{2,4}{1,5}\nbput[npos=.66]{\id\tensor\assoc}
\end{psmatrix}
\end{displaymath}
\myitem{3} \emph{Unit}.  An object $\unit\in \C$.
\myitem{4} \emph{Unit isomorphisms}.  
Natural isomorphisms\footnote{Conventionally the unit isomorphisms are typed
  $A\tensor \unit\to A$ and $\unit\tensor A\to A$\; \cite{Mac71}.  We
  reverse them to ease the definition of split star-autonomous
  category below.}
\begin{displaymath}
  \begin{array}{rcccl}
    \unitl_A&:&A&\to&\unit\tensor A\\
    \unitr_A&:&A&\to&A\tensor\unit\\
  \end{array}
\end{displaymath}
%Natural isomorphisms $\unitl_A:A\to\unit\tensor A$ and $\unitr_A:A\to A\tensor\unit$, 
natural in the object $A\in\C$, such that the following triangle
commutes:
\begin{displaymath}
\begin{psmatrix}[ref=c,nodesep=1ex,colsep=1ex]
{} & A\tensor B \\
(A\tensor\unit)\tensor B && A\tensor (\unit\tensor B) 
\ncline{->}{2,1}{2,3}\nbput{\assoc}
\ncline{->}{1,2}{2,1}\nbput[npos=.57]{\unitr\tensor\id}
\ncline{->}{1,2}{2,3}\naput[npos=.57]{\id\tensor\unitl}
\end{psmatrix}
\end{displaymath}
%% \myitem{4} \emph{Unit isomorphisms}.  Natural isomorphisms
%% $\unitl_A:\unit\tensor A\to A$ and $\unitr_A:A\tensor\unit\to A$,
%% natural the object $A\in\C$, such that the following triangle
%% commutes for all $B$:\vspace{0.4ex}
%% \begin{displaymath}
%% \begin{psmatrix}[ref=c,nodesep=1ex,colsep=1ex]
%% (A\tensor\unit)\tensor B && A\tensor (\unit\tensor B) \\
%% {} & A\tensor B
%% \ncline{->}{1,1}{1,3}\naput{\assoc}
%% \ncline{->}{1,1}{2,2}\nbput[npos=.57]{\unitr\tensor\id}
%% \ncline{->}{1,3}{2,2}\naput[npos=.57]{\id\tensor\unitl}
%% \end{psmatrix}
%% \end{displaymath}
\myitem{5} \emph{Symmetry}.
A natural isomorphism $\sym_{A,B}:A\tensor B\to B\tensor A$, natural in
objects $A,B\in\C$, such that 
the following diagrams commute:\vspace{.8ex}
\begin{displaymath}\begin{psmatrix}[ref=c,nodesep=1ex,colsep=1ex]
A\tensor B  && A\tensor B \\
{} & B\tensor A  
\ncline{->}{1,1}{1,3}\naput{\id}
\ncline{->}{1,1}{2,2}\nbput{\sym}
\ncline{->}{2,2}{1,3}\nbput{\sym}
\end{psmatrix}
\hspace{5ex}
\begin{psmatrix}[ref=c,nodesep=1ex,colsep=8.5ex]
(A\tensor B)\tensor C 
&
A\tensor (B\tensor C)
& 
(B\tensor C)\tensor A
\\
(B\tensor A)\tensor C
&
B\tensor (A\tensor C)
&
B\tensor(C\tensor A)
\ncline{->}{1,1}{1,2}\naput{\assoc}
\ncline{->}{1,2}{1,3}\naput{\sym}
\ncline{->}{1,3}{2,3}\naput{\assoc}
\ncline{->}{1,1}{2,1}\nbput{\sym\tensor\id}
\ncline{->}{2,1}{2,2}\nbput{\assoc}
\ncline{->}{2,2}{2,3}\nbput{\id\tensor\sym}
\end{psmatrix}\end{displaymath}
\myitem{6} \emph{Involution}. A full and faithful functor $(-)\perp:\C\op\to\C$.
\myitem{7} \emph{Closure}. A natural isomorphism $\C(A\tensor B,C\perp)\to\C(A,(B\tensor C)\perp)$, natural 
in objects $A,B,C\in\C$.
\end{enumerate}
Axioms (1)--(4) define a monoidal category and (1)--(5) define a
symmetric monoidal category \cite{Mac71}.\footnote{In \cite{Mac71} Mac
  Lane demands $\unitr_I=\unitl_I$ for a monoidal category and
  $\unitl_A;\sym_{\unit,A}=\unitr_A$ for a symmetric monoidal
  category, which are superfluous \cite{Kel64,JS93}.}  The above is
not the original definition of star-autonomous category, but (modulo
our slightly different presentation of symmetric monoidal category) is
equivalent \cite{Bar79}.

Define a \defn{split star-autonomous category} by relaxing (4): demand
only that the natural transformations $\unitl_A:A\to\unit\tensor A$
and $\unitr_A:A\to A\tensor\unit$ be split monomorphisms (sections),
rather than than isomorphisms.  Thus we require for each $A$ the
existence\footnote{An anonymous referee noted that if we include a
  specific choice of retractions for each $A$ in the definition of
  split star-autonomous category, the definition becomes monadic over
  the category of categories.}  of retractions $\unitll_A:\unit\tensor
A\to A$ and $\unitrr_A:A\tensor \unit\to A$ such that
$\unitl_A;\unitll_A=\id_A$ and
$\unitr_A;\unitrr_A=\id_A$.\vspace{.8ex}%
\begin{displaymath}
\begin{psmatrix}[ref=c,nodesep=1ex,colsep=4ex]
A  && A \\
{} & \unit\tensor A  
\ncline{->}{1,1}{1,3}\naput{\id}
\ncline{->}{1,1}{2,2}\nbput{\unitl}
\ncline{->}{2,2}{1,3}\nbput{\unitll}
\end{psmatrix}
\hspace{13ex}
\begin{psmatrix}[ref=c,nodesep=1ex,colsep=4ex]
A  && A \\
{} & A\tensor\unit  
\ncline{->}{1,1}{1,3}\naput{\id}
\ncline{->}{1,1}{2,2}\nbput{\unitr}
\ncline{->}{2,2}{1,3}\nbput{\unitrr}
\end{psmatrix}
\makebox[0pt][l]{\raisebox{5ex}{\hspace*{10.8ex}(required)}}
\end{displaymath}\vspace{.8ex}%
We drop the requirement that $\unitll_A;\unitl_A=\id_{\unit\tensor A}$
and $\unitrr_A;\unitr_A=\id_{A\tensor\unit}$.\vspace{0.8ex}%
\begin{displaymath}
\begin{psmatrix}[ref=c,nodesep=1ex,colsep=4ex]
\unit\tensor A  && \unit\tensor A \\
{} & A  
\ncline{->}{1,1}{1,3}\naput{\id}
\ncline{->}{1,1}{2,2}\nbput{\unitll}
\ncline{->}{2,2}{1,3}\nbput{\unitl}
\end{psmatrix}
\hspace{10ex}
\begin{psmatrix}[ref=c,nodesep=1ex,colsep=4ex]
A\tensor\unit  && A\tensor\unit \\
{} & A
\ncline{->}{1,1}{1,3}\naput{\id}
\ncline{->}{1,1}{2,2}\nbput{\unitrr}
\ncline{->}{2,2}{1,3}\nbput{\unitr}
\end{psmatrix}
\makebox[0pt][l]{\raisebox{5ex}{\hspace*{8.8ex}(dropped)}}
\end{displaymath}
Although $\A$-linkings have yet to be defined formally, the following
example should nonetheless help to motivate the definition.
\begin{displaymath}
\begin{psmatrix}[labelsep=2.2ex,rowsep=1.8\baselineskip,mcol=r,colsep=2.5\baselineskip]% mcol=r => right-justify nodes
\rnode{aaa}{a} && \rnode{aaaa}{a} &&& \one{u}\tensor\rnode{a}{a} &&\one{uu}\tensor\rnode{aa}{a} 
\\
\one{u'}\tensor\rnode{aaa'}{a} & = &&&&\rnode{a'}{a} & \neq 
\\
\rnode{aaa''}{a} && \rnode{aaaa''}{a} &&& \one{u''}\tensor\rnode{a''}{a} &&\one{uu''}\tensor\rnode{aa''}{a}
\dvec{a}{a'}\naput{\unitll}
\dvec{a'}{a''}\naput{\unitl}
\dvec{aa}{aa''}
\dvec{uu}{uu''}\nbput{\id}
\dvec{aaa}{aaa'}\naput{\unitl}
\dvec{aaa'}{aaa''}\naput{\unitll}
\dvec{aaaa}{aaaa''}\nbput{\id}
\vecangles{u}{a'}{-75}{125}
\vecangles{u'}{aaa''}{-75}{125}
\end{psmatrix}
\end{displaymath}\vspace{.1ex}

\section{Linkings}\label{linkings}

This section presents the split star-autonomous category $\linka$ of
$\A$-linkings over a category $\A$.  Each linking is two-sided, being
a morphism $S\to T$ between two star-autonomous shapes $S$ and $T$,
analogous to the original Eilenberg-Kelly-MacLane graphs
\cite{EK66,KM71}.  Section~\ref{one-sided} introduces one-sided
linkings, more analogous to the graphs in \cite{KL80} for compact
closed categories.  Auxiliary one-sided linkings will facilitate later
proofs.

\subsection{The category $\lf$ of partial leaf functions between signed sets}\label{lf}

Define the category $\lf$ as follows.  An object is a \defn{signed
  set} $X$, whose elements we shall call \defn{leaves}, each signed
either \defn{positive} or \defn{negative}.  (Thus a signed set is a
set $X$ equipped with a function $X\to\{+,-\}$.)  A morphism $X\to Y$
is a \defn{partial leaf function}: a partial function from
$X\posrestr+Y\negrestr$ to $X\negrestr+Y\posrestr$, where $+$ is
disjoint union and $(-)\posrestr$ (resp.\ $(-)\negrestr$) is the
operation which restricts a signed set to its positive
(resp.\ negative) leaves.  For example,
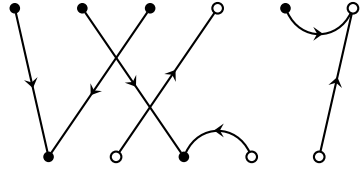
\begin{center}\vspace{.2ex}\begin{math}
\psset{nodesep=-.2mm}
\begin{psmatrix}
\posb{u}\bgap\posb{uu}\bgap\posb{a}\bgap\negb{aa}\bgap\posb{uuu}\bgap\negb{uuuu} \\
\posb{a1}\bgap\negb{aa1}\bgap\posb{b1}\bgap\negb{bb1}\bgap\negb{u1}
\linevec{u}{a1}
\linevecpos{uu}{b1}{.53}
\linevecpos{a}{a1}{.59}
\linevecpos{aa1}{aa}{.59}
\uloopvecanglesheight{bb1}{b1}{115}{65}{1}
\dloopvecanglesheight{uuu}{uuuu}{-65}{-115}{1}
\linevec{u1}{uuuu}
\end{psmatrix}
\end{math}\vspace{.4ex}\end{center}
is a (total) morphism from the upper signed set, 4 positive $\bullet$
and 2 negative $\circ$ leaves, to the lower one, 2 positive and 3
negative leaves.  Composition is simply finite (directed) path
composition:
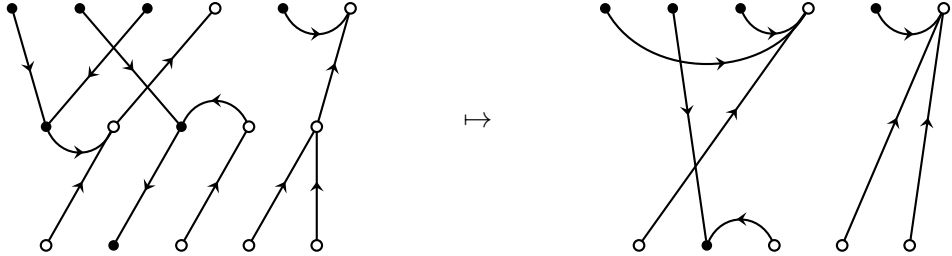
\begin{center}\label{goi-comp}\vspace{.8ex}
\psset{rowsep=1.1cm,nodesep=-.2mm}\begin{psmatrix}
$\posb{u}\bgap\posb{uu}\bgap\posb{a}\bgap\negb{aa}\bgap\posb{uuu}\bgap\negb{uuuu}$ 
&&
$\posb{u'}\bgap\posb{uu'}\bgap\posb{a'}\bgap\negb{aa'}\bgap\posb{uuu'}\bgap\negb{uuuu'}$
\\
$\posb{a1}\bgap\negb{aa1}\bgap\posb{b1}\bgap\negb{bb1}\bgap\negb{u1}$
&
$\mapsto$ \\
$\negb{u2}\gposb{b2}\gnegb{bb2}\gnegb{uu2}\gnegb{uuu2}$
&&
$\negb{u2'}\gposb{b2'}\gnegb{bb2'}\gnegb{uu2'}\gnegb{uuu2'}$
\linevec{u}{a1}
\linevecpos{uu}{b1}{.53}
\linevecpos{a}{a1}{.59}
\linevecpos{aa1}{aa}{.59}
\uloopvecanglesheight{bb1}{b1}{115}{65}{1}
\dloopvecanglesheight{uuu}{uuuu}{-65}{-115}{1}
\linevec{u1}{uuuu}
\dloopvecanglesheight{a1}{aa1}{-65}{-115}{1}
\linevec{u2}{aa1}
\linevec{b1}{b2}
\linevec{bb2}{bb1}
\linevec{uu2}{u1}
\linevec{uuu2}{u1}
\vecanglesposheight{u'}{aa'}{-60}{-120}{.58}{.8}
\linevecpos{uu'}{b2'}{.45}
\dloopvecanglesheight{a'}{aa'}{-65}{-120}{1}
\dloopvecanglesheight{uuu'}{uuuu'}{-65}{-115}{1}
\linevecpos{u2'}{aa'}{.58}
\uloopvecanglesheight{bb2'}{b2'}{115}{65}{1}
\linevec{uu2'}{uuuu'}
\linevec{uuu2'}{uuuu'}
\end{psmatrix}\vspace{1.4ex}\end{center}
Formally, given partial leaf functions $f:X\to Y$ and $g:Y\to Z$
define $\seqcomp f g\,:\,X\to Z$ by
$(\seqcomp f g)(l)=l'$ iff there is a finite directed path from $l$ to
$l'$ in the union of $f$ and $g$, viewed as a directed graph on
$X+Y+Z$.  The following proposition guarantees that $\seqcomp f g$ is
well-defined (single-valued).
\begin{proposition}[Unique Path Property]\label{upp}
If $\seqcomp f g$ is a composite partial leaf function containing the edge
$\langle l,l'\rangle$ (\ie, $(\seqcomp f g)(l)=l'$), then a unique path
$ll_0\ldots l_nl'$ gave rise to $\langle l,l'\rangle$ during
composition.
\end{proposition}
\begin{proof}
Suppose $ll_0'\ldots l_m'l'$ were an alternative path.  Let $i$ be
minimal with $l_i\neq l'_i$.  Then $f(l_{i-1})=l_i$ and
$f(l_{i-1})=l_i'$, or $g(l_{i-1})=l_i$ and $g(l_{i-1})=l_i'$,
contradicting (partial) functionality.
\end{proof}
The category $\lf$ is $\intsetp$, the result of applying the feedback
construction $\int$ \cite{JSV96} (or geometry of interaction
construction \cite{Abr96,Gir89}) to the traced monoidal category
$\setp$ of sets and partial functions, with tensor as coproduct.  Thus
$\lf$ is compact closed, with tensor as disjoint union and the dual of
a signed set obtained by reversing signs.  The subcategory of $\lf$
whose objects are finite and morphisms are bijections is the category
of involutions defined in \cite[\S3]{KL80}, a modified presentation of
the fixed-point free involutions of \cite{EK66,KM71}.

\subsection{The category $\linka$ of $\A$-linkings over a set $\A$}\label{switchings}

Fix a set $\A=\{a,b,c,\ldots\}$ of generators. 
An \defn{atom} is any generator in $\A$ or the constant $\unit$.
An \defn{$\A$-shape} is an expression generated from atoms by binary
tensor $\tensor$ and unary dual $(-)\perp$, \eg\ $\big(a\tensor
(b\perp\tensor a\perpp)\big)\perpp{}\perp\tensor
\unit\perp$.  The \defn{sign} of an atom or tensor in a shape is \defn{positive} $+$ iff
it is under an even number of duals $(-)\perp$, otherwise
\defn{negative} $-$.  Here is a shape with signs subscripted:
\raisebox{2pt}{\footnotesize$\big(\underset{-}{\unit}\perpp\underset{-}{\tensor}
(\underset{+}{a}\perpp{}\perp\underset{-}{\tensor}(\underset{+}{\unit}\underset{+}{\tensor}
\underset{+}{b})\perp)\big)\perp$}\:.
Write $\leaves{S}$ for the underlying signed set of a shape $S$,
obtained from its leaves, \ie, its occurrences of atoms.  
A \defn{leaf function} $X\to Y$ between signed sets is a partial leaf
function $X\to Y$ which is total (\ie, defined on the whole of
$X^++Y^-$).
A \defn{leaf function} $f:S\to T$ between shapes is a leaf function
$f:\leaves{S}\to\leaves{T}$ between the underlying signed sets.  The
\defn{graph of $f$} is the disjoint union of the underlying parse
trees of $S$ and $T$ (trees labelled with atoms at the leaves and
$\tensor$ or $\ast$ at internal vertices) together with the edges of
$f$, undirected.  For example, if $f:S\to T$ is
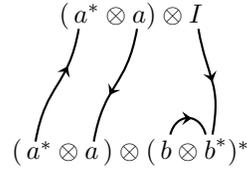
\begin{center}\vspace{-2.7ex}\small
\begin{psmatrix}[rowsep=3.3\baselineskip]
$(\,\rnode{a}{a}\perp\tensor\rnode{aa}{a})\tensor \one{uu}$ \\
$(\,\rnode{a1}{a}\perp\tensor\rnode{aa1}{a}\,) \tensor (\,\rnode{b1}{b}\tensor\rnode{bb1}{b}\perp)\perp$
\vecanglespos{uu}{bb1}{-80}{80}{.5}
\vecanglespos{aa}{aa1}{-100}{80}{.59}
\vecanglespos{a1}{a}{80}{-100}{.67}
\uloopvec{b1}{bb1}
\end{psmatrix}
\end{center}
then the graph of $f$ is shown below-left:
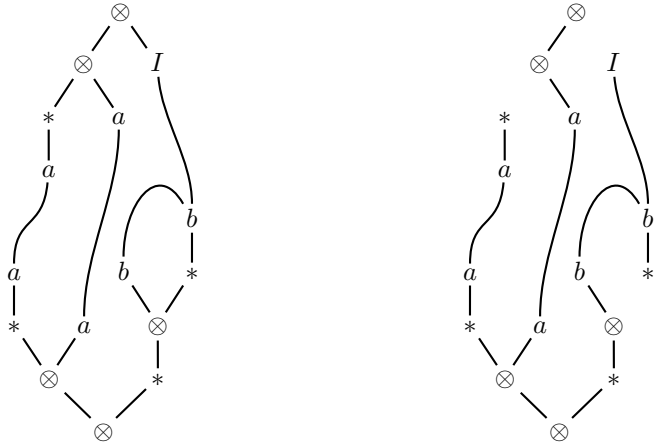
\begin{center}\vspace{-.4ex}\small\begin{math}\begin{array}{c}
\psset{colsep=3cm,framesep=0pt,nodesep=0pt,rowsep=1cm}
\psset{treemode=D,nodesep=2pt,levelsep=20pt}
  \pstree{\TR[name=t11]{\tensor}}
    {
     \pstree{\TR[name=t112]{\tensor}}
      {
       \pstree{\TR[name=ast1]{\ast}}
         {\TR[name=aa]{a}
         }
       \TR[name=a]{a}
      }
      {
       \TR[name=uu]{\unit} 
      }
    }
  
\\\\
\psset{treemode=U,nodesep=2pt,levelsep=20pt}
\pstree{\TR[name=tt1]{\tensor}}
  {\pstree{\TR[name=tt11]{\tensor}}
    {
     \pstree{\TR[name=ast'1]{\ast}}
       {\TR[name=aa1]{a}
       }
     \TR[name=a1]{a}
    }
     {\pstree{\TR[name=ast'2]{\ast}}
       {\pstree{\TR[name=tt121]{\tensor}}
          {\TR[name=bb1]{b}
           \pstree{\TR[name=ast'3]{\ast}}
            {\TR[name=b1]{b}
            }
           
          }
       }
     }
  }
\lineangles{uu}{b1}{-85}{90}
\lineangles{a}{a1}{-90}{90}
\lineanglesheight{aa1}{aa}{90}{-95}{1.1}
\lineanglesheight{b1}{bb1}{120}{90}{1.3}
\end{array}
\hspace{20ex}
\begin{array}{c}
\psset{colsep=3cm,framesep=0pt,nodesep=0pt,rowsep=1cm}
\psset{treemode=D,nodesep=2pt,levelsep=20pt}
  \pstree{\TR[name=t11]{\tensor}}
    {
     \pstree{\TR[name=t112]{\tensor}}
      {
       \pstree{\psset{linestyle=none}\TR[name=ast1]{\ast}}
         {\TR[name=aa]{a}
         }
       \TR[name=a]{a}
      }
      {
       \psset{linestyle=none}\TR[name=uu]{\unit} 
      }
    }
  
\\\\
\psset{treemode=U,nodesep=2pt,levelsep=20pt}
\pstree{\TR[name=tt1]{\tensor}}
  {\pstree{\TR[name=tt11]{\tensor}}
    {
     \pstree{\TR[name=ast'1]{\ast}}
       {\TR[name=aa1]{a}
       }
     \TR[name=a1]{a}
    }
     {\pstree{\TR[name=ast'2]{\ast}}
       {\pstree{\TR[name=tt121]{\tensor}}
          {\TR[name=bb1]{b}
           \pstree{\psset{linestyle=none}\TR[name=ast'3]{\ast}}
            {\TR[name=b1]{b}
            }
           
          }
       }
     }
  }
\lineangles{uu}{b1}{-85}{90}
\lineangles{a}{a1}{-90}{90}
\lineanglesheight{aa1}{aa}{90}{-95}{1.1}
\lineanglesheight{b1}{bb1}{120}{90}{1.3}
\end{array}\end{math}\end{center}
A \defn{switching} of a leaf function $f:S\to T$ between shapes is any
subgraph of the graph of $f$ obtained by deleting one of the two
argument edges of each positive tensor in $S$ and negative tensor in
$T$.  See above-right for an example.
A leaf function $f:S\to T$ is an \defn{$\A$-linking}
if it satisfies:
\begin{itemize}\label{criterion}
\myitem{1} \Matching. 
Restricting $f$ to $a$-labelled leaves (both in $S$ and in $T$) yields
a bijection for each generator $a\in\A$
\myitem{2} \Switching. 
Every switching of $f$ is a tree.
\end{itemize}
The leaf function above-left is an $\A$-linking: its switching shown
above-right is a tree, as are its seven other switchings.  Other
examples of $\A$-linkings are depicted on page~\ref{firstexample}.

The linking criterion (conditions (1) and (2)) is the analogue of
allowability in \cite{EK66,KM71}, presented combinatorially rather
than inductively/syntactically.  It is traditional in linear logic
\cite{Gir87} to present allowability combinatorially.  The criterion
above derives from a standard one for multiplicative proof nets
\cite{DR89,Hug12}.

\begin{proposition}\label{linear-time}
Verifying that a leaf function is a linking is linear time in
the number of leaves.
\end{proposition}
Thus the linking criterion (allowability) is very simple.
This proposition is proved in Section~\ref{proofs}, using one-sided
linkings.

Given linkings $f:S\to T$ and $g:T\to U$ define $\seqcomp f g:S\to U$ as their
path composite (in the category $\lf$ of partial leaf functions
between signed sets, defined in Section~\ref{lf}).
\begin{proposition}\label{link-comp}
The composite of two linkings is a linking.
\end{proposition}
The proof is in Section~\ref{proofs}, using one-sided linkings.

Write $\linka$ for the category of $\A$-linkings between $\A$-shapes.
Identities are inherited from $\lf$: the identity $S\to S$ has an edge
between the $i\nth$ leaf in the input shape and the $i\nth$ leaf in
the output shape.  The identity is a well-defined linking (\ie, every
switching is a tree) by a simple induction on the number of tensors in
$S$.

\parag{Compatibility}

Given leaf functions $f:X\to Y$ and $g:Y\to Z$ between signed sets,
write $f+g$ for the disjoint union of $f$ and $g$, viewed as a simple
directed graph on $X+Y+Z$.\footnote{A \defn{simple directed graph} on
a set $V$ (\cf\ \cite{Bol02}) is a set of edges on $V$, where an
\defn{edge on $V$} is an ordered pair $vw$ of distinct elements of
$V$.}
The following theorem will not be used in the sequel; we present it
since it is the analogue of Theorem~2.1 in \cite{KM71}.
\begin{theorem}[Compatibility]\label{compatibility}
If $f:S\to T$ and $g:T\to U$ are linkings between shapes, then $f+g$
contains no cycle.
\end{theorem}
By a \defn{cycle} we mean an undirected graph \cite{Bol02} on a vertex
set $\{v_1,\ldots,v_n\}$, all $v_i$ distinct, $n\ge 3$, and with an
edge $v_iv_j$ iff $j=i+1$ mod $n$.  The proof of the Compatibility
Theorem is in Section~\ref{proofs}, using one-sided linkings.

\subsection{$\linka$ is split star-autonomous}\label{la-split-star}
The category $\linka$ of $\A$-linkings is split star-autonomous, as
defined in Section~\ref{split-star}.  

Tensor and dual act symbolically on objects, \ie, the tensor of shapes
$S$ and $T$ is the shape $S\tensor T$, and the dual of $S$ is
$S\perp$.  Tensor acts as disjoint union on morphisms, hence is
functorial.  Given linkings $f:S\to T$ and $f':S'\to T'$, the leaf
function $f\tensor f':S\tensor S'\to T\tensor T'$ is a well-defined
linking since every switching of $f\tensor f'$ is a disjoint union of
switchings of $f$ and $f'$, connected at the tensor of $T\tensor T'$,
together with one of the two argument edges of the tensor of $S\tensor
S'$.  The dual $f\perp:T\perp\to S\perp$ of a linking $f:S\to T$ has
the same underlying directed graph as $f$, hence $(-)\perp$ is
functorial, full and faithful.

Associativity and symmetry are the obvious bijective leaf functions
(exactly the associativity and symmetry involutions of
\cite[\S3]{KM71}): associativity $(S\tensor T)\tensor U\to
S\tensor(T\tensor U)$ has edges between the $i\nth$ leaf of the input
and the $i\nth$ leaf of the output, and symmetry $S\tensor T\to
T\tensor S$ has edges between the $i\nth$ leaf of $S$ in $S\tensor T$
and the $i\nth$ leaf of $S$ in $T\tensor S$, and similarly for $T$.
Both are well-defined linkings by a simple induction (analogous to the
well-definedness of the identity).

The natural isomorphism $\linka(S\tensor
T,U\perp)\cong\linka(S,(T\tensor U)\perp)$ is the restriction of the
corresponding natural isomorphism
$\lf(\leaves{S}+\leaves{T},\leaves{U}\perp)\cong\lf(\leaves{S},(\leaves{Y}+
\leaves{Z})\perp)$ in the underlying compact closed category $\lf$ of leaf functions 
between signed sets.
This restriction is well-defined since switchings
$S\tensor T\to U\perp$ and 
$S\to (T\tensor U)\perp$ are in bijection.

Define $\unitl_S:S\to\unit\tensor S$ as the identity $S\to S$ together with
$\unit\tensor(-)$ added to the syntax of the output shape, and define
$\unitr_S:S\to S\tensor\unit$ similarly.  Since the added $\unit$ and $\tensor$
are positive, the switchings of $\unitl_S$ and $\unitr_S$ are in bijection with
those of the identity, hence $\unitl_S$ and $\unitr_S$ are well-defined
linkings.  The requisite triangle commutes since the edge between the
distinguished $\unit$'s in associativity $\assoc:(S\tensor \unit)\tensor T\to
S\tensor(\unit\tensor T)$ does not connect to an edge of $\unitr\tensor\id$
during composition.
Naturality of $\unitl_S$ and $\unitr_S$ holds because there is no edge to
the added $\unit$.

Define the retraction $\unitll_S:\unit\tensor S\to S$ from the
identity $S\to S$ by\label{split-unit-arrows} adding $\unit\tensor(-)$
to the input shape together with an edge from the added $\unit$ to an
arbitrary positive leaf of the $S$ on the right of the arrow
$\unit\tensor S\to S$ or a negative leaf of the $S$ on the left of the
arrow.  
This is a well-defined linking since every switching is a switching of
the identity $S\to S$ together with two new edges and two new
vertices, arranged so that the graph remains a tree, irrespective of
whether the added tensor is switched left or right.
We have $\seqcomp{\unitl_S}{\unitll_S}\;=\;\id_S$ 
since the edge of $\unitll_S$ from the added $\unit$ meets no edge of
$\unitl_S$ during the composition $\seqcomp{\unitl_S}{\unitll_S}$.
The retraction $\unitrr_S:S\tensor\unit\to S$ analogous.\footnote{Note
  that one cannot choose natural retractions: there are two candidates
  for $\unitll_{a\tensor a}$, and in either case, naturality
  $\:\seqcomp \unitll \sigma\:=\:\seqcomp{(\id\tensor
    \sigma)}{\unitll}\;$ fails for the symmetry map
  $\sigma\,:\,a\tensor a\,\to\,a\tensor a$.}

\subsection{The category $\linka$ of $\A$-linkings over an arbitrary base category $\A$}\label{arbitrary}

This section generalises $\A$-linkings from discrete $\A$ to an
arbitrary category $\A$.  Shapes are generated from the objects of
$\A$ as before.  A \defn{leaf function} $S\to T$ between shapes is a
partial leaf function $\leaves{S}\to\leaves{T}$ between the underlying
signed sets which is total, equipped with a labelling: every edge from
a leaf labelled by a generator (object of $\A$) is labelled with a
morphism of $\A$.  A leaf function $f:S\to T$ is an
\defn{$\A$-linking} if it satisfies:\pagebreak[1]
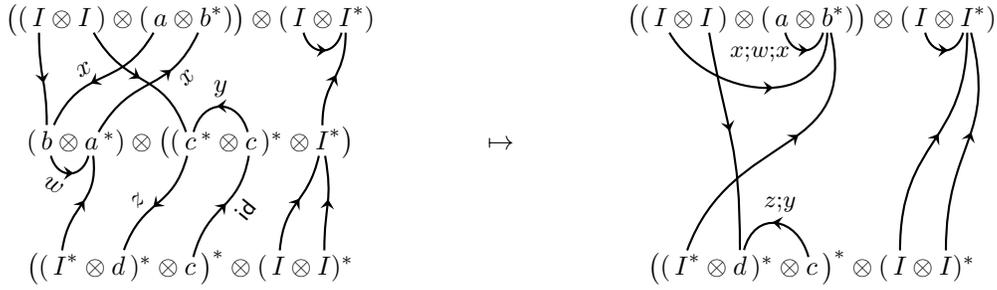
\begin{itemize}\label{arb-matching}
\myitem{1a} \Bijection. Restricting $f$ to $\A$-labelled leaves (both in $S$ and in $T$) yields a bijection.
\myitem{1b} \Labelling. If $x$ is the label of an edge 
from a leaf labelled $a$ to a leaf labelled $b$, then $x:a\to b$ is a
morphism in $\A$.
\myitem{2} \Switching. Every switching of $f$ is a tree.
\end{itemize}
(In forming the switchings of $f$, ignore edge labels.) Conditions
(1a) and (1b) reduced to (1) \Matching{} in the discrete case (since
all labels are identities).

For example, if $x:a\to b$ and $y:c\to c$ in $\A$, then
\begin{center}
\small
\begin{psmatrix}[rowsep=3.3\baselineskip]
$\big((\,\one{u}\tensor\one{uu}\,)\tensor(\,\rnode{a}{a}\tensor\rnode{b}{b}\perp)\big)
\tensor (\,\one{uuu}\tensor \one{uuuu}\perp)$ \\
$(\,\rnode{b1}{b}\tensor\rnode{a1}{a}\,\perp) \tensor \big((\,\rnode{c1}{c}\,\perp\tensor\rnode{cc1}{c}\,)\perp
\tensor \one{u1}\perp\big)$
\dvec{u}{b1}
\vecanglespos{uu}{c1}{-60}{110}{.47}
\vecanglespos{a}{b1}{-120}{65}{.59}\nbput[npos=.59,nrot=:D]{$x$}
\vecanglespos{a1}{b}{60}{-115}{.67}\nbput[npos=.76,nrot=:U]{$x$}
\uloopvecleft{cc1}{c1}\nbput[npos=.5,nrot=:D]{$y$}
\dloopvec{uuu}{uuuu}
\uvec{u1}{uuuu}
\end{psmatrix}
\end{center}
is a linking from the upper to the lower shape.

Composition is path composition, as in the discrete case, but
simultaneously collecting labels along each path and composing them in
$\A$.  More precisely, if $l_1\ldots l_n$ is a path traversed during
(underlying discrete) composition, resulting in the edge $\langle
l_0,l_n\rangle$ in the (underlying discrete) composite, then:
\begin{itemize}
\item if every edge $\langle l_{i-1},l_{i}\rangle$ is
labelled with a morphism $x_i$ in $\A$, the composite edge $\langle
l_0,l_n\rangle$ is labelled by the composite $x_n\ldots x_1$ in $\A$,
\item otherwise $\langle l_0,l_n\rangle$ is unlabelled.
\end{itemize}
Thus an edge in the composite is labelled iff every edge along the
path giving rise to it is labelled.  Here is an example of
composition:
\vspace{2ex}\begin{center}\small
\begin{psmatrix}[rowsep=2.8\baselineskip]
$\big((\,\one{u}\tensor\one{uu}\,)\tensor(\,\rnode{a}{a}\tensor\rnode{b}{b}\perp)\big)
\tensor (\,\one{uuu}\tensor \one{uuuu}\perp)$
&&
$\big((\,\one{u'}\tensor\one{uu'}\,)\tensor(\,\rnode{a'}{a}\tensor\rnode{b'}{b}\perp)\big)
\tensor (\,\one{uuu'}\tensor \one{uuuu'}\perp)$
\\
$(\,\rnode{b1}{b}\tensor\rnode{a1}{a}\,\perp) \tensor \big((\,\rnode{c1}{c}\,\perp\tensor\rnode{cc1}{c}\,)\perp
\tensor \one{u1}\perp\big)$
&
$\mapsto$ \\
$\big((\,\one{u2}\perp\tensor \rnode{d2}{d}\,)\perp\tensor \rnode{c2}{c}\,\big)\perp
\tensor (\,\one{uu2} \tensor \one{uuu2})\perp$
&&
$\big((\,\one{u2'}\perp\tensor \rnode{d2'}{d}\,)\perp\tensor \rnode{c2'}{c}\,\big)\perp
\tensor (\,\one{uu2'} \tensor \one{uuu2'})\perp$
\dvec{u}{b1}
\vecanglespos{uu}{c1}{-60}{110}{.47}
\vecanglespos{a}{b1}{-120}{65}{.59}\nbput[npos=.6,nrot=:D]{$x$}
\vecanglespos{a1}{b}{60}{-115}{.67}\nbput[npos=.76,nrot=:U]{$x$}
\uloopvecleft{cc1}{c1}\nbput[npos=.5,nrot=:D]{$y$}
\dloopvec{uuu}{uuuu}
\uvec{u1}{uuuu}
\dloopvec{b1}{a1}\nbput[npos=.32,nrot=:U]{$w$}
\vecangles{u2}{a1}{82}{-83}
\vecangles{c1}{d2}{-105}{70}\nbput[npos=.62,nrot=:U]{$z$}
\vecangles{c2}{cc1}{75}{-100}\nbput[npos=.63,nrot=:U]{\footnotesize$\id$}
\vecangles{uu2}{u1}{85}{-103}
\vecangles{uuu2}{u1}{90}{-80}
\vecanglesposheight{u'}{b'}{-60}{-90}{.58}{1}
\vecanglespos{uu'}{d2'}{-80}{90}{.45}
\dloopvec{a'}{b'}\nbput[npos=.2,labelsep=1pt]{\footnotesize$x;\!w;\!x\mkern-2mu$}
\dloopvec{uuu'}{uuuu'}
\vecanglespos{u2'}{b'}{70}{-77}{.58}
\uloopvecleft{c2'}{d2'}\nbput[npos=.5]{\footnotesize$z;\!y$}
\vecangles{uu2'}{uuuu'}{90}{-90}
\vecangles{uuu2'}{uuuu'}{90}{-76}
\end{psmatrix}\end{center}\vspace{2ex}
By the Unique Path Property (Proposition~\ref{upp}), composition is
well-defined: there is no ambiguity in constructing the labels on the
edges of a composite.
Composition is associative since concatenation of paths is
associative, and composition in $\A$ is associative.
The identity linking has all labels identities in $\A$.

\parag{Split star-autonomy}  
The split star-autonomous structure carries over from the discrete
case.  Labels on linking edges do not interfere: every label of a
canonical map (associativity, symmetry or unit map) is an identity.

\section{Nets}

This section quotients the category $\linka$ of $\A$-linkings by
Trimble rewiring \cite{Tri94,BCST96}, yielding the free star-autonomous
category generated by $\A$.  Define two $\A$-linkings as
\defn{similar} if one can be obtained from the other by re-targeting
an edge from an $\unit$,
\eg
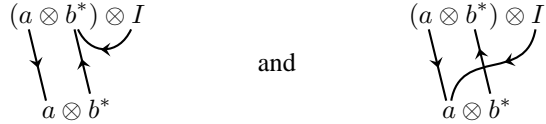
\begin{center}
\small
\begin{psmatrix}[rowsep=.5\baselineskip]
$(\rnode{a}{a}\tensor\rnode{b}{b}\perp)\tensor\one{u}$ &&
$(\rnode{a'}{a}\tensor\rnode{b'}{b}\perp)\tensor\one{u'}$ \\ &
and \\ $\rnode{a1}{a}\tensor\rnode{b1}{b}\perp$ &&
$\rnode{a1'}{a}\tensor\rnode{b1'}{b}\perp$
\linevec{a}{a1}
\linevec{b1}{b}
\nccurve[ncurv=1.1,angleA=-115,angleB=-65]{u}{b}
\lput{:0}{\psline{->}(.05,.005)(.1,.005)}
\linevec{a'}{a1'}
\linevecpos{b1'}{b'}{.7}
\vecanglespos{u'}{a1'}{-100}{80}{.4}
\end{psmatrix}
\end{center}
in the discrete case, or, if $x:a\to b$ and $y:c\to c$ are morphisms in $\A$,
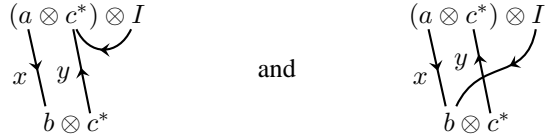
\begin{center}
\small
\begin{psmatrix}[rowsep=.7\baselineskip]
$(\rnode{a}{a}\tensor\rnode{b}{c}\perp)\tensor\one{u}$ &&
$(\rnode{a'}{a}\tensor\rnode{b'}{c}\perp)\tensor\one{u'}$ \\ &
and \\ $\rnode{a1}{b}\tensor\rnode{b1}{c}\perp$ &&
$\rnode{a1'}{b}\tensor\rnode{b1'}{c}\perp$
\linevec{a}{a1}\nbput[npos=.59]{$x$}
\linevec{b1}{b}\naput[npos=.5]{$y$}
\nccurve[ncurv=1.1,angleA=-115,angleB=-65]{u}{b}
\lput{:0}{\psline{->}(.05,.005)(.1,.005)}
\linevec{a'}{a1'}\nbput[npos=.59]{$x$}
\linevecpos{b1'}{b'}{.7}\naput[npos=.66]{$y$}
\vecanglespos{u'}{a1'}{-100}{65}{.4}
\end{psmatrix}
\end{center}
An \defn{$\A$-net} is an equivalence class of $\A$-linkings modulo
similarity (\ie, modulo the equivalence relation generated by
similarity).

\begin{theorem}[Net Compositionality]\label{net-comp}
Composition of $\A$-linkings respects equivalence.
\end{theorem}
In other words, if $f,f':S\to T$ and $g,g':T\to U$ are linkings with
$f$ equivalent to $f'$ and $g$ equivalent to $g'$, then the composite
linkings $\seqcomp f g$ and $\seqcomp{f'}{g'}:S\to U$ are equivalent.  This
theorem is proved in Section~\ref{net-comp-proof} via one-sided nets.
By the theorem, composition of $\A$-nets is well-defined: given proof
nets $[f]:S\to T$, and $[g]:T\to U$, where $[h]$ denotes the
equivalence class of a linking $h$, define
$\seqcomp{[f]}{[g]}=[\seqcomp f g]:S\to U$.  Write $\netsa$ for the
category of $\A$-nets.  Typically we abbreviate a net $[f]$ to $f$,
when it is clear from context that we are dealing with a net rather
than a linking.

\parag{Star autonomy} 
The split star-autonomous structure of the category $\linka$ of
$\A$-linkings respects equivalence, yielding
a star-autonomous structure on $\netsa$.

On morphisms, tensor in $\linka$ is disjoint union, hence respects
similarity in each argument, \ie, if linkings $f$ and $f'$ are
similar, then $f\tensor g$ and $f'\tensor g$ are similar, as are
$g\tensor f$ and $g\tensor f'$.  Duality on $\linka$ respects
similarity (if $f$ and $f'$ are similar then $f\perp$ and $f'{}\perp$
are similar) since
it acts trivially on the graph of a morphism (so re-targeting an edge
from an $\unit$ amounts to the same thing before and after dualising).
Similarly, the natural isomorphism $\linka(S\tensor
T,U\perp)\cong\linka(S,(T\tensor U)\perp)$ respects similarity, since
a linking $S\tensor T\to U\perp$ has the same underlying directed
graph as its transpose $S\to (T\tensor U)\perp$.

The split monomorphisms $\unitl_S:S\to \unit\tensor S$ and $\unitr_S:S\to S\tensor\unit$ in
$\linka$ become isomorphisms upon quotienting.
%
%Choose a left inverse $\unitll_S:\unit\tensor S\to S$ for $\unitl_S$.
The composite\vspace*{-1ex}%
\begin{displaymath}
  \begin{psmatrix}[colsep=6.5ex]
  \Rnode{a}{\unit\tensor S} & \Rnode{b}{S} & \Rnode{c}{\unit\tensor S} 
  \ncLine{->}{a}{b}\naput[npos=.6]{\unitll_S}\ncLine{->}{b}{c}\naput[npos=.4]{\unitl_S}
  \end{psmatrix}
\end{displaymath}
in $\linka$ differs from the identity $\unit\tensor S\to\unit \tensor
S$ by just one edge, the edge from the distinguished $\unit$ on the
left of the arrow, hence is similar to the identity.  Thus
$\seqcomp{\unitll_S}{\unitl_S}\;=\;\id_S$ in $\netsa$.  Similarly,
$\seqcomp{\unitrr_S}{\unitr_S}\;=\;\id_S$ in $\netsa$.

The associativity, symmetry and unit coherence diagrams commute in
$\linka$, hence also in $\netsa$.

\subsection{Free star-autonomous category and full coherence}\label{freestar}
\begin{theorem}[freeness]\label{free}
For any category $\A$, the category $\netsa$ of $\A$-nets 
is the free star-autono\-mous category generated by $\A$.
\end{theorem}
The proof is the subject of Section~\ref{mainproof}.
By finiteness of equivalence classes, we have:
\begin{theorem}[Full Coherence]
Equality of morphisms in the free star-autonomous category generated
by a category is decidable.
\end{theorem}
This theorem was first proved using the rewiring approach in
\cite{BCST96} (in a more general form, over a polygraph with
equations, and with \emph{star-autonomous category} axiomatised as a
symmetric linearly distributive category with negation).

Two examples were given in the Introduction
(p.~\pageref{intro-rewiring-examples}).  Three more are provided
below.

\emph{Example 3: identity = twist on $\unit$.}\label{example3}
Example~1 (p.~\pageref{intro-rewiring-example1}) proved
$\id_{\bot\tensor\bot}\neq\tw_{\bot\tensor\bot}:\bot\tensor\bot\to\bot\tensor\bot$.
The following pair of rewirings shows
$\id_{\unit\tensor\unit}=\tw_{\unit\tensor\unit}:\unit\tensor\unit\to\unit\tensor\unit$.
\begin{center}\vspace{3ex}\small\psset{rowsep=4\baselineskip}
\renewcommand{\tensor}{\mkern1mu\otimes\mkern1mu}
\newcommand{\core}[2]{%
\rnode{#1}{\psframebox[linestyle=none]{%
\psset{nodesep=.5ex,rowsep=2.5\baselineskip}\begin{psmatrix}
$\one{a}\tensor\one{b}$
\\
$\one{aa}\tensor\one{bb}$
#2
\end{psmatrix}}}}
\begin{psmatrix}
  \core{D}{%
    \dvec{a}{aa}
    \dvec{b}{bb}
  }
  & \hspace{-3ex} &
  \core{C}{%
  \vecanglespos{a}{bb}{-75}{105}{.65}
  \dvec{b}{bb}
  }
  & \hspace{-4ex} &
  \core{B}{%
  \vecanglespos{a}{bb}{-75}{105}{.78}
  \vecanglespos{b}{aa}{-105}{75}{.78}
  }
\end{psmatrix}
  \psset{linestyle=none}
  \ncline{B}{C}\ncput{$\leftrightarrow$}
  \ncline{C}{D}\ncput{$\leftrightarrow$}
\vspace{3ex}\end{center} Compare this with \cite[end of \S3.1]{BCST96},
which proves the dual result:
%$\id_{\unit\tensor\unit}=\tw_{\unit\tensor\unit}$ 
$\id_{\bot\scriptparr\bot}=\tw_{\bot\scriptparr\bot}$.

\emph{Example 4: Triple-dual problem.}\label{example4} We show that
the following diagram commutes (like Example~2 on
page~\pageref{intro-rewiring-example2}, an instance of the triple-dual
problem: see \cite[\S1, diagram\ (1.4)]{KM71})\vspace{2ex}
\begin{equation}\label{triangle2}
\psset{nodesep=1ex}
\vspace{2ex}\begin{psmatrix}[ref=c]
\big((\unit\multimap\unit)\multimap\unit\big)\multimap\unit &
\unit\multimap\unit\rule{0pt}{3ex} \\
& 
\big((\unit\multimap\unit)\multimap\unit\big)\multimap\unit
\ncline{->}{1,1}{1,2}\naput{k_\unit\multimap \id}
\ncline{->}{1,2}{2,2}\naput{k_{(\unit\multimap\mkern2mu\unit)}}
\ncline{->}{1,1}{2,2}\nbput{\id}
\end{psmatrix}\vspace{1ex}
\end{equation}
where 
$A\multimap B=(A\tensor B\perp)\perp$ and
$k_A:A\;\to\;(A\mkern-2mu\multimap\mkern-2mu
A)\mkern-2mu\multimap\mkern-2mu A$ is the canonical map of its type.
Each path in the triangle determines a corresponding linking
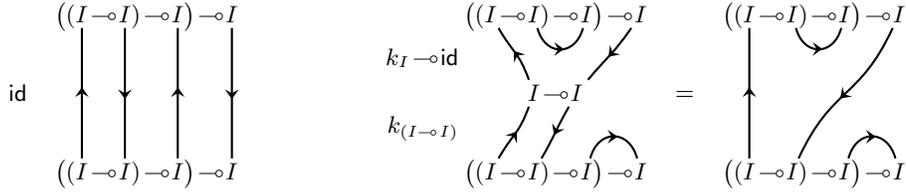
\begin{center}\vspace{3ex}\footnotesize\begin{math}
\newcommand{\tightmap}{\mkern-3mu\multimap\mkern-4mu}
\psset{rowsep=1.7\baselineskip,nodesep=2pt}
\rule{.3cm}{0pt}\begin{psmatrix}
\big((\one{a}\tightmap\one{b})\tightmap\one{bb}\big)\tightmap\one{bbb}\\\\
\big((\one{a1}\tightmap\one{b1})\tightmap\one{bb1}\big)\tightmap\one{bbb1}
\uvec{a1}{a}
\dvec{b}{b1}
\uvec{bb1}{bb}
\dvec{bbb}{bbb1}
\blabel{a}{a1}{\id}{.7cm} 
\end{psmatrix}
\rule{3cm}{0cm}
\rnode{A}{\begin{psmatrix} 
\big((\one{a}\tightmap\one{b})\tightmap\one{bb}\big)\tightmap\one{bbb}\\
\one{am}\tightmap\one{bm} \\
\big((\one{a1}\tightmap\one{b1})\tightmap\one{bb1}\big)\tightmap\one{bbb1}
\blabelpos{a}{a1}{k_\unit\tightmap \id}{.4cm}{.23} 
\blabelpos{a}{a1}{k_{(\unit\multimap\mkern2mu\unit)}}{.4cm}{.77} 
\vecangles{am}{a}{110}{-60}
\vecangles{a1}{am}{60}{-110}
\dloopvec{b}{bb}
\vecangles{bbb}{bm}{-125}{50}
\vecangles{bm}{b1}{-120}{60}
\uloopvec{bb1}{bbb1}
\end{psmatrix}}
\rule{1cm}{0cm}
\rnode{B}{\begin{psmatrix} 
\big((\one{a}\tightmap\one{b})\tightmap\one{bb}\big)\tightmap\one{bbb}\\
\\
\big((\one{a1}\tightmap\one{b1})\tightmap\one{bb1}\big)\tightmap\one{bbb1}
\uvec{a1}{a}
\dloopvec{b}{bb}
\vecangles{bbb}{b1}{-115}{65}
\uloopvec{bb1}{bbb1}
\end{psmatrix}}
\clabel{A}{B}{=}
\end{math}\vspace{3ex}\end{center}
equivalent via five rewirings:
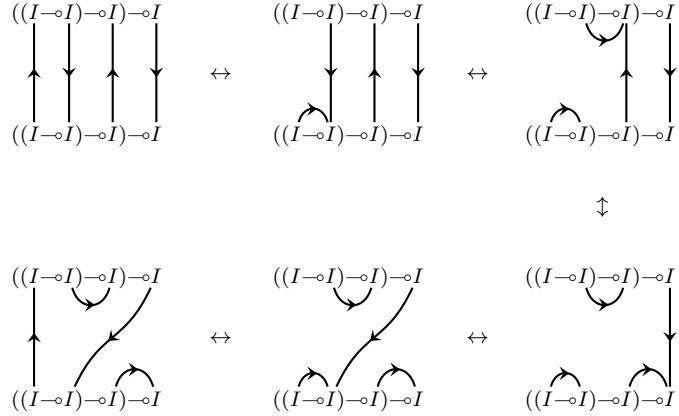
\begin{center}\vspace{3ex}\scriptsize
\psset{nodesep=.5ex,rowsep=1.3cm}
\newcommand{\tightmap}{\mkern-5mu\multimap\mkern-5mu}
\newcommand{\core}[2]{%
\rnode{#1}{\psframebox[linestyle=none]{%
\psset{nodesep=.5ex,rowsep=1.3cm}\begin{psmatrix}
$((\one{a}\tightmap\one{b})\tightmap\one{bb})\tightmap\one{bbb}\;$
\\
$((\one{a1}\tightmap\one{b1})\tightmap\one{bb1})\tightmap\one{bbb1}\;$
#2
\end{psmatrix}}}}
\begin{psmatrix}
  \core{D}{%
    \uvec{a1}{a}
    \dvec{b}{b1}
    \uvec{bb1}{bb}
    \dvec{bbb}{bbb1}
  }
  & \hspace{-3ex} &
  \core{C}{%
    \uloopvec{a1}{b1}
    \dvec{b}{b1}
    \uvec{bb1}{bb}
    \dvec{bbb}{bbb1}
  }
  & \hspace{-4ex} &
  \core{B}{%
    \uloopvec{a1}{b1}
    \dloopvec{b}{bb}
    \uvec{bb1}{bb}
    \dvec{bbb}{bbb1}
  }
\\
  \core{X}{%
	\uvec{a1}{a}
	\dloopvec{b}{bb}
	\vecangles{bbb}{b1}{-115}{65}
	\uloopvec{bb1}{bbb1}
  }
  & \hspace{-4ex} &
  \core{Y}{%
	\uloopvec{a1}{b1}
	\dloopvec{b}{bb}
	\vecangles{bbb}{b1}{-115}{65}
	\uloopvec{bb1}{bbb1}
  }
  & \hspace{-4ex} &
  \core{Z}{%
	\uloopvec{a1}{b1}
	\dloopvec{b}{bb}
        \dvec{bbb}{bbb1}
	\uloopvec{bb1}{bbb1}
  }
\end{psmatrix}
  \psset{linestyle=none}
  \ncline{B}{C}\ncput{$\leftrightarrow$}
  \ncline{C}{D}\ncput{$\leftrightarrow$}
  \ncline{B}{Z}\ncput[npos=.5,nrot=:D]{$\leftrightarrow$}
  \ncline{X}{Y}\ncput{$\leftrightarrow$}
  \ncline{Y}{Z}\ncput{$\leftrightarrow$}
\vspace{3ex}\end{center} We conclude that triangle (\ref{triangle2})
commutes in every star-autonomous category.  Compare this with
\cite[Fig.~3]{BCST96}. As with Example~2, a key advantage here is that
the composite $\:k_\unit\mkern-3mu\multimap\mkern-3mu\id\:;\:
k_{(\unit\multimap\unit)}\:$ is immediate on inspection, since we have
simple path composition.

\emph{Example 5: more path composition.}\label{path-power}
The following example illustrates a larger path composition.
Let $A\,\smallparr\, B\:=\:(A\perp\tensor B\perp)\perp$ and $\bot=I\perp$.
The path composition of linkings
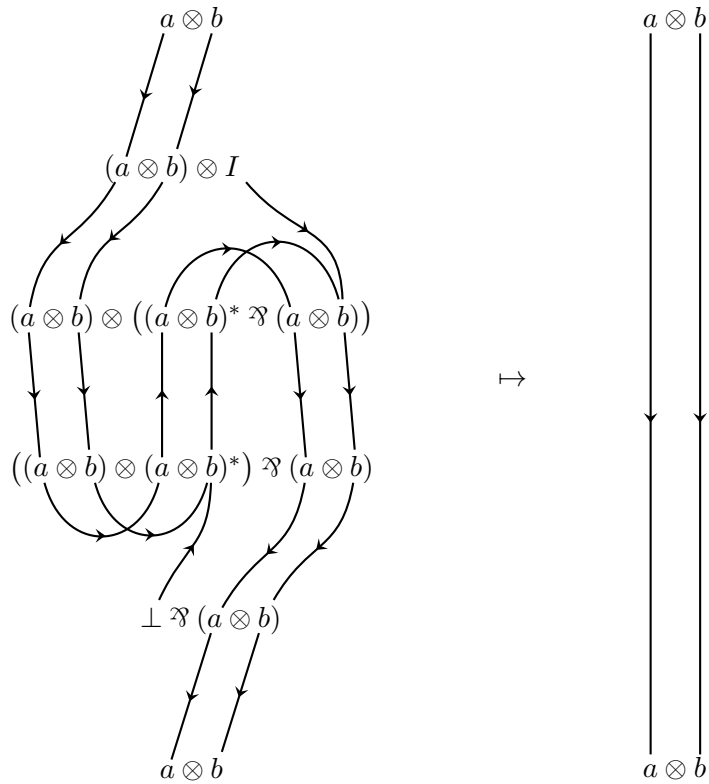
\begin{center}\vspace{3ex}\hspace{10ex}\psset{nodesep=.5ex,yunit=2cm}\begin{pspicture}[.5](0,-.2)(5,5.3)
\rput(0,0){\hspace{0ex}$\rnode{a}{a} \tensor \rnode{b}{b}$}
\rput(0,1){\hspace{3ex}$\rnode{u1}{\bot} 
\,\smallparr\, (\rnode{a1}{a} \tensor \rnode{b1}{b})$}
\rput(0,2){$\big((\rnode{aa2}{a}\tensor \rnode{bb2}{b})\tensor(\rnode{aaa2}{a}\tensor\rnode{bbb2}{b})\perp\big)
\,\smallparr\, (\rnode{a2}{a} \tensor \rnode{b2}{b})$}
\rput(0,3){$(\rnode{aa3}{a}\tensor \rnode{bb3}{b}) \tensor \big((\rnode{aaa3}{a}\tensor\rnode{bbb3}{b})\perp
\,\smallparr\, (\rnode{a3}{a} \tensor \rnode{b3}{b})\big)$}
\rput(0,4){$(\rnode{a4}{a} \tensor \rnode{b4}{b}) \tensor \rnode{u4}{\unit}$\hspace{3ex}}
\rput(0,5){$\rnode{a5}{a} \tensor \rnode{b5}{b}$\hspace{0ex}}
\vecanglespos{u4}{b3}{-50}{92}{.5}% top I downwards arrow
\vecanglespos{u1}{bbb2}{65}{-92}{.5}% bottom _|_ upwards arrow
\linevec{a1}{a}
\linevec{b1}{b}
\vecangles{a2}{a1}{-100}{60}
\vecangles{b2}{b1}{-100}{60}
\dloopvec{aa2}{aaa2}
\dloopvec{bb2}{bbb2}
\linevec{aa3}{aa2}
\linevec{bb3}{bb2}
\linevec{aaa2}{aaa3}
\linevec{bbb2}{bbb3}
\linevec{a3}{a2}
\linevec{b3}{b2}
\uloopvec{aaa3}{a3}
\uloopvec{bbb3}{b3}
\vecangles{a4}{aa3}{-120}{80}
\vecangles{b4}{bb3}{-120}{80}
\linevec{a5}{a4}
\linevec{b5}{b4}
\end{pspicture}
\hspace{8ex}
$\mapsto$
\hspace{3ex}
\begin{pspicture}[.5](0,-.2)(5,5.3)
\psset{nodesep=.5ex}
\rput(2.5,0){$\rnode{a}{a\strut} \tensor \rnode{b}{b\strut}$}
\rput(2.5,5){$\rnode{a5}{a} \tensor \rnode{b5}{b}$}
\psset{nodesepB=-1pt}
\linevec{a5}{a}
\linevec{b5}{b}
\end{pspicture}\vspace{2ex}\end{center}
shows that the following diagram commutes in every star-autonomous
category, where each map is canonical at its type:
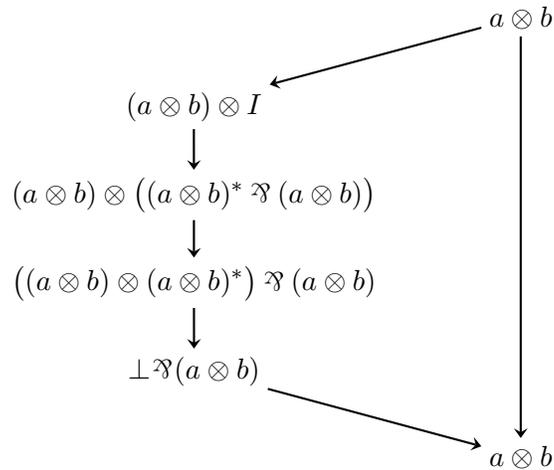
\begin{center}\begin{math}
\begin{psmatrix}[ref=c,nodesep=.6ex,rowsep=4ex]
&
a\tensor b
\\
(a\tensor b)\tensor \unit
\\
(a\tensor b)\tensor \big((a\tensor b)\perp\,\smallparr\,(a\tensor b)\big)
\\
\big((a\tensor b)\tensor(a\tensor b)\perp \big) \,\smallparr\, (a\tensor b)
\\
\bot\smallparr(a\tensor b)
\\
& a\tensor b
\ncline{->}{1,2}{6,2}
\ncline{->}{1,2}{2,1}
\ncline{->}{2,1}{3,1}
\ncline{->}{3,1}{4,1}
\ncline{->}{4,1}{5,1}
\ncline{->}{5,1}{6,2}
\end{psmatrix}\end{math}\end{center}
This example is, of course, rather contrived.  (What other canonical
map $a\tensor b\to a\tensor b$ could there be?)  The point is to give
an example of path composition in a large diagram.

\section{One-sided linkings}\label{one-sided}

Eilenberg-Kelly-MacLane graphs \cite{EK66,KM71} are \emph{two-sided}
in the sense that they are between a source and a target.  Taking
advantage of duality, Kelly-Laplaza graphs \cite{KL80} for compact
closed categories are \emph{one-sided}: the authors define a
(fixed-point free) involution on a single signed set, hence on a
single shape\footnote{We refer to words built freely from generators
  and the constant $\unit$ by binary $\tensor$ and unary $(-)\perp$ as
  \emph{shapes}, following the Eilenberg-Kelly-MacLane terminology for
  similar freely generated expressions.  These shapes are exactly the
  words defined in \cite[\S3]{KL80}, objects of the free compact
  closed category.}.  They then define a two-sided involution, \ie, a
morphism $S\to T$, as a (one-sided) involution on $S\perp\tensor T$.

Since star-autonomous categories are more general than compact closed
categories, to obtain a directly analogous one-sided representation
one would define a two-sided linking $S\to T$ as a one-sided linking
on the shape $(S\tensor T\perp)\perp$.  However, to avoid the extra
baggage of two auxiliary duals $(-)\perp$ and a tensor $\tensor$, we
shall instead define one-sided linkings on \emph{sequents}
$\;S_1,\ldots,S_n\;$ of shapes $S_i$, and define a two-sided linking
$S\to T$ as a one-sided linking on the two-shape sequent $\hspace{1.5pt}S\perp\!,T\hspace{1.5pt}$.  Using
the sequent calculus style, we can also define explicit cuts between
shapes, facilitating an inductive proof that two-sided linkings
compose, \ie, that the path composite of two (two-sided) linkings is a
well-defined (two-sided) linking.

The material here on one-sided linkings amounts to the multiplicative
proof nets with units in \cite{Hug12}, with explicit negation
$(-)\perp$, and with axiom links between dual occurrences of
generators in $\A$ generalised to morphisms of $\A$.  However, we
shall not require any familiarity on the part of the reader with proof
nets or linear logic \cite{Gir87}.

\parag{Sequents} 
The following definitions are a mild generalisation of those in
Section~2 of \cite{Hug12}.

Fix a set $\A=\{a,b,\ldots\}$ of generators.  Henceforth identify a
shape (generated from $\A$) with its parse tree, a tree labelled with
atoms (generators or the constant $\unit$) at the leaves, and
$\tensor$ and $\ast$ at internal vertices.  (Examples of parse trees
were shown in Section~\ref{switchings} (page~\pageref{switchings}), in
the example of a switching.)  A \defn{sequent} is a non-empty disjoint
union of shapes.  Thus a sequent is a particular kind of labelled
forest.  We take $S,T,\ldots$ to range over shapes, and
$\Gamma\mkern-3mu,\Delta,\ldots$ to range over (possibly empty)
disjoint unions of shapes.  A \defn{cut pair} $\cutpair{S}$ is a
disjoint union of a shape $S$ with its dual $S\perp$
%complementary shapes $S$ and $S\perp$ (the latter
%shape equal to the former with an extra root $\ast$) 
together with an undirected edge, a \defn{cut}, between the root of
$S$ and the root of $S\perp$, \eg
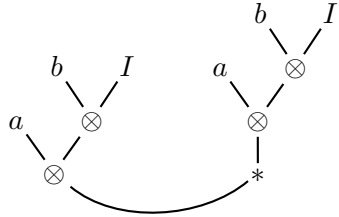
\begin{displaymath}\hspace{1ex}\psset{colsep=2.5cm,framesep=0pt,nodesep=0pt,rowsep=1cm}
\renewcommand{\S}[1]{\pstree[treemode=U,nodesep=2pt,levelsep=20pt]{\TR[name=#1]{\tensor}}
  {\TR{a}\pstree{\TR{\tensor}}
      {\TR{b}\TR{\unit}}}}
\S{t}
\hspace{6ex}
\pstree[treemode=U,nodesep=2pt,levelsep=20pt]{\TR[name=ast]{\ast}}{\S{tt}}
\nccurve[nodesepA=2pt,nodesepB=1pt,angleA=-40,angleB=-140]{t}{ast}
\vspace{3ex}
\end{displaymath}
for $S=a\tensor(b\tensor\unit)$.
A \defn{cut sequent} is a disjoint
union of a sequent and zero or more cut pairs.  A \defn{switching} of
a cut sequent is any subgraph obtained by deleting one of the two
argument edges of each negative $\tensor$
(\cf\ Section~\ref{switchings}).  We use comma to denote disjoint
union, \ie\ $\;S_1,\ldots,S_n\;$ is the disjoint union of the shapes
$S_i$.

\parag{Linkings}  We begin with the discrete case, $\A$ a set of generators.  For the more general 
case of $\A$ an arbitrary category, see Section~\ref{arb} below.

A \defn{leaf function} on a cut sequent is a function from its
negative leaves to its positive leaves.  An \defn{$\A$-linking} on a
cut sequent $\Gamma$ is a leaf function $f$ on $\Gamma$ satisfying:
\begin{itemize}\label{onesided-linking-criterion}
\myitem{1} \Matching. For any generator $a\in\A$, the restriction of $f$ to $a$-labelled leaves
is a bijection between the negative $a$-labelled leaves of $\Gamma$
and the positive $a$-labelled leaves of $\Gamma$.
\myitem{2} \Switching. For any switching $\Gamma'$ of $\Gamma$, the undirected graph obtained by adding the 
edges of $f$ to $\Gamma'$ is a tree.
\end{itemize}
Thus two-sided linkings $S\to T$ between shapes, as defined in
Section~\ref{linkings}, are in bijection with one-sided linkings on
the two-shape sequent $S\perp,T$.

\subsection{Cut elimination}
Let $f$ be a linking on the cut sequent $\Gamma,\cutpair{S}$.
The result $f'$ of \defn{eliminating} the cut in $\cutpair{S}$
is:
\begin{itemize}
\item \emph{Atom.}\label{atom-cut}
Suppose $S$ is an atom $\alpha$ (a generator $a\in\A$ or the constant $\unit$).
Thus (in parse tree terms) the cut pair $\cutpair{S}$ comprises leaves
$l^+,l^-$ labelled $\alpha$,
a vertex $v$ labelled $\ast$, an argument edge from $v$ to $l^-$, and a
cut edge between $v$ and $l^+$.
Delete the cut pair (\ie, the vertices $l^+,l^-,v$ and associated
edges), and reset every $f$-edge to $l^+$ to target $f(l^-)$ instead,
\ie, for all negative leaves $n$ of the sequent such that $f(n)=l^+$,
set $f(n)=f(l^-)$.  
Schematically,
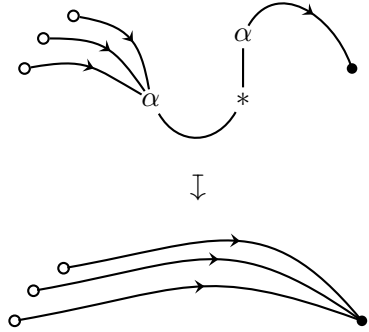
\begin{center}
\hspace{8ex}\psset{colsep=2.5cm,framesep=0pt,nodesep=0pt,rowsep=1.2cm}
\begin{psmatrix}
\Cnode(-3,.5){a}
\Cnode(-2.5,.9){b}
\Cnode(-1.7,1.2){c}
\circlenode[linestyle=none]{A}{$\alpha$}\rule{0pt}{7.2ex}
\hspace{4ex}
\circlenode[linestyle=none]{ast}{$\ast$}
\nput[labelsep=3ex]{90}{ast}{\circlenode[linestyle=none]{AA}{$\alpha$}}\ncline{ast}{AA}
{\psset{nodesep=1pt}\nccurve[ncurv=1,angleA=-60,angleB=-120]{A}{ast}}
\psset{nodesepB=-1pt}
\vecangles{a}{A}{10}{150}
\vecangles{b}{A}{0}{125}
\vecangles{c}{A}{-10}{100}
\hspace{4ex}
\Cnode*(1,.5){d}
\psset{nodesepB=-1pt,nodesepA=-1pt}
\vecangles{AA}{d}{60}{110}
\\
\rput{270}(-.8,.6){$\mapsto$}
\\
\Cnode(-3,.5){a}
\Cnode(-2.5,.9){b}
\Cnode(-1.7,1.2){c}
\hspace{14ex}
\psset{nodesepB=-1pt}
\Cnode*(1,.5){d}
\psset{nodesepB=-1pt,nodesepA=-1pt}
\vecangles{a}{d}{10}{160}
\vecangles{b}{d}{10}{145}
\vecangles{c}{d}{10}{130}
\end{psmatrix}
\end{center}
where the vertices $\circ$ represent all instances of $n$ and
$\bullet$ is $f(l^-)$.  (The left-most $\alpha$ in the picture is the
leaf $l^+$, the right-most $\alpha$ is the leaf $l^-$, and the $\ast$
is the vertex $v$.)
\item \emph{Tensor.} Suppose $S=T\tensor U$.  
Replace $\cutpair{S}$ by $\cutpair{T},\cutpair{U}$.  The leaves, and
$f$, remain unchanged (identifying the leaves of $T$ in $T\tensor U$
with $T$ in $\cutpair{T}$, \etc).
Schematically:
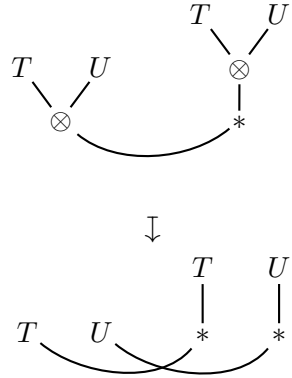
\begin{center}
\hspace{-6ex}\psset{colsep=2.5cm,framesep=0pt,nodesep=0pt,rowsep=1cm}
\begin{psmatrix}
\raisebox{1ex}{$
\pstree[treemode=U,nodesep=2pt,levelsep=20pt]
{\TR[name=t]{\tensor}}
  {\TR[name=b]{T}
   \TR[name=c]{U}}
\hspace{6ex}
\pstree[treemode=U,nodesep=2pt,levelsep=20pt]{\TR[name=ast]{\ast}}
  {\pstree{\TR[name=p]{\tensor}}
     {\TR[name=bb]{T}
      \TR[name=cc]{U}
     }
  }
$}
\nccurve[nodesepA=2pt,nodesepB=1pt,angleA=-40,angleB=-140]{t}{ast}
\\
\rput{270}(0,.2){$\mapsto$}
\\
\raisebox{5ex}{$
\hspace{.5ex}\rnode{b}{T}
\hspace{4ex}
\rnode{c}{U}
\hspace{6ex}
\circlenode[linestyle=none]{ast1}{\ast}
\nput[labelsep=3ex]{90}{ast1}{\circlenode[linestyle=none]{bb}{T}}\ncline{ast1}{bb}
\hspace{4ex}
\circlenode[linestyle=none]{ast2}{\ast}
\nput[labelsep=3ex]{90}{ast2}{\circlenode[linestyle=none]{cc}{U}}\ncline{ast2}{cc}
\hspace{-2ex}
\nccurve[nodesepA=1pt,nodesepB=2pt,angleA=-40,angleB=-130]{b}{ast1}
\nccurve[nodesepA=1pt,nodesepB=2pt,angleA=-40,angleB=-130]{c}{ast2}
$}
\end{psmatrix}
\end{center}
\item \emph{Dual.} Suppose $S=T\perp$.
Replace the cut pair $\dualcutpair{T}$ by $\revcutpair{\,T}$.  The
leaves, and $f$, remain unchanged.
Schematically,
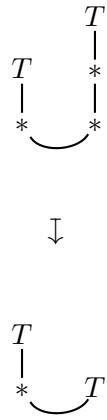
\begin{center}
\psset{colsep=2.5cm,framesep=0pt,nodesep=0pt,rowsep=1cm}
\begin{psmatrix}
\raisebox{1ex}{\hspace{-5ex}$
\pstree[treemode=U,nodesep=2pt,levelsep=20pt]
{\TR[name=ast]{\ast}}
  {\TR[name=t]{T}}
\hspace{4ex}
\pstree[treemode=U,nodesep=2pt,levelsep=20pt]{\TR[name=ast1]{\ast}}
  {\pstree{\TR[name=ast2]{\ast}}
     {\TR[name=tt]{T}
     }
  }
$}
\nccurve[nodesepA=2pt,nodesepB=1pt,angleA=-60,angleB=-120]{ast}{ast1}
\\
\rput{270}(-1,.2){$\mapsto$}
\\
\raisebox{4ex}{\hspace{-5ex}$
\pstree[treemode=U,nodesep=2pt,levelsep=20pt]
{\TR[name=ast]{\ast}}
  {\TR[name=t]{T}}
\hspace{4ex}
\pstree[treemode=U,nodesep=2pt,levelsep=20pt]{\TR[name=tt]{T}}{}
$}
\nccurve[nodesepA=2pt,nodesepB=1pt,angleA=-60,angleB=-120]{ast}{tt}
\end{psmatrix}
\end{center}
\end{itemize}
\begin{theorem}\label{cut-elim}
Cut elimination is well-defined: eliminating a cut from a linking on a
cut sequent yields a linking on a cut sequent.
\end{theorem}
\begin{proof}
The atomic and dual cases are trivial, since switchings correspond
before and after the elimination.  The tensor case is a simple
combinatorial argument: any cycle in a switching after elimination
induces a cycle in a switching before elimination.  (This is a
standard combinatorial argument for tensor elimination in
multiplicative proof nets in linear logic.  See
\cite{Gir87,DR89,Gir96} for details (\cf\ \cite{Hug12}, Theorem~2)).
\end{proof}
\begin{proposition}\label{confluence}
Cut elimination is locally confluent.
\end{proposition}
\begin{proof}
The only non-trivial case is a pair of atomic eliminations.  
This case
is clear from the following schematic involving two interacting atomic
cut redexes $\cutpair{\alpha}$ and $\cutpair{\beta}$.
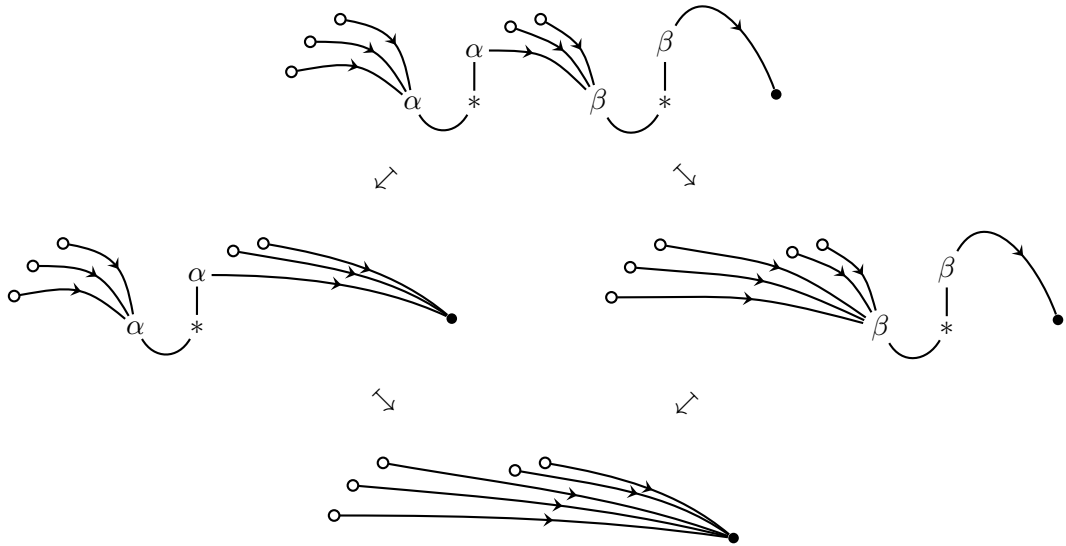
\begin{center}\psset{nodesep=-.2pt}
\newcommand{\abc}{\Cnode(-3,.5){a}
\Cnode(-2.5,.9){b}
\Cnode(-1.7,1.2){c}}
\newcommand{\bbcc}{\Cnode(-2,1.1){bb}
\Cnode(-1.2,1.2){cc}}
\newcommand{\abcA}{{\psset{nodesepB=2pt}\vecangles{a}{A}{10}{140}
\vecangles{b}{A}{0}{120}
\vecangles{c}{A}{-10}{100}}}
\newcommand{\bbccB}{{\psset{nodesepB=2pt}\vecangles{bb}{B}{-20}{120}
\vecangles{cc}{B}{-30}{105}}}
\newcommand{\maingap}{\hspace{8ex}}
\newcommand{\cutgap}{\hspace{3ex}}
\newcommand{\leftcut}{%\Cnode*{S}
$\circlenode[linestyle=none,framesep=-2pt]{A}{\alpha}$
\cutgap
\rnode{ast}{$\ast$}
\nput[labelsep=3ex]{90}{ast}{$\rnode{AA}{\alpha}$}\ncline[nodesepA=2pt,nodesepB=2pt]{ast}{AA}
\nccurve[ncurv=1,nodesep=2pt,angleA=-60,angleB=-120]{A}{ast}}%\naput{$c$}}
\newcommand{\rightcut}{%\Cnode*{B}
$\circlenode[linestyle=none,framesep=-2pt]{B}{\beta}$
\cutgap
\rnode{astC}{$\ast$}
\nput[labelsep=3ex]{90}{astC}{$\rnode{BB}{\beta}$}\ncline[nodesepA=2pt,nodesepB=2pt]{astC}{BB}
\nccurve[ncurv=1,nodesep=2pt,angleA=-60,angleB=-120]{B}{astC}}%\naput{$d$}}
\newcommand{\midjump}{{\psset{nodesep=2pt}\vecanglespos{AA}{B}{0}{135}{.4}}}
\newcommand{\abcB}{{\psset{nodesepB=2pt}\vecangles{a}{B}{0}{165}
\vecangles{b}{B}{-5}{155}
\vecangles{c}{B}{-10}{145}}}
\renewcommand{\C}{\maingap\Cnode*(0,.2){C}}
\newcommand{\rightjump}{{\psset{nodesepA=2pt}\vecangles{BB}{C}{60}{110}}}
\newcommand{\bbccC}{\vecanglesheight{bb}{C}{-10}{145}{.7}
\vecanglesheight{cc}{C}{-10}{140}{.7}}
\newcommand{\AAbbccC}{{\psset{nodesepA=2pt}\vecanglesheight{AA}{C}{0}{155}{.65}}\bbccC}
\newcommand{\abcC}{\vecanglesheight{a}{C}{0}{172}{.7}
\vecanglesheight{b}{C}{-5}{166}{.7}
\vecanglesheight{c}{C}{-10}{160}{.7}}
\newcommand{\finaledges}{\bbccC\abcC}
\begin{pspicture}(-3,.4)(0,7.4)
\rput(0,6){%
\abc
\leftcut
\abcA
\maingap
\bbcc
\rightcut
\midjump
\bbccB
\C
\rightjump
}
\rput(-8,3){%
\abc
\leftcut
\abcA
\maingap
\bbcc
\cutgap
\C
\AAbbccC
}
\rput(8,3){%
\abc
\cutgap
\maingap
\bbcc
\rightcut
\abcB
\bbccB
\C
\rightjump
}
\rput(0,0){%
\abc
\cutgap
\maingap
\bbcc
\cutgap
\C
\finaledges
}
\rput{-135}(-5.5,5){$\mapsto$}
\rput{-45}(2.5,5){$\mapsto$}
\rput{-135}(2.5,2){$\mapsto$}
\rput{-45}(-5.5,2){$\mapsto$}
\end{pspicture}\vspace{1ex}
\end{center}
(This is adapted from the proof of Proposition~1 of \cite{Hug12}.)
\end{proof}
\begin{theorem}\label{strong-norm}
Cut elimination is strongly normalising.
\end{theorem}
\begin{proof}
It is locally confluent, and eliminating a cut reduces the number of
vertices of the cut sequent.
\end{proof}
\parag{Turbo cut elimination}
Analogous to \cite{Hug12}, normalisation can be completed in a
single step.  For $l$ the $i\nth$ leaf of a shape $S$ in a cut pair
$\cutpair{S}$, let $l\perp$ denote the $i\nth$ leaf of $S\perp$.
The \defn{normal form} of a cut sequent $\Gamma$ is the sequent
$\reduc{\Gamma}$ obtained by deleting all cut pairs. Given a linking
$f$ on $\Gamma$, its \defn{normal form} $\reduc{f}$ is the linking on
$\reduc{\Gamma}$ obtained by replacing every set of directed edges
$\langle l_0,l_1\rangle,\langle l_1\perp,l_2\rangle,\langle
l_2\perp,l_3\rangle,\ldots,
\langle l_{n-1}\perp,l_n\rangle$ in $f$ in which only $l_0$
and $l_n$ occur in $\reduc{\Gamma}$ by the single directed edge
$\langle l_0,l_n\rangle$.  By a simple induction on the number of
vertices in cut sequents, $\reduc{f}$ is precisely the normal form of
$f$ under one-step cut elimination.  (In particular, this implies the
turbo normal form $\reduc{f}$ is indeed a linking, \ie, turbo cut
elimination is well defined on linkings.)

\subsection{Arbitrary base category}\label{arb}
So far we have only presented one-sided linkings over a set $\A$ of
generators.  When $\A$ is a category, as in the two-sided case we add
labels to every edge from a generator (object of $\A$), and define an
\defn{$\A$-linking} on a cut sequent $\Gamma$ as a leaf function $f$ on $\Gamma$
satisfying:
\begin{itemize}\label{onesided-arb-matching}
\myitem{1a} \Bijection. Restricting $f$ to generator-labelled leaves of
$\Gamma$ yields a bijection.
\myitem{1b} \Labelling. If $x$ is the label of an edge 
from a leaf labelled $a$ to a leaf labelled $b$, then $x:a\to b$ is a
morphism in $\A$.
\myitem{2} \Switching. For any switching $\Gamma'$ of $\Gamma$, the undirected graph obtained by adding the 
edges of $f$ to $\Gamma'$ is a tree.
\end{itemize}
For example, if $x:a\to b$
and $y:c\to c$ are morphisms in $\A$, then
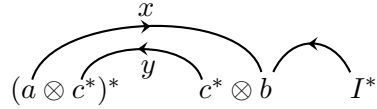
\begin{center}\begin{math}
(\rnode{a}{a}\tensor\rnode{cc}{c}\perp)\perp
\hspace{6ex}\rule{0pt}{6ex}
\rnode{c}{c}\perp\tensor\rnode{b}{b}
\hspace{6ex}\one{u}\perp
\uloopvecheight{a}{b}{.6}\naput[npos=.5]{x}
\uloopvecleftheight{c}{cc}{.7}\naput[npos=.45]{y}
\nccurve[ncurv=1.1,angleA=115,angleB=65]{u}{b}
\lput{:0}{\psline{->}(.05,.005)(.1,.005)}
\end{math}
\end{center}
is a (one-sided) $\A$-linking on the three-shape sequent $(a\tensor c\perp)\perp,\,c\perp\tensor b,\,\unit\perp$.
As in the two-sided case, conditions (1a) and (1b) reduce to (1)
\Matching{} (page~\pageref{onesided-linking-criterion}) in the
discrete case, since every label is an identity.

Atomic cut elimination (page~\pageref{atom-cut}) incorporates labels
as follows: when re-setting an $f$-edge $n\to l^+$ to target $f(l^-)$
instead, let the $\A$-morphism $x$ be the label of the edge $n\to l^+$
(if any) and let the $\A$-morphism $y$ be the label of the edge
$l^-\to f(l^-)$ (if any); if both $x$ and $y$ are present, label the
output edge $n\to f(l^-)$ by the composite $\seqcomp y x$ in $\A$,
otherwise leave $n\to f(l^-)$ unlabelled.  Correspondingly, we adjust
the definition of turbo cut elimination: the output edge $\langle
l_0,l_n\rangle$ is labelled by the composite $\A$-morphism $x_1;\ldots;
x_n$ iff every $\langle l_{i-1},l_i\rangle$ is labelled by an
$\A$-morphism $x_i$ (\cf\ path composition of two-sided linkings
defined in Section~\ref{arbitrary}).

The properties of cut elimination (Theorem~\ref{cut-elim},
Proposition~\ref{confluence} and Theorem~\ref{strong-norm}) are
unaffected by the presence of labels.  

\subsection{From one-sided linkings to two-sided linkings}\label{proofs}

By design, two-sided linkings $S\to T$ between shapes, as defined in
Section~\ref{linkings}, are in bijection with one-sided linkings on
the two-shape sequent $S\perp,T$.  Via this correspondence, we can
take care of the proof obligations remaining from
Section~\ref{linkings}.

\begin{proofof}{Proposition~\ref{linear-time} (linear time verification of the linking criterion)}
Every star-autonomous shape $S$ induces a formula $\mll{S}$ of multiplicative linear logic:
replace negative tensors in $S$ by pars $\parr$, replace negative
generators $a$ by $a^\bot$, replace negative $\unit$'s by $\bot$, and
delete all duals $(-)\perp$.  For example, $\big((1\tensor
a\perp)\perp\tensor 1\big)\perp$ becomes $(1\tensor a^\bot)\parr\bot$.
Thus every shape sequent $\Gamma=S_1,\ldots,S_n$ induces a formula
sequent $\mll{\Gamma}$, namely $\mll{S_1},\ldots,\mll{S_n}$.
The formula $\mll{S}$ has the same leaf vertices as the original shape $S$,
and switchings of $\mll{S}$ are in bijection with switchings of $S$.
Thus a leaf function on a shape sequent $\Gamma$ is a linking iff it
constitutes a multiplicative proof net on the formula sequent
$\mll{\Gamma}$, in the sense of \cite{Hug12}.
Via the translation $\Gamma\mapsto\mll{\Gamma}$,
Proposition~\ref{linear-time} becomes a corollary of Theorem~4 of
\cite{Hug12} (which, in turn, is not much more than a corollary of
linear time verification of the proof net criterion for unit-free
multiplicative nets \cite{Gue99,MO00}).
\end{proofof}
Given linkings $f:S\to T$ and $g:T\to U$, since the definition of
turbo cut elimination is precisely path composition, the composite
$\seqcomp f g:S\to U$ corresponds to the normal form of the one-sided linking
$f\cup g$, the disjoint union of $f$ and $g$, on the cut sequent
$S\perp,\cutpair{T},U$.  
\begin{proofof}{Proposition~\ref{link-comp} (two-sided linkings compose)}
Via the correspondence just described, this is a corollary of
Theorem~\ref{cut-elim} (and the correspondence between turbo cut
elimination of one-sided linkings and normalisation by one-step cut
elimination).
\end{proofof}
\begin{proofof}{Theorem~\ref{compatibility} (compatibility: cycles do 
not arise during composition)} Suppose $f:S\to T$ and $g:T\to U$ are
linkings such that the disjoint union $f+g$, a directed graph on
$\leaves{S}+\leaves{T}+\leaves{U}$, contains a cycle.  Let $f\cup g$
be the corresponding one-sided linking on the cut sequent
$S\perp,\cutpair{T},U$, and let $\overline{f\cup g}$ be the result of
eliminating all tensor $\tensor$ and dual $(-)\perp$ cuts from $f\cup
g$, a one-sided linking on the cut sequent
\raisebox{.3ex}{\small$S,\:\cutpair{\alpha_1\!},\ldots,\:\cutpair{\alpha_n\!},U$},
where $\alpha_i$ is the atom labelling the $i\nth$ leaf of $T$.  Had
the directed graph $f+g$ a cycle $C$, then $\overline{f\cup g}$ would
contain a cycle in every switching (the edges of the cycle $C$
alternating with cut edges between the $\alpha_i$), and therefore fail
to be a linking, contradicting Theorem~\ref{cut-elim}.
\end{proofof}
See \cite{Blu93} for more on the relationship between compatibility
and the multiplicative proof net criterion, in the unit-free case.

\section{One-sided nets}\label{one-sided-nets}\label{net-comp-proof}

We define a one-sided net as a one-sided linking modulo Trimble
rewiring \cite{Tri94,BCST96}.

One-sided $\A$-linkings $f$ and $g$ on a cut sequent $\Gamma$ are
\defn{similar} if they differ on a negative $\unit$, \ie, one can be
obtained from the other by re-targeting one edge from an $\unit$.  A
one-sided \defn{$\A$-net} on a cut sequent $\Gamma$ is an equivalence
class of $\A$-linkings on $\Gamma$ modulo similarity (\ie, modulo the
transitive closure of similarity).  Similarity here coincides with the
earlier two-sided case: two-sided $\A$-linkings $f,g:S\to T$ are
similar iff the corresponding one-sided linkings on the two-shape
sequent $\;S\perp,T\;$ are similar in the sense just defined.
\begin{theorem}\label{net-cut-elim}
Cut elimination respects equivalence, \ie, cut elimination is
well-defined on $\A$-nets.  More precisely: if $f$ and $g$ are
equivalent $\A$-linkings on the cut sequent $\Gamma$ containing a cut
pair $\cutpair{S}$, and the $\A$-linkings $f'$ and $g'$ result from
eliminating $\cutpair{S}$ from $f$ and $g$, respectively, then $f'$
and $g'$ are equivalent.
\end{theorem}
Before proving theorem (below), we illustrate it with an example.
(The notation in the example corresponds to the notation in the
proof.)  Similar (hence equivalent) linkings $f$ and $g$ are shown top
and bottom, with adjacent normal forms $f'$ and $g'$; the sequence
$h_0,h_1,h_2,h_3$ of pairwise similar linkings witnesses the equivalence of $f'$ and $g'$.
%$h_0$, $\,h_1$, $\,h_2$, $\,h_3$ provide a sequence of pairwise similar linkings
%witnessing the equivalence of $f'$ and $g'$.
\vspace*{5ex}\begin{displaymath}
\newcommand{\gap}{\hspace{4ex}}
\newcommand{\nl}{\\[9ex]}
\renewcommand{\u}[1]{\circlenode[framesep=0ex,linestyle=none]{#1}{\unit}}
\renewcommand{\perp}{\makebox[0pt]{$\rule{0ex}{0ex}^*$}}
\newcommand{\cutverts}{
  \u{i}
  \gap
  \u{n}\perp
  \lineanglesheight{i}{n}{-45}{-135}{.8}
}
\newlength{\cutwidth}
\settowidth{\cutwidth}{$\cutverts$}
\newcommand{\cutvertsgone}{\makebox[\cutwidth]{}}
\newcommand{\leftverts}{
  \u{l1}\perp
  \gap
  \u{l2}\perp
  \gap
  \u{l3}\perp
  \gap
}
\newcommand{\rightverts}{  
  \gap
  \u{a}
  \tensor
  \u{b}
}
\newcommand{\verts}[2]{\makebox[0pt][r]{$#1$}:\hspace{8ex} \leftverts #2 \rightverts}
\newcommand{\allverts}[1]{\verts{#1}{\cutverts}
  \vecanglesheight{l1}{i}{55}{120}{0.5}
  \vecanglesheight{l2}{i}{50}{130}{0.5}
  \vecanglesheight{l3}{i}{45}{140}{0.5}}
\newcommand{\normverts}[1]{\verts{#1}{\cutvertsgone}}
\newcommand{\edges}[3]{
  \vecanglesheight{l1}{#1}{55}{120}{0.3}
  \vecanglesheight{l2}{#2}{50}{130}{0.3}
  \vecanglesheight{l3}{#3}{45}{140}{0.3}
}
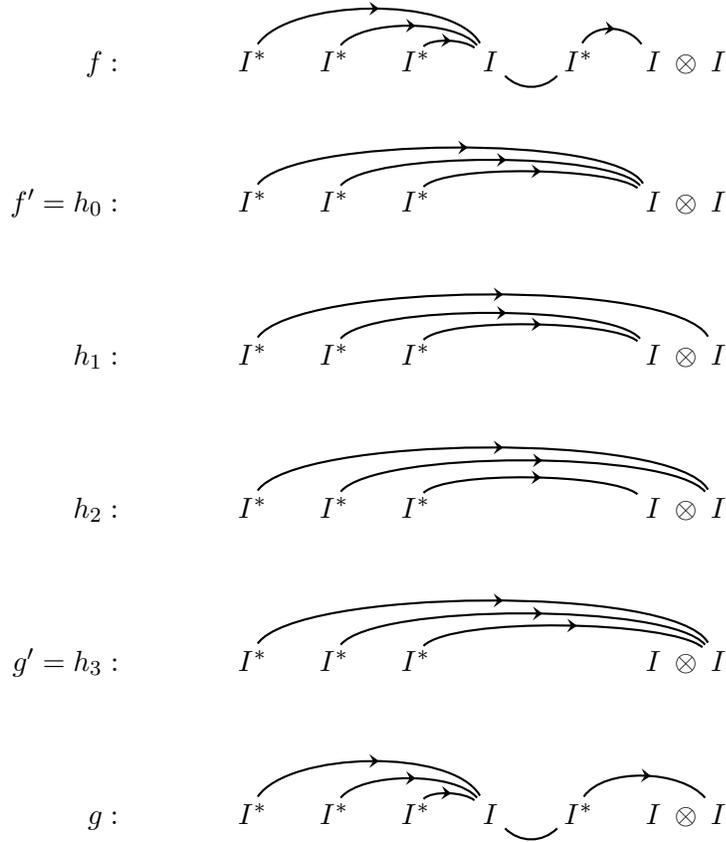
\begin{array}{c}
  \allverts{f} \vecanglesheight{n}{a}{60}{125}{.8}
  \\[8ex]
  \normverts{f'=h_0}\edges a a a
  \nl
  \normverts{h_1}\edges b a a
  \nl
  \normverts{h_2}\edges b b a
  \nl
  \normverts{g'=h_3}\edges b b b
  \nl
  \allverts{g}\vecanglesheight{n}{b}{55}{130}{.6}
  \\[3ex]
\end{array}\end{displaymath}
\begin{proof}
If $\cutpair{S}$ is a tensor or dual cut, the result is trivial, since
leaves and linkings are untouched by eliminating the cut.  Suppose the
cut is atomic.  Let $l^-$ be the negative leaf of 
$\cutpair{S}$.  If $l^-$ is not the leaf $l$ on which $f$ and $g$
differ, then $f'$ and $g'$ are similar or equal (hence equivalent, as
desired) since they differ on at most $l$ after eliminating the cut.
So assume $l^-=l$.  
Thus we have the situation
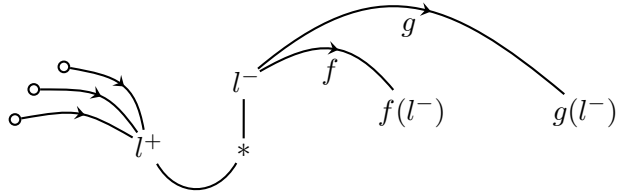
\begin{center}
\hspace{8ex}\psset{colsep=2.5cm,framesep=0pt,nodesep=0pt,rowsep=1.2cm}
\Cnode(-3,.5){a}
\Cnode(-2.5,.9){b}
\Cnode(-1.7,1.2){c}
\circlenode[linestyle=none]{A}{$l^+\!$}\rule{0pt}{12ex}
\hspace{4ex}
\circlenode[linestyle=none]{ast}{$\ast$}
\nput[labelsep=3ex]{90}{ast}{\circlenode[linestyle=none]{AA}{$l^-\!$}}\ncline{ast}{AA}
{\psset{nodesep=1pt}\nccurve[ncurv=1,angleA=-60,angleB=-120]{A}{ast}}
\psset{nodesepB=-1pt}
\vecangles{a}{A}{10}{150}
\vecangles{b}{A}{0}{125}
\vecangles{c}{A}{-10}{100}
\hspace{4ex}
\rput{0}(2,.6){\circlenode[linestyle=none]{d}{$f(l^-)$}}
\hspace{7ex}
\rput{0}(4,.6){\circlenode[linestyle=none]{dd}{$g(l^-)$}}
\psset{nodesepB=-1pt,nodesepA=0pt,nodesepB=-4pt}
\vecangles{AA}{d}{30}{130}\nbput{$f$}
\vecangles{AA}{dd}{40}{140}\nbput{$g$}
\hspace{5ex}\vspace{2ex}
\end{center}
where the leaf vertices $l^-$ and $l^+$ are labelled by $\unit$, the
unlabelled directed edges are present in both $f$ and $g$, and the
edge labelled $f$ (resp.\ $g$) is present in $f$ (resp.\ $g$) only.
The vertices $\circ$ schematically represent the negative leaves
$l_1\ldots l_n$ of $\Gamma$ whose edge targets $l^+$ (in both $f$ and
$g$), \ie, such that $f(l_i)=g(l_i)=l^+$.  Note that, since $l^+$,
$l^-$ and the $l_i$ are labelled $\unit$, none of the edges from them
(\ie, the directed edges depicted above) is labelled by an
$\A$-morphism, so the case of non-discrete $\A$ coincides with the
case of $\A$ a set.

Let $h$ be the (partial) leaf function obtained from $f$ (or
equivalently $g$) by deleting the edge from $l^-$.  Since $f$ and $g$
are linkings, every switching $\sigma$ of $h$ is a disjoint union of
two trees, with all $l_i$ in one tree and $f(l^-)$ and $g(l^-)$ in the
other (otherwise the corresponding switching of one of $f$ or $g$
would contain a cycle, via the edge between $l^-$ and $f(l^-)$ or
$g(l^-)$, respectively).  Let $h^-$ be the result of deleting from $h$
all the edges from the $l_i$, and for $0\le j\le n$ construct $h_j$
from $h^-$ by adding an edge from $l_i$ to $g(l^-)$ for $1\le i\le j$,
and from $l_i$ to $f(l^-)$ for $j<i\le n$, and also an edge from $l^-$
to $f(l^-)$ (or $g(l^-)$; the choice is arbitrary).
Since for any switching of $h$, the $l_i$ are in one tree and $f(l^-)$
and $g(l^-)$ are in the other, $h_i$ is a well-defined linking.
Let $h_i'$ be result of eliminating the atomic cut from $h_i$. By
design, $f'=h_0'$ and $g'=h_n'$.  By Theorem~\ref{cut-elim}, each
$h'_i$ is a linking.  Since $h_i'$ and $h_{i-1}'$ differ on just one
$\unit$ (namely $l_i$), they are similar, so $f'$ and $g'$ are
equivalent, via the $h'_i$.
\end{proof}
By the correspondence between (turbo) cut elimination of one-sided
linkings and path composition of two-sided linkings, we obtain
Theorem~\ref{net-comp} (compositionality of two-sided nets).

\section{Proof that the category $\netsa$ of $\A$-nets is free star-autonomous}\label{mainproof}

This section proves the Freeness Theorem, Theorem~\ref{free} on page~\pageref{free}:
\textsl{For any category $\A$, the category $\netsa$ of $\A$-nets 
is the free star-autonomous category generated by $\A$}.

\subsection{Lax linkings}

Suppose $\A$ is a set.
Given $\A$-shapes $S$ and $T$, define a \defn{lax leaf function} $S\to
T$ as the relaxation of a leaf function (as defined in
Section~\ref{one-sided}) obtained by permitting edges from $\unit$'s to
target any vertex of $S$ or $T$ (viewed as a parse trees).  For
example, here is a lax leaf function from $(a\perp\tensor \unit)\tensor \unit$
to $(a\perp \tensor I)\tensor (b\tensor b\perp)\perp$, drawn in
parse-tree form on the left, and compact in-line form on the right:
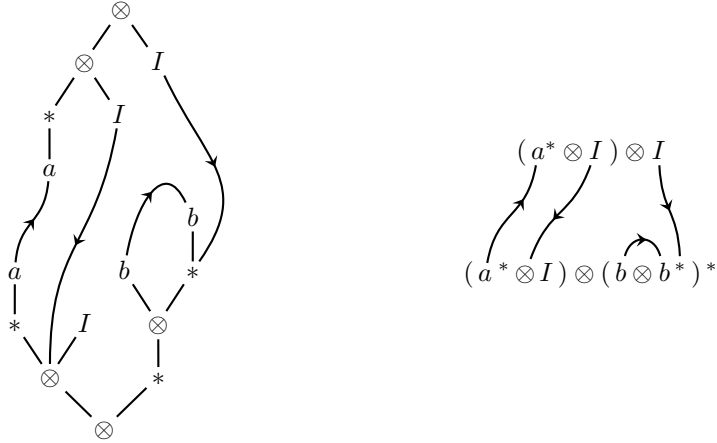
\begin{center}\small\begin{math}\begin{array}{c}
\psset{colsep=3cm,framesep=0pt,nodesep=0pt,rowsep=1cm}
\psset{treemode=D,nodesep=2pt,levelsep=20pt}
  \pstree{\TR[name=t11]{\tensor}}
    {
     \pstree{\TR[name=t112]{\tensor}}
      {
       \pstree{\TR[name=ast1]{\ast}}
         {\TR[name=aa]{a}
         }
       \TR[name=a]{\unit}
      }
      {
       \TR[name=uu]{\unit} 
      }
    }
  
\\\\
\psset{treemode=U,nodesep=2pt,levelsep=20pt}
\pstree{\TR[name=tt1]{\tensor}}
  {\pstree{\TR[name=tt11]{\tensor}}
    {
     \pstree{\TR[name=ast'1]{\ast}}
       {\TR[name=aa1]{a}
       }
     \TR[name=a1]{\unit}
    }
     {\pstree{\TR[name=ast'2]{\ast}}
       {\pstree{\TR[name=tt121]{\tensor}}
          {\TR[name=bb1]{b}
           \pstree{\TR[name=ast'3]{\ast}}
            {\TR[name=b1]{b}
            }
           
          }
       }
     }
  }
\vecanglespos{uu}{ast'3}{-65}{55}{.5}
\vecanglespos{a}{tt11}{-100}{90}{.5}
\vecanglesheight{aa1}{aa}{85}{-95}{1}
{\psset{nodesepA=2pt}\vecanglesposheight{bb1}{b1}{80}{115}{.6}{1.3}}
\end{array}
\hspace{20ex}
\begin{array}{c}
\begin{psmatrix}[rowsep=3.3\baselineskip]
(\,\rnode{a}{a}\perp\tensor\rnode{aa}{I}\:)\tensor \one{uu} \\[8ex]
(\,\rnode{a1}{a}\:\perp\,\rnode{aa1}{\tensor}\,{I}\,) \tensor (\,\rnode{b1}{b}\tensor\rnode{bb1}{b}\,\rnode{perp}{{}\perp}\,)\,\perp
\vecanglespos{uu}{perp}{-87}{93}{.5}
\vecanglespos{aa}{aa1}{-110}{70}{.59}
\uloopvec{b1}{bb1}
\vecanglespos{a1}{a}{80}{-100}{.67}
\end{psmatrix}\end{array}
\end{math}\end{center}
Define the graph of a lax leaf function $f:S\to T$ by analogy with the
original non-lax case: the undirected graph which is the disjoint
union of the two parse trees $S$ and $T$, together with the edges of
$f$, undirected.  Define a switching as before.  A \defn{lax linking}
is a lax leaf function satisfying the linking criterion, \ie,
conditions (1) \Matching{} and (2) \Switching{} on
page~\pageref{criterion}.  For example, the lax leaf function depicted
below-left is a lax linking, since both its switchings are trees, but
the lax leaf function below-right is not, since the switching in which
the upper tensor chooses its right argument has a cycle.
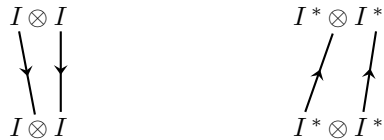
\begin{center}\vspace{1ex}
\small\psset{rowsep=2.5\baselineskip}$\begin{psmatrix}
\one{b}\,\rnode{t}{\tensor}\:\one{bb} \\ 
\one{b1}\,\rnode{t1}{\tensor}\:\one{bb1}
\linevec{b}{t1}\linevec{bb}{bb1}
\end{psmatrix}
\rule{3cm}{0cm}
\begin{psmatrix} 
\one{b}\,\perp\,\rnode{t}{\tensor}\:\one{bb}\,\rnode{perp}{{}\perp} \\
\one{b1}\,\perp\,\rnode{t1}{\tensor}\:\one{bb1}\,\rnode{perp1}{{}\perp}\;
\linevec{b1}{t}\linevec{bb1}{perp}
\end{psmatrix}$\vspace{3ex}\end{center}

We make the corresponding lax definitions in the one-sided case.
Given a cut sequent $\Gamma$ over $\A$, define a \defn{lax leaf
function} on $\Gamma$ as the variant of a leaf function obtained by
permitted edges from $\unit$'s to target any vertex of $\Gamma$.
Switchings again generalise to lax leaf functions in the obvious way.
Define a \defn{lax linking} on $\Gamma$ as a lax leaf function on
$\Gamma$ which satisfies the linking criterion, \ie, conditions (1)
\Matching{} and (2) \Switching{} on
page~\pageref{onesided-linking-criterion}.
For example, here is a lax linking on the three-shape sequent
$(I\tensor a)\perp,\,a\tensor I,\,I\perp$\,:
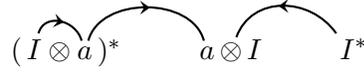
\begin{center}\vspace{3ex}\begin{math}
(\,\rnode{aa}{I}\tensor\rnode{bb}{a}\,)\perp
\hspace{6ex}\rule{0pt}{6ex}
\rnode{b}{a}\,\rnode{t}{\tensor}\,\rnode{a}{I}
\hspace{6ex}\one{u}\perp
\uloopvecheight{aa}{bb}{1}
\uloopvecheight{bb}{b}{.8}
\nccurve[ncurv=.9,angleA=125,angleB=70,nodesepB=2pt]{u}{t}
\lput{:0}{\psline{->}(.05,.005)(.1,.005)}
\end{math}
\vspace{3ex}\end{center}

When $\A$ is an arbitrary category, generalise the definitions of lax
linking exactly as in the original non-lax case: add $\A$-morphisms as
labels on edges between generators (objects of $\A$), then replace
condition (1) \Matching{} by conditions (1a) \Bijection{} and (1b)
\Labelling{} (page~\pageref{arb-matching} in the two-sided case, and
page~\pageref{onesided-arb-matching} in the one-sided case).
Extending the non-lax case, lax linkings $S\to T$ correspond to lax
linkings on the two-shape sequent $\;S\perp,T\:$ such that no edge
targets the distinguished (\ie\ outermost) $\ast$-vertex of $S\perp$.

\subsection{Lax equivalence}

Define two lax linkings as \defn{similar} if they differ by a single
edge from an $\unit$, and \defn{lax equivalent} if they are equivalent
modulo similarity on lax linkings (\ie, modulo the
reflexive-transitive closure of similarity).  To help avoid ambiguity,
while dealing with lax linkings in this section we shall refer to the
original non-lax notion of a linking as a \defn{standard linking}, and
the original non-lax notion of equivalence between standard linkings
as \defn{standard equivalence}.  The following Lemma, proved in
section~\ref{proof-lax-rewiring}, is the key technical step in
the proof that the category $\netsa$ is free.
\begin{lemma}[Lax rewiring]\label{lax-rewiring}
Standard linkings are standard equivalent iff they are lax equivalent.
\end{lemma}
Thus if the standard linking $f$ can be rewired to the standard
linking $g$ along a sequence $f=h_1\ldots h_n=g$ of lax linkings with
$h_i$ similar to $h_{i-1}$ for $1<i\le n$, then $f$ can be rewired to
$g$ along a sequence of \emph{standard} linkings: there exists a
sequence $f=k_1\ldots k_m=g$ of standard linkings with $k_i$ similar
to $k_{i-1}$ for $1<i\le n$.  In other words, adding lax linkings has no
impact on equivalence of standard linkings: no additional standard
linkings are identified when we permit `lax rewiring', via lax
linkings.

\subsubsection{Proof of the Lax Rewiring Lemma}\label{proof-lax-rewiring}

Since two-sided lax linkings $S\to T$ are in bijection with certain
one-sided lax linkings on the two-formula sequent $S\perp,T$,
henceforth we shall assume all lax linkings are one-sided.

An \defn{atomic linking} is a lax linking whose every edge targets an
atom.  (Note that an atomic linking need not be a standard linking
since an edge from an $\unit$ may target a negative leaf.)  Define
atomic linkings $f$ and $g$ as \defn{atomic equivalent} if they are
lax equivalent via atomic linkings: there exists a sequence
$f=h_1\ldots h_n=g$ of atomic linkings with $h_i$ similar to $h_{i-1}$
for $1<i\le n$.  The following lemma reduces atomic equivalence to standard
equivalence.
\begin{lemma}\label{lemma-reduce-atomic}
Standard linkings are standard equivalent iff they are atomic equivalent.
\end{lemma}
\begin{proof}
Define the \defn{depth} of an atomic linking as the number of its
edges which target negative leaves.  Thus an atomic linking is a
standard linking iff it has depth $0$.  Define the depth of a sequence
of lax linkings $f_1\ldots f_m$ as the maximum of the depths of the
$f_i$.  Suppose $f_1\ldots f_m$ has depth $d$.  The \defn{size} of
$f_1\ldots f_m$ is the number of $f_i$ of depth $d$.

Suppose $f$ and $g$ are standard linkings which are atomic equivalent
via a sequence $f=h_1\ldots h_n=g$ of atomic linkings with $h_i$ similar to
$h_{i-1}$ for $1<i\le n$.
We proceed by a primary induction on the depth of $h_1\ldots h_n$, and
a secondary induction on its size.

\begin{itemize}
\item Primary induction base: $h_1\ldots h_n$ has depth $0$.  
Then all $h_i$ are standard linkings, so $f$ and $g$ 
are already standard equivalent.
\item Primary induction step: $h_1\ldots h_n$ has depth $d>0$.
Let $pqr$ be a three-element subsequence of $h_1\ldots h_n$ with $q$
of depth $d$ (\ie, $q$ is of maximum depth in the sequence) and $r$ of
depth $d-1$.  (Such a subsequence exists since $h_n=g$ has depth $0$,
and consecutive elements in the sequence differ by at most one in depth.)

Let $l$ be the negative leaf such that $q(l)$ is a negative leaf
$l^-$, and $r(l)$ is a positive leaf $l^+$, and otherwise $q$ and $r$
are identical.  (The leaf $l$ must exist, since $q$ has depth $d$
while $r$ has depth $d-1$.)  
Let $l'$ be the negative leaf such that
$p(l')\neq q(l')$, and let $l''=p(l')$.  If $l\neq l'$
then $p$, $q$ and $r$ are:
\begin{equation}
\addtolength{\arraycolsep}{.4ex}
\begin{array}{r@{\;\;\;\;}cccccc}
\rule{0pt}{4ex}
p: & \rnode{l}{l} & \rnode{l-}{l}^-\hspace{-1.5ex} & & \rnode{l'}{l'} & \rnode{l''}{l''}\hspace{-1.5ex} \\[1.7ex]
\uloopvec{l}{l-}\uloopvec{l'}{l''}
q: & \rnode{l}{l} & \rnode{l-}{l}^-\hspace{-1.5ex} & & \rnode{l'}{l'} & & \rnode{l'''}{l'''} \\[1.7ex]
\uloopvec{l}{l-}\psset{nodesepB=0pt}\vecangles{l'}{l'''}{35}{145}
r: & \rnode{l}{l} & & \rnode{l+}{l}^+ & \rnode{l'}{l'} & & \rnode{l'''}{l'''} 
\vecangles{l}{l+}{35}{145}\psset{nodesepB=0pt}\vecangles{l'}{l'''}{35}{145}
\end{array}
\end{equation}
where $l'''=q(l')=r(l')$, and if $l=l'$ then $p$, $q$ and $r$ are:
\begin{equation}\label{l-eq-l'}
\addtolength{\arraycolsep}{.4ex}
\begin{array}{r@{\;\;\;\;}cccccc}
\rule{0pt}{4ex}
p: & \rnode{l}{l} & \rnode{l''}{l''}\hspace{-1ex} \\[1.7ex]
\uloopvec{l}{l''}
q: & \rnode{l}{l} & & \rnode{l-}{l}^-\hspace{-1.5ex} \\[1.7ex]
\vecangles{l}{l-}{35}{145}
r: & \rnode{l}{l} & & & \rnode{l+}{l}^+
\vecangles{l}{l+}{35}{145}
\end{array}
\end{equation}
\begin{itemize}
\myitem{a} Secondary induction base: $h_1\ldots h_n$ has size $1$.  Thus $p$ has depth $d-1$.  Therefore the 
rewirings $p\mapsto q$ and $q\mapsto r$ each re-target the same
negative leaf, \ie, $l=l'$, as in diagram~(\ref{l-eq-l'}) above, so $p$
and $r$ are similar, without need for $q$ as an intermediate.
Delete $q$ from $h_1\ldots h_n$, and appeal to the primary induction
hypothesis with this new sequence of depth $d-1$.
\myitem{b} Secondary induction step: $h_1\ldots h_n$ has size $s>1$.  If $p$ has depth $d-1$, we can delete $q$ exactly 
as in the previous case, reducing the size of $h_1\ldots h_n$ by $1$,
and appeal to the secondary induction hypothesis.  So assume $p$ has depth $d$.
Let $l^{\star+}$ be the leaf at the end of the path in (the directed graph)
$p$ which begins at $l^-$.  Thus $l^{\star+}$ is positive.
\begin{itemize}
\myitem{i} Case $l^-\neq l'$.  We have:
\begin{equation}
\addtolength{\arraycolsep}{.4ex}
\begin{array}{r@{\;\;\;\;}ccccccccc@{\;\;\;\;\;\;\;\;\;\;}l}
\rule{0pt}{4ex}
p: & \rnode{l}{l} & & \rnode{lminus}{l}^-\hspace*{-.5ex}\pnode(0,.33){e0}
& \pnode(0,.33){e1} & \pnode(0,.33){e2} \raisebox{2ex}{\small$\,\cdots$} \pnode(0,.33){e3} 
& \pnode(0,.33){e4}\rnode{lstar+}{l^\star}^+ & \rnode{l'}{l'} & \rnode{l''}{l''}
&& \text{(depth $d$)}\\[1.7ex]
\vecangles{l}{lminus}{35}{145}
\uloopvec{l'}{l''}
\psset{nodesep=0pt}
\uloopvec{e0}{e1}
\uloopvec{e1}{e2}
\uloopvec{e3}{e4}
q: & \rnode{l}{l} & & \rnode{lminus}{l}^-\hspace*{-.5ex}\pnode(0,.33){e0}
& \pnode(0,.33){e1} & \pnode(0,.33){e2} \raisebox{2ex}{\small$\,\cdots$} \pnode(0,.33){e3} 
& \pnode(0,.33){e4}\rnode{lstar+}{l^\star}^+ & \rnode{l'}{l'} && \rnode{l'''}{l'''}
& \text{(depth $d$)}\\[1.7ex]
\vecangles{l'}{l'''}{35}{145}
\vecangles{l}{lminus}{35}{145}\psset{nodesepA=5pt,nodesepB=0pt}
\psset{nodesep=0pt}
\uloopvec{e0}{e1}
\uloopvec{e1}{e2}
\uloopvec{e3}{e4}
r: & \rnode{l}{l} & \rnode{l+}{l}^+ &
&&&
& \rnode{l'}{l'} && \rnode{l'''}{l'''}
\vecangles{l'}{l'''}{35}{145}
\uloopvec{l}{l+}\psset{nodesepA=5pt,nodesepB=0pt}
& \text{(depth $d-1$)}
\end{array}
\end{equation}
Define $p'$ as $p$ but for $p'(l)=l^{\star+}$ (versus $p(l)=l^-$) and
similarly, define $q'$ as $q$ but for $q'(l)=l^{\star+}$ (versus
$q(l)=l^-$). Substitute $p'q'$ for $q$:
\begin{equation}
\addtolength{\arraycolsep}{.4ex}
\begin{array}{r@{\;\;\;\;}ccccccccc@{\;\;\;\;\;\;\;\;\;\;}l}
\rule{0pt}{4ex}
p: & \rnode{l}{l} & & \rnode{lminus}{l}^-\hspace*{-.5ex}\pnode(0,.33){e0}
& \pnode(0,.33){e1} & \pnode(0,.33){e2} \raisebox{2ex}{\small$\,\cdots$} \pnode(0,.33){e3} 
& \pnode(0,.33){e4}\rnode{lstar+}{l^\star}^+ & \rnode{l'}{l'} & \rnode{l''}{l''}
&& \text{(depth $d$)}\\[1.7ex]
\vecangles{l}{lminus}{35}{145}
\uloopvec{l'}{l''}
\psset{nodesep=0pt}
\uloopvec{e0}{e1}
\uloopvec{e1}{e2}
\uloopvec{e3}{e4}
p': & \rnode{l}{l} & & \rnode{lminus}{l}^-\hspace*{-.5ex}\pnode(0,.33){e0}
& \pnode(0,.33){e1} & \pnode(0,.33){e2} \raisebox{2ex}{\small$\,\cdots$} \pnode(0,.33){e3} 
& \pnode(0,.33){e4}\rnode{lstar+}{l^\star}^+ & \rnode{l'}{l'} & \rnode{l''}{l''}
&& \text{(depth $d-1$)}\\[1.7ex]
\vecanglesheight{l}{lstar+}{35}{95}{.5}
\uloopvec{l'}{l''}
\psset{nodesep=0pt}
\uloopvec{e0}{e1}
\uloopvec{e1}{e2}
\uloopvec{e3}{e4}
q': & \rnode{l}{l} & & \rnode{lminus}{l}^-\hspace*{-.5ex}\pnode(0,.33){e0}
& \pnode(0,.33){e1} & \pnode(0,.33){e2} \raisebox{2ex}{\small$\,\cdots$} \pnode(0,.33){e3} 
& \pnode(0,.33){e4}\rnode{lstar+}{l^\star}^+ & \rnode{l'}{l'} && \rnode{l'''}{l'''}
& \text{(depth $d-1$)}\\[1.7ex]
\vecangles{l'}{l'''}{35}{145}
\vecanglesheight{l}{lstar+}{35}{95}{.5}
\psset{nodesepA=5pt,nodesepB=0pt}
\psset{nodesep=0pt}
\uloopvec{e0}{e1}
\uloopvec{e1}{e2}
\uloopvec{e3}{e4}
r: & \rnode{l}{l} & \rnode{l+}{l}^+ &
&&&
& \rnode{l'}{l'} && \rnode{l'''}{l'''}
\vecangles{l'}{l'''}{35}{145}
\uloopvec{l}{l+}\psset{nodesepA=5pt,nodesepB=0pt}
& \text{(depth $d-1$)}
\end{array}
\end{equation}
Note that $p'$ has depth $d-1$ since $p$ has depth $d$ and $l^{\star+}$ is
positive whereas $l^-$ is negative; similarly, $q'$ has depth $d-1$.

The lax leaf functions $p'$ and $q'$ are lax linkings (\ie, satisfy the
linking correctness criterion),
since $p$ and $q$ are lax linkings: the edge from $l$ targets $l^-$ in $p$
and $q$, and $l^{\star+}$ in $p'$ and $q'$; since $l$ and $l^{\star+}$
are connected along the edges of the lax leaf function, a switching of $p$
(resp.\ $q$) is a tree iff the corresponding switching of $p'$ (resp.\
$q'$) is a tree.

By construction, the pairs $p\leftrightarrow p'$, $p'\leftrightarrow
q'$ and $q'\leftrightarrow r$ are similar.  Thus we can appeal to the
inductive hypothesis with the original sequence $h_1\ldots h_n$ with
$p'q'$ substituted for $q$, which has strictly smaller size than the
original (since $p'$ and $q'$ each have depth $d-1$, whereas $q$ has
depth $d$).
\myitem{ii} Case: $l^-=l'$.  We have:
\begin{equation}
\addtolength{\arraycolsep}{.4ex}
\begin{array}{r@{\;\;\;\;}cccccccc@{\;\;\;\;\;\;\;\;\;\;}l}
\rule{0pt}{4ex}
p: & \rnode{l}{l} & & \rnode{lminus}{l}^- & & \rnode{l''}{l''}\hspace*{-.5ex}\pnode(0,.33){e0}
& \pnode(0,.33){e1} & \pnode(0,.33){e2} \raisebox{2ex}{\small$\,\cdots$} \pnode(0,.33){e3} 
& \pnode(0,.33){e4}\rnode{lstar+}{l^\star}^+
& \text{(depth $d$)}\\[1.7ex]
\vecangles{l}{lminus}{35}{145}\psset{nodesepA=7pt}\vecangles{lminus}{l''}{25}{155}
\psset{nodesep=0pt}
\uloopvec{e0}{e1}
\uloopvec{e1}{e2}
\uloopvec{e3}{e4}
q: & \rnode{l}{l} & & \rnode{lminus}{l}^- & \rnode{l'''}{l'''}
&&&
&& \text{(depth $d$)}\\[1.7ex]
\vecangles{l}{lminus}{35}{145}\psset{nodesepA=5pt,nodesepB=0pt}
\vecangles{lminus}{l'''}{40}{140}
r: & \rnode{l}{l} & \rnode{l+}{l}^+  & \rnode{lminus}{l}^- & \rnode{l'''}{l'''} 
\uloopvec{l}{l+}\psset{nodesepA=5pt,nodesepB=0pt}
\vecangles{lminus}{l'''}{40}{140}
&&&
&& \text{(depth $d-1$)}
\end{array}
\end{equation}
(Note that $l''=l^{\star+}$ is possible.)
Define $p'$ as $p$ but for $p'(l)=l^{\star+}$ (versus $p(l)=l^-$) and define $r'$ as $r$ but for $r'(l)=l^{\star+}$ (versus $r(l)=l^+$):
\begin{equation}
\addtolength{\arraycolsep}{.4ex}
\begin{array}{r@{\;\;\;\;}cccccccc@{\;\;\;\;\;\;\;\;\;\;}l}
\rule{0pt}{4ex}
p: & \rnode{l}{l} & & \rnode{lminus}{l}^- & & \rnode{l''}{l''}\hspace*{-.5ex}\pnode(0,.33){e0}
& \pnode(0,.33){e1} & \pnode(0,.33){e2} \raisebox{2ex}{\small$\,\cdots$} \pnode(0,.33){e3} 
& \pnode(0,.33){e4}\rnode{lstar+}{l^\star}^+
& \text{(depth $d$)}\\[3ex]
\vecangles{l}{lminus}{35}{145}\psset{nodesepA=7pt}\vecangles{lminus}{l''}{25}{155}
\psset{nodesep=0pt}
\uloopvec{e0}{e1}
\uloopvec{e1}{e2}
\uloopvec{e3}{e4}
p': & \rnode{l}{l} & & \rnode{lminus}{l}^- & & \rnode{l''}{l''}\hspace*{-.5ex}\pnode(0,.33){e0}
& \pnode(0,.33){e1} & \pnode(0,.33){e2} \raisebox{2ex}{\small$\,\cdots$} \pnode(0,.33){e3} 
& \pnode(0,.33){e4}\rnode{lstar+}{l^\star}^+
& \text{(depth $d-1$)}\\[3ex]
\vecanglesheight{l}{lstar+}{25}{95}{.4}\psset{nodesepA=7pt}\vecangles{lminus}{l''}{25}{160}
\psset{nodesep=0pt}
\uloopvec{e0}{e1}
\uloopvec{e1}{e2}
\uloopvec{e3}{e4}
r': & \rnode{l}{l} & \rnode{l+}{l}^+  & \rnode{lminus}{l}^- & \rnode{l'''}{l'''} & \rnode{l''}{l''}\hspace*{-.5ex}\pnode(0,.33){e0}
& \pnode(0,.33){e1} & \pnode(0,.33){e2} \raisebox{2ex}{\small$\,\cdots$} \pnode(0,.33){e3} 
& \pnode(0,.33){e4}\rnode{lstar+}{l^\star}^+
& \text{(depth $d-1$)}\\[2.5ex]
\vecanglesheight{l}{lstar+}{25}{95}{.4}\psset{nodesepA=5pt,nodesepB=0pt}\vecangles{lminus}{l'''}{40}{140}
\psset{nodesep=0pt}
\uloopvec{e0}{e1}
\uloopvec{e1}{e2}
\uloopvec{e3}{e4}
r: & \rnode{l}{l} & \rnode{l+}{l}^+  & \rnode{lminus}{l}^- & \rnode{l'''}{l'''} 
\uloopvec{l}{l+}\psset{nodesepA=5pt,nodesepB=0pt}
\vecangles{lminus}{l'''}{40}{140}
&&&&& \text{(depth $d-1$)}
\end{array}
\end{equation}
Note that $p'$ has depth $d-1$ since $p$ has depth $d$ and $l^{\star+}$ is
positive whereas $l^-$ is negative, and that $r'$ has depth $d-1$
since $r$ has depth $d-1$ and both $l^{\star+}$ and $l^+$ are positive.

The lax leaf function $p'$ is a lax linking (\ie, satisfies the linking
correctness criterion), 
since the targets $l^-$ and $l^{\star+}$ of the edge from $l$, in $p$ and $p'$
respectively, are connected along the edges of the lax leaf function.

\emph{Claim}: $r'$ is a lax linking. 
\emph{Proof.}  Let $e$ be the edge of $r'$ from $l$ to $l^{\star+}$, and 
let $e'$ be the edge of $r'$ from $l^-$ to $l'''$.
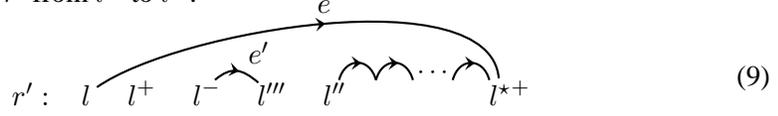
\begin{equation}
\addtolength{\arraycolsep}{.4ex}
\begin{array}{r@{\;\;\;\;}cccccccc@{\;\;\;\;\;\;\;\;\;\;}l}
\rule{0pt}{5ex}
r': & \rnode{l}{l} & \rnode{l+}{l}^+  & \rnode{lminus}{l}^- & \rnode{l'''}{l'''} & \rnode{l''}{l''}\hspace*{-.5ex}\pnode(0,.33){e0}
& \pnode(0,.33){e1} & \pnode(0,.33){e2} \raisebox{2ex}{\small$\,\cdots$} \pnode(0,.33){e3} 
& \pnode(0,.33){e4}\rnode{lstar+}{l^\star}^+
\vecanglesheight{l}{lstar+}{35}{95}{.5}\naput{e}\psset{nodesepA=5pt,nodesepB=0pt}\vecangles{lminus}{l'''}{40}{140}\naput[npos=.75]{e'}
\psset{nodesep=0pt}
\uloopvec{e0}{e1}
\uloopvec{e1}{e2}
\uloopvec{e3}{e4}
\end{array}
\end{equation}
Suppose $C$ is a cycle in a switching $\sigma$ of $r'$.  The cycle $C$
must traverse both $e$ and $e'$, for if it did not traverse $e'$ it
would be contained in the corresponding switching of $p'$ (already
proved to be a lax linking, and differing from $r'$ only on $l^-$)
whilst if it did not traverse $e$ it would be contained in the
corresponding switching of $r$ (differing from $r'$ only on $l$).
\begin{itemize}
\item[-] Case: $C$ (oriented one way or the other) has the form $e\pi e' \pi'$ for 
sequences of edges $\pi$ and $\pi'$. Let $e''$ be the edge from
$l$ to $l^-$ in $q$:
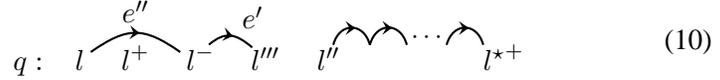
\begin{equation}
\addtolength{\arraycolsep}{.4ex}
\begin{array}{r@{\;\;\;\;}cccccccc}
\rule{0pt}{5ex}
q: & \rnode{l}{l} & \rnode{l+}{l}^+  & \rnode{lminus}{l}^- & \rnode{l'''}{l'''} & \rnode{l''}{l''}\hspace*{-.5ex}\pnode(0,.33){e0}
& \pnode(0,.33){e1} & \pnode(0,.33){e2} \raisebox{2ex}{\small$\,\cdots$} \pnode(0,.33){e3} 
& \pnode(0,.33){e4}\rnode{lstar+}{l^\star}^+
\vecangles{l}{lminus}{35}{145}\naput{e''}\psset{nodesepA=5pt,nodesepB=0pt}\vecangles{lminus}{l'''}{40}{140}\naput[npos=.75]{e'}
\psset{nodesep=0pt}
\uloopvec{e0}{e1}
\uloopvec{e1}{e2}
\uloopvec{e3}{e4}
\end{array}
\end{equation}
We obtain a cycle of the form $e''e'\pi'$ in the corresponding switching of $q$,
contradicting the fact that $q$ is a lax linking.
\item[-] Case: $C$ (oriented one way or the other) has the form $e\pi \overline{e}'\pi'$ for sequences 
of edges $\pi$ and $\pi'$, where $\overline{e}'$ is $e'$ traversed in
the direction opposite to its given orientation:
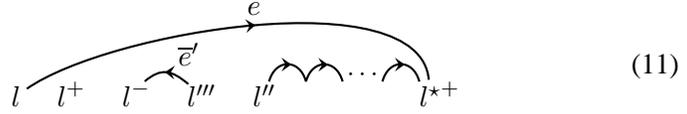
\begin{equation}
\addtolength{\arraycolsep}{.4ex}
\begin{array}{cccccccc}
\rule{0pt}{7ex}
\rnode{l}{l} & \rnode{l+}{l}^+  & \rnode{lminus}{l}^- & \rnode{l'''}{l'''} & \rnode{l''}{l''}\hspace*{-.5ex}\pnode(0,.33){e0}
& \pnode(0,.33){e1} & \pnode(0,.33){e2} \raisebox{2ex}{\small$\,\cdots$} \pnode(0,.33){e3} 
& \pnode(0,.33){e4}\rnode{lstar+}{l^\star}^+
\vecanglesheight{l}{lstar+}{35}{95}{.5}\naput{e}\psset{nodesepA=0pt,nodesepB=5pt}\vecangles{l'''}{lminus}{140}{40}
\nbput[npos=.25]{\overline{e}'}
\psset{nodesep=0pt}
\uloopvec{e0}{e1}
\uloopvec{e1}{e2}
\uloopvec{e3}{e4}
\end{array}
\end{equation}
Again, let $e''$ be the edge from
$l$ to $l^-$ in $q$:
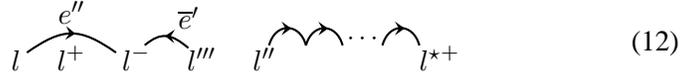
\begin{equation}
\addtolength{\arraycolsep}{.4ex}
\begin{array}{cccccccc}
\rule{0pt}{5ex}
\rnode{l}{l} & \rnode{l+}{l}^+  & \rnode{lminus}{l}^- & \rnode{l'''}{l'''} & \rnode{l''}{l''}\hspace*{-.5ex}\pnode(0,.33){e0}
& \pnode(0,.33){e1} & \pnode(0,.33){e2} \raisebox{2ex}{\small$\,\cdots$} \pnode(0,.33){e3} 
& \pnode(0,.33){e4}\rnode{lstar+}{l^\star}^+
\vecangles{l}{lminus}{35}{145}\naput{e''}\psset{nodesepA=0pt,nodesepB=5pt}\vecangles{l'''}{lminus}{140}{40}
\nbput[npos=.25]{\overline{e}'}
\psset{nodesep=0pt}
\uloopvec{e0}{e1}
\uloopvec{e1}{e2}
\uloopvec{e3}{e4}
\end{array}
\end{equation}
We obtain a cycle $e''\pi'$ in the corresponding switching of $q$,
contradicting the fact that $q$ is a lax linking.
\end{itemize}
\emph{QED Claim.}

By construction, the pairs $p\leftrightarrow p'$, $p'\leftrightarrow
r'$ and $r'\leftrightarrow r$ are similar, and $p'$ and $r'$ are
lax linkings.  Thus we can appeal to the inductive hypothesis with the
original sequence $h_1\ldots h_n$ with $p'r'$ substituted for $q$,
which has strictly smaller size than the original (since $p'$ and $r'$
each have depth $d-1$, whereas $q$ has depth $d$).
\end{itemize}
\end{itemize}
\end{itemize}
\end{proof}

\begin{lemma}\label{shift-to-arg}
Let $f$ be a lax linking, let $l$ be a negative $\unit$-labelled leaf
whose $f$-edge targets a vertex $v$ with an argument-vertex $a$ (thus
$v$ is either $\tensor$ or $\ast$).
Let $f'$ be the result of retargeting the edge from $l$ to point to
$a$ instead of $v$ (\ie, the lax leaf function $f'$ is $f$ but with
$f'(l)=a$ versus $f(l)=v$).
Then $f'$ is a lax linking.
\end{lemma}
\begin{proof}
If $v$ is a $\ast$-vertex, the result is immediate since switchings of
$f$ correspond to switchings of $f'$.  Thus assume $v$ is a
$\tensor$-vertex.  Without loss of generality, $a$ is the left
argument of $v$.  If $v$ is positive (hence retains both argument
edges in every switching), then the result is immediate.  So assume
$v$ is negative.\footnote{Note for readers familiar with linear
  distributivity in linear logic: the remainder of this proof is
  simply the usual argument (in disguise) that applying linear
  distributivity $\,A\tensor(B\,\footnoteparr\, C)\:\to\:(A\tensor
  B)\,\footnoteparr\, C\,$ preserves MLL proof net correctness.  The
  correspondence is as follows.  Let $S$ be the (sub)shape rooted at
  the negative $\tensor$-vertex $v$ in the proof in the main text;
  since $v$ is negative, think of it as a par.  Thus we think of $S$
  as $S_1\footnoteparr S_2$.  Substitute $I\tensor S$ for $S$, and
  retarget the $f$-edge from $l$ to point to the new $\unit$; call
  this the lax linking $\widehat{f}$. Now apply linear distributivity
  (and ignore $\ast$-vertices, which are irrelevant for cycles in
  switchings), yielding a lax linking $\widehat{f}'$ on $(I\tensor
  S_1)\footnoteparr S_2$.  This lax linking $\widehat{f}'$ corresponds
  to $f'$ just as $\widehat{f}$ corresponded to $f$; hence $f'$ is a
  well-defined lax linking.}  Towards a contradiction, suppose
$\sigma$ is a switching of $f'$ containing a cycle $C$.  Let $e$ be
the edge in $f$ from $l$ to $v$, let $e'$ be the edge in $f'$ from $l$
to $a$, let $e_1$ be the edge between $v$ and its left argument, and
let $e_2$ be the edge between $v$ and its right argument.
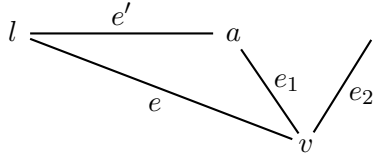
\begin{center}\vspace{2ex}
  \begin{math}
    %% \circlenode[linestyle=none]{l}{$l$}\hspace{6ex}
    %% \circlenode[linestyle=none[{a}{$a$}\hspace{3ex}
    %% \circlenode[linestyle=none[{b}{}\hspace{3ex}
    \rnode{l}{l\:\strut}\hspace{15ex}
    \rnode{a}{\:a\strut}\hspace{10ex}
    \rnode{b}{\;\strut}
    \ncline[linestyle=none]{a}{b}\nbput[labelsep=7ex]{\rnode{v}{\strut v}}
    \ncline[nodesep=0pt]{v}{a}\nbput[labelsep=.1ex]{e_1}
    \ncline[nodesep=0pt]{v}{b}\nbput[labelsep=.3ex]{e_2}
    \ncline{l}{v}\nbput[labelsep=.5ex]{e}
    \ncline{l}{a}\naput[labelsep=.5ex]{e'}
  \end{math}\vspace*{9ex}
\end{center}
%% \hspace{-6ex}\psset{colsep=2.5cm,framesep=0pt,nodesep=0pt,rowsep=1cm}
%% \begin{psmatrix}
%% \raisebox{5ex}{$
%% \pstree[treemode=U,nodesep=2pt,levelsep=20pt]
%% {\TR[name=t]{\tensor}}
%%   {\TR[name=b]{T}
%%    \TR[name=c]{U}}
%% \hspace{6ex}
%% \pstree[treemode=U,nodesep=2pt,levelsep=20pt]{\TR[name=ast]{\ast}}
%%   {\pstree{\TR[name=p]{\tensor}}
%%      {\TR[name=bb]{T}
%%       \TR[name=cc]{U}
%%      }
%%   }
%% $}
%% \nccurve[nodesepA=2pt,nodesepB=1pt,angleA=-40,angleB=-140]{t}{ast}
%% \end{psmatrix}
%
We may assume $C$ contains $e'$ (otherwise $C$ is a cycle in the
corresponding switching $\sigma_f$ of $f$) and $e_2$ (otherwise we can
assume $e_1$ is in $\sigma$ (by substituting $e_1$ for $e_2$ in
$\sigma$ if necessary), and we obtain a cycle in $\sigma_f$ by
replacing $e'$ in $C$ by $e$ and $e_1$).  Thus, oriented one way or
the other, $C$ has the form $e'\pi e_2\pi'$, where $e'$ is traversed
from $l$ to $a$, so $e\pi'$ (resp.\ $ee_2\pi'$) is a cycle if $e_2$ is
oriented towards (resp.\ away from) $v$.
%% (its natural orientation) and
%% the undirected argument edge $e_2$ may be traversed either (a) towards
%% $v$, or (b) away from $v$.  In case (a), $e\pi $ is a cycle in the
%% corresponding switching of $f$, and in case (b), $ee_2\pi'$ is a
%% cycle in the corresponding switching of $f$.
%; either case contradicts
%the fact that $f$ is a lax linking.
\end{proof}
\begin{corollary}\label{shift-to-leaf}
Let $f$ be a lax linking, let $l$ be a negative $\unit$-labelled leaf
whose $f$-edge targets a vertex $v$, and let $l'$ be one of the leaves
above $v$ (\ie, a leaf which is a hereditary argument of $v$).  The
result of retargeting the $f$-edge $l\to v$ to $l\to l'$ is a lax
linking.
\end{corollary}
\begin{proof}
Iterate Lemma~\ref{shift-to-arg}.
\end{proof}
\begin{lemma}\label{atomic-lemma}
Standard linkings are atomic equivalent iff they are lax equivalent.
\end{lemma}
\begin{proof}
A three-level induction.  Define the \defn{volume} of a shape as its
number of vertices.  The volume of a vertex $v$ in a shape is the
volume of the (sub)shape rooted at $v$.  For example, the volumes of
vertices have been subscripted on the following shape:
\raisebox{2pt}{\footnotesize$\big(\underset{1}{\unit}\,\underset{2}{{}\perp}\,\underset{3}{{}\perp}\,\underset{13}{\tensor}
((\underset{1}{a}\underset{4}{\tensor}\underset{1}{a}\,\underset{2}{{}\perp}\,)\underset{9}{\tensor}(\underset{1}{\unit}
\underset{3}{\tensor}\underset{1}{b})\,\underset{4}{{}\perp}\,)\big)\,\underset{14}{{}\perp}\,$}\:.
Define the volume of a negative $I$ in a lax linking $f$ as the volume
of its target under $f$, and the volume of $f$ as the maximum of the
volumes of its negative $I$'s.  Thus $f$ has volume $1$ iff it is an
atomic lax linking (\ie, every $f$-edge from a negative $I$ targets a
leaf).  The volume of a sequence $f_1\ldots f_m$ of lax linkings is
the maximum of the volumes of the $f_i$.  Let $V$ be the volume of
$f_1\ldots f_m$.  The \defn{depth} of $f_i$ (with respect to
$f_1\ldots f_m$) is the number of negative $I$'s in $f_i$ which have
volume $V$.  Define the depth of a sequence of linkings $f_1\ldots
f_m$ as the maximum of the depths of the $f_i$.  Suppose $f_1\ldots
f_m$ has depth $d$.  The
\defn{size} of $f_1\ldots f_m$ is the number of $f_i$ of depth $d$.

Let $f$ and $g$ be standard linkings, lax equivalent via a sequence
$f=h_1\ldots h_n=g$ of lax linkings with $h_i$ similar to $h_{i-1}$
for $1<i\le n$.  We proceed by a primary induction on the volume $V$ of
$h_1\ldots h_n$, a secondary induction on its depth $d$, and a
tertiary induction on its size $s$.

If $V=1$, then all the $h_i$ are already atomic lax linkings.

Suppose $V>1$.  Let $pqr$ be a subsequence of $h_1\ldots h_n$ in which
$q$ has depth $d>0$ (hence volume $V$) and $r$ has depth $e<d$.  (Such a
subsequence exists since $h_n=g$ has volume $1<V$.)  Let $c$ be the
depth of $p$.  There are two cases:
\begin{enumerate}
\item \emph{Case $c<d$.}  Since $c<d>e$ the rewirings $p\leftrightarrow q$ and $q\leftrightarrow r$ both rewire an edge
from the same negative $I$.  Delete $q$ from $h_1\ldots h_n$.  This
reduces at least one of the volume, depth or size of $h_1\ldots h_n$.
\item \emph{Case $c=d$.}  
If $l=l'$, then we can simply delete $q$, as in the previous case.
Thus assume $l\neq l'$.
We have
\begin{equation}
\addtolength{\arraycolsep}{.4ex}
\begin{array}{r@{\;\;\;\;}cccccc@{\;\;\;\;\;\;\;\;\;\;}l}
\rule{0pt}{4ex}
p: & \rnode{l}{l} & \rnode{v_1}{v_1}\hspace{-1.5ex} & & \rnode{l'}{l'}
& \rnode{w_1}{w_1}\hspace{-1.5ex} && \text{(depth $d$)}\\[5ex]
\uloopvec{l}{v_1}
\uloopvec{l'}{w_1}
q: & \rnode{l}{l} & \rnode{v_1}{v_1}\hspace{-1.5ex} & & \rnode{l'}{l'}
& & \rnode{w_2}{w_2} & \text{(depth $d$)}\\[4.2ex]
\uloopvec{l}{v_1}
\naput{V}
\vecangles{l'}{w_2}{35}{145}
r: & \rnode{l}{l} & & \rnode{v_2}{v_2} & \rnode{l'}{l'} & &
\rnode{w_2}{w_2}
\vecangles{l}{v_2}{35}{145}
\naput{<\mkern-4mu V}
\vecangles{l'}{w_2}{35}{145}
& \text{(depth $e<d$)}
\end{array}
\end{equation}
where $l$ and $l'$ are negative $I$-labelled leaves, the $v_i$ and
$w_j$ are vertices of unspecified type, and the volume of $l$ is $V$
in $p$ and $q$ and $<\mkern-4mu V$ in $r$.

Let $a$ be a leaf above $v_1$.  Define $p'$ and $q'$ from $p$ and $q$
by re-targeting the edge from $l$ to target $a$ instead of $v_1$.
Each of $p'$ and $q'$ is a well-defined lax linking by
Corollary~\ref{shift-to-leaf}.  Substitute $p'q'$ for $q$ in
$h_1\ldots h_n$.  This reduces at least one of the volume, depth or
size of $h_1\ldots h_n$, since the volume of $l$ in $p'$ and $q'$ is
strictly less than $V$.
\end{enumerate}\vspace{-4.5ex}
\end{proof}
\begin{proofof}{Lemma~\ref{lax-rewiring} 
(the Lax Rewiring Lemma: standard linkings are standard equivalent iff
they are lax equivalent)} By Lemma~\ref{lemma-reduce-atomic} standard
linkings are standard equivalent iff they are atomic equivalent, and
by Lemma~\ref{atomic-lemma} they are atomic equivalent iff they are
lax equivalent.
\end{proofof}

\subsection{Main freeness proof}

We are now ready to prove the Freeness Theorem (Theorem~\ref{free},
page~\pageref{free}): for any category $\A$, the category $\netsa$ of
$\A$-nets is the free star-autonomous category generated by $\A$.
Rather than prove the theorem from scratch, we show that $\netsa$ is
isomorphic to a full subcategory of the circuit category
$\polynetstara$ of \cite{BCST96}, where $(\mathcal{C}_\A,E_\A)$ is the
polygraph representing $\A$, with typed components $\mathcal{C}_\A$
and equations $E_\A$.  By Theorem~5.1 of \cite{BCST96},
$\polynetstara$ is the free symmetric linearly (=weakly) distributive
category with negation generated by $\A$.  Write $\netstara$ for the
full subcategory of $\polynetstara$ whose circuits are cotensor-free
and cotensor-unit-free (\ie\ par-free and $\bot$-free, in linear logic
terminology \cite{Gir87}).  Thus the objects of $\netstara$ are in
bijection with $\A$-shapes.  By the equivalence between symmetric
linearly distributive categories with negation and star-autonomous
categories in \cite{CS97}, $\netstara$ is the free star-autonomous
category generated by $\A$, as a corollary of Theorem~5.1 of
\cite{BCST96}.  Thus our Freeness Theorem, Theorem~\ref{free}, follows
from:
\begin{proposition}\label{iso-prop}
$\netsa$ is isomorphic to $\netstara$.
\end{proposition}
\subsection{Proof of Proposition~\ref{iso-prop}}

In this paper we have defined a \emph{net} as an equivalence class of
\emph{linkings}.  For $\netstara$ we shall use a similar two-level
convention: henceforth \defn{circuit-net} refers to an equivalence
class (a morphism of $\netstara$) and \defn{circuit} refers to a
representative.\footnote{Thus, in particular, we shall always assume a
  circuit satisfies the correctness criterion.}
%
%Given $\A$-shapes $S$ and $T$ define a \defn{normal circuit of type}
%$S\to T$ as any normalized\footnote{\cite{BCST96} was forced to work
%  with equivalence classes including un-normalized circuits since
%  thinning links could block redexes.} circuit representing a
%circuit-net in $\netstara$ of type $S\to T$ (with a single input wire
%of type $S$ and a single output wire of type $T$).  
For example, given morphisms $x:a\to b$ and $y:c\to d$ in $\A$, here
is a normal (\ie, redex-free\footnote{\cite{BCST96} was forced to work
  with equivalence classes including un-normalized circuits since
  thinning links could block redexes.}) circuit
$\unit\tensor\big((b\tensor c\perp)\perp\tensor I\big)\to (a\tensor
d\perp)\perp$:
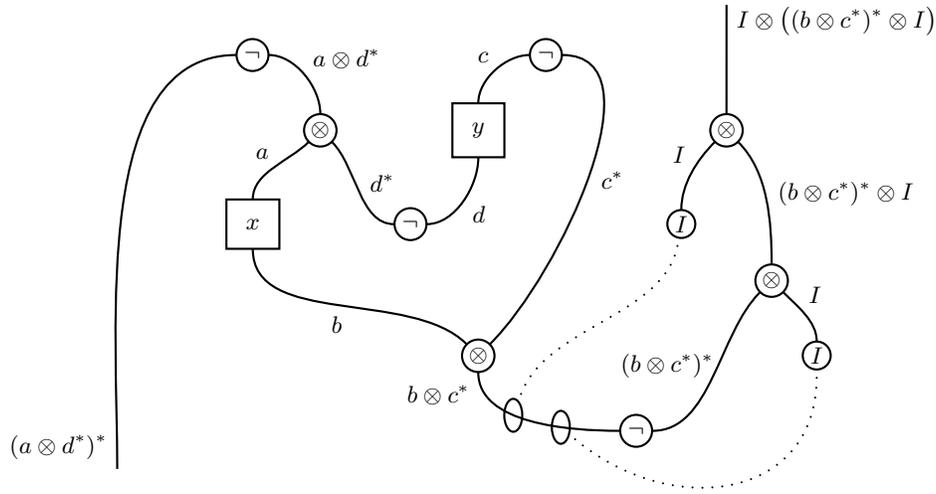
\begin{center}\footnotesize\psset{xunit=1.1cm,yunit=.5cm}
\newcommand{\blanklink}[3]{\pnode(#1,#2){#3}}
\newcommand{\unitlink}[3]{\cnodeput[framesep=1pt](#1,#2){#3}{$\unit$}}
\newcommand{\neglink}[3]{\cnodeput[framesep=2pt](#1,#2){#3}{$\neg$}}
\newcommand{\tlink}[3]{\cnodeput[framesep=1pt](#1,#2){#3}{$\tensor$}}
\newcommand{\boxlink}[4]{\rput(#1,#2){\rnode{#3}{\psframebox[framesep=7pt]{$#4$}}}}
\newcommand{\conc}{3}
\newcommand{\rowzero}{4}
\newcommand{\rowone}{6}
\newcommand{\rowtwo}{8}
\newcommand{\rowthree}{9.5}
\newcommand{\rowfour}{12}
\newcommand{\rowfive}{14}
\newcommand{\hypoth}{15.5}
\newcommand{\colone}{-4}
\newcommand{\coltwo}{2}
\newcommand{\colthree}{2}
\newcommand{\colfour}{5}
\newcommand{\colfive}{8}% redundant
\newcommand{\colsix}{9}
\newcommand{\colseven}{12}
\newcommand{\coleight}{15}
\newcommand{\colnine}{19}
\newcommand{\colten}{21}
\newcommand{\coleleven}{23}
\newcommand{\coltwelve}{25}
\newcommand{\colthirteen}{27}
\newcommand{\wire}[4]{\nccurve[nodesep=-.2pt,angleA=#3,angleB=#4]{#1}{#2}}
\newcommand{\down}[2]{\ncline[nodesepA=-.2pt,angleA=-90,angleB=90]{#1}{#2}}
\newcommand{\up}[2]{\ncline[nodesepA=-.2pt,angleA=90,angleB=-90]{#1}{#2}}
\newcommand{\leftdown}[2]{\wire{#1}{#2}{180}{90}}
\newcommand{\rightdown}[2]{\wire{#1}{#2}{0}{90}}
\newcommand{\leftup}[2]{\wire{#1}{#2}{180}{-90}}
\newcommand{\rightup}[2]{\wire{#1}{#2}{0}{-90}}
\newcommand{\leftdiagdown}[2]{\wire{#1}{#2}{-135}{90}}
\newcommand{\rightdiagdown}[2]{\wire{#1}{#2}{-45}{90}}
\newcommand{\leftdiagup}[2]{\wire{#1}{#2}{135}{-90}}
\newcommand{\rightdiagup}[2]{\wire{#1}{#2}{45}{-90}}
\begin{pspicture}(2,3)(4,15)
\psset{xunit=.3cm,yunit=1}
\blanklink{\colone}{\conc}{conc}
\neglink{\coltwo}{\rowfive}{neg}
\tlink{\colfour}{\rowfour}{tensor}
\boxlink{\colthree}{\rowthree}{x}{x}
\neglink{\colsix}{\rowthree}{neg1}
\boxlink{\colseven}{\rowfour}{y}{y}
\neglink{\coleight}{\rowfive}{neg2}
\tlink{\colseven}{\rowone}{tensor2}
\neglink{\colnine}{\rowzero}{neg3}
\tlink{\coltwelve}{\rowtwo}{tensor3}
\unitlink{\colthirteen}{\rowone}{I}
\tlink{\coleleven}{\rowfour}{tensor4}
\unitlink{\colten}{\rowthree}{J}
\blanklink{\coleleven}{\hypoth}{hypoth}
\leftdown{neg}{conc}
\nbput[npos=.95]{$(a\tensor d\perp)\perp$}
\rightdown{neg}{tensor}
\naput[npos=.5]{$a\tensor d\perp$}
\leftdiagdown{tensor}{x}
\nbput[npos=.5]{$a$}
\wire{tensor}{neg1}{-45}{180}
\naput[npos=.5]{$d\perp$}
\rightup{neg1}{y}
\nbput[npos=.5]{$d$}
\leftdown{neg2}{y}
\nbput[npos=.5]{$c$}
\wire{neg2}{tensor2}{0}{45}
\naput[npos=.5]{$c\perp$}
\leftdiagup{tensor2}{x}
\naput[npos=.5]{$b$}
\leftup{neg3}{tensor2}
\naput[npos=.9]{$b\tensor c\perp$}
\ncput[npos=.33]{\ovalnode{j1}{\rule{0pt}{4pt}}}
\ncput[npos=.6]{\ovalnode{j2}{\rule{0pt}{4pt}}}
\wire{neg3}{tensor3}{0}{-135}
\naput[npos=.5]{$\!\!(b\tensor c\perp)\perp\!\!\!$}
\rightdiagdown{tensor3}{I}
\naput[npos=.4]{$\unit$}
\rightdiagdown{tensor4}{tensor3}
\naput[npos=.5]{$(b\tensor c\perp)\perp\tensor\unit$}
\leftdiagdown{tensor4}{J}
\nbput[npos=.4]{$\unit$}
\up{tensor4}{hypoth}
\nbput[npos=.85]{$\unit\tensor\big((b\tensor c\perp)\perp\tensor\unit\big)$}
\psset{linestyle=dotted}
\nccurve[angleA=-100,angleB=70]{J}{j2}
\nccurve[angleA=-90,angleB=-40,ncurv=1]{I}{j1}
\end{pspicture}
\end{center}
(We write $\unit$ for the tensor unit, denoted $\top$ in
\cite{BCST96}.)  Define a \defn{canonical circuit} as a normal circuit
modulo the ordering of thinning links attached along each wire.  For
example, the following normal circuit denotes the same canonical
circuit as the normal circuit above:
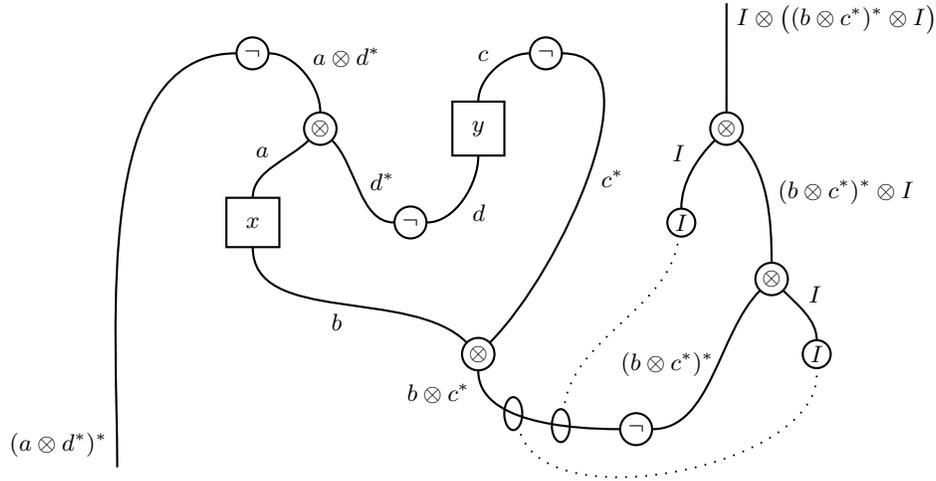
\begin{center}\footnotesize\psset{xunit=1.1cm,yunit=.5cm}
\newcommand{\blanklink}[3]{\pnode(#1,#2){#3}}
\newcommand{\unitlink}[3]{\cnodeput[framesep=1pt](#1,#2){#3}{$\unit$}}
\newcommand{\neglink}[3]{\cnodeput[framesep=2pt](#1,#2){#3}{$\neg$}}
\newcommand{\tlink}[3]{\cnodeput[framesep=1pt](#1,#2){#3}{$\tensor$}}
\newcommand{\boxlink}[4]{\rput(#1,#2){\rnode{#3}{\psframebox[framesep=7pt]{$#4$}}}}
\newcommand{\conc}{3}
\newcommand{\rowzero}{4}
\newcommand{\rowone}{6}
\newcommand{\rowtwo}{8}
\newcommand{\rowthree}{9.5}
\newcommand{\rowfour}{12}
\newcommand{\rowfive}{14}
\newcommand{\hypoth}{15.5}
\newcommand{\colone}{-4}
\newcommand{\coltwo}{2}
\newcommand{\colthree}{2}
\newcommand{\colfour}{5}
\newcommand{\colfive}{8}% redundant
\newcommand{\colsix}{9}
\newcommand{\colseven}{12}
\newcommand{\coleight}{15}
\newcommand{\colnine}{19}
\newcommand{\colten}{21}
\newcommand{\coleleven}{23}
\newcommand{\coltwelve}{25}
\newcommand{\colthirteen}{27}
\newcommand{\wire}[4]{\nccurve[nodesep=-.2pt,angleA=#3,angleB=#4]{#1}{#2}}
\newcommand{\down}[2]{\ncline[nodesepA=-.2pt,angleA=-90,angleB=90]{#1}{#2}}
\newcommand{\up}[2]{\ncline[nodesepA=-.2pt,angleA=90,angleB=-90]{#1}{#2}}
\newcommand{\leftdown}[2]{\wire{#1}{#2}{180}{90}}
\newcommand{\rightdown}[2]{\wire{#1}{#2}{0}{90}}
\newcommand{\leftup}[2]{\wire{#1}{#2}{180}{-90}}
\newcommand{\rightup}[2]{\wire{#1}{#2}{0}{-90}}
\newcommand{\leftdiagdown}[2]{\wire{#1}{#2}{-135}{90}}
\newcommand{\rightdiagdown}[2]{\wire{#1}{#2}{-45}{90}}
\newcommand{\leftdiagup}[2]{\wire{#1}{#2}{135}{-90}}
\newcommand{\rightdiagup}[2]{\wire{#1}{#2}{45}{-90}}
\begin{pspicture}(2,3)(4,15)
\psset{xunit=.3cm,yunit=1}
\blanklink{\colone}{\conc}{conc}
\neglink{\coltwo}{\rowfive}{neg}
\tlink{\colfour}{\rowfour}{tensor}
\boxlink{\colthree}{\rowthree}{x}{x}
\neglink{\colsix}{\rowthree}{neg1}
\boxlink{\colseven}{\rowfour}{y}{y}
\neglink{\coleight}{\rowfive}{neg2}
\tlink{\colseven}{\rowone}{tensor2}
\neglink{\colnine}{\rowzero}{neg3}
\tlink{\coltwelve}{\rowtwo}{tensor3}
\unitlink{\colthirteen}{\rowone}{I}
\tlink{\coleleven}{\rowfour}{tensor4}
\unitlink{\colten}{\rowthree}{J}
\blanklink{\coleleven}{\hypoth}{hypoth}
\leftdown{neg}{conc}
\nbput[npos=.95]{$(a\tensor d\perp)\perp$}
\rightdown{neg}{tensor}
\naput[npos=.5]{$a\tensor d\perp$}
\leftdiagdown{tensor}{x}
\nbput[npos=.5]{$a$}
\wire{tensor}{neg1}{-45}{180}
\naput[npos=.5]{$d\perp$}
\rightup{neg1}{y}
\nbput[npos=.5]{$d$}
\leftdown{neg2}{y}
\nbput[npos=.5]{$c$}
\wire{neg2}{tensor2}{0}{45}
\naput[npos=.5]{$c\perp$}
\leftdiagup{tensor2}{x}
\naput[npos=.5]{$b$}
\leftup{neg3}{tensor2}
\naput[npos=.9]{$b\tensor c\perp$}
\ncput[npos=.33]{\ovalnode{j1}{\rule{0pt}{4pt}}}
\ncput[npos=.6]{\ovalnode{j2}{\rule{0pt}{4pt}}}
\wire{neg3}{tensor3}{0}{-135}
\naput[npos=.5]{$\!\!(b\tensor c\perp)\perp\!\!\!$}
\rightdiagdown{tensor3}{I}
\naput[npos=.4]{$\unit$}
\rightdiagdown{tensor4}{tensor3}
\naput[npos=.5]{$(b\tensor c\perp)\perp\tensor\unit$}
\leftdiagdown{tensor4}{J}
\nbput[npos=.4]{$\unit$}
\up{tensor4}{hypoth}
\nbput[npos=.85]{$\unit\tensor\big((b\tensor c\perp)\perp\tensor\unit\big)$}
\psset{linestyle=dotted}
\nccurve[angleA=-100,angleB=80]{J}{j1}
\nccurve[angleA=-90,angleB=-70,ncurv=.6]{I}{j2}
\end{pspicture}
\end{center}
We render a canonical circuit uniquely by superimposing the attachment points
of thinning links, for example, drawing the above canonical circuit as
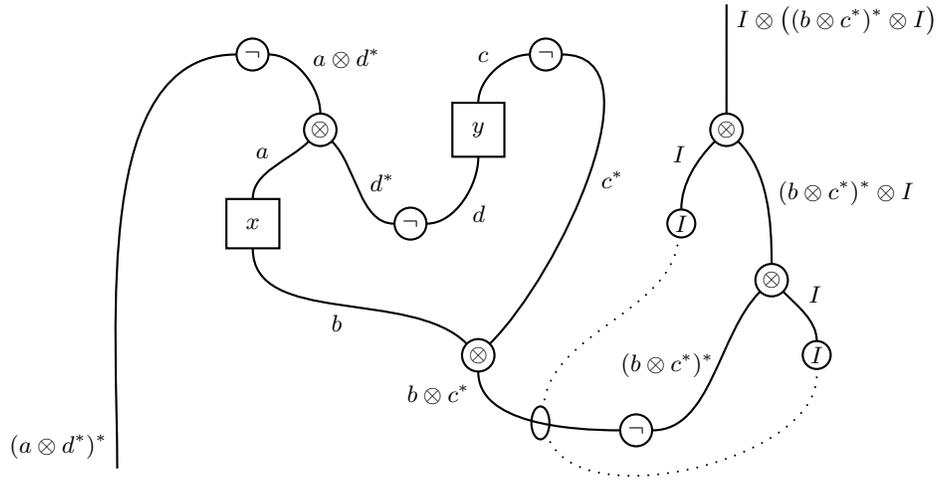
\begin{center}\label{canonical-circuit}\footnotesize\psset{xunit=1.1cm,yunit=.5cm}
\newcommand{\blanklink}[3]{\pnode(#1,#2){#3}}
\newcommand{\unitlink}[3]{\cnodeput[framesep=1pt](#1,#2){#3}{$\unit$}}
\newcommand{\neglink}[3]{\cnodeput[framesep=2pt](#1,#2){#3}{$\neg$}}
\newcommand{\tlink}[3]{\cnodeput[framesep=1pt](#1,#2){#3}{$\tensor$}}
\newcommand{\boxlink}[4]{\rput(#1,#2){\rnode{#3}{\psframebox[framesep=7pt]{$#4$}}}}
\newcommand{\conc}{3}
\newcommand{\rowzero}{4}
\newcommand{\rowone}{6}
\newcommand{\rowtwo}{8}
\newcommand{\rowthree}{9.5}
\newcommand{\rowfour}{12}
\newcommand{\rowfive}{14}
\newcommand{\hypoth}{15.5}
\newcommand{\colone}{-4}
\newcommand{\coltwo}{2}
\newcommand{\colthree}{2}
\newcommand{\colfour}{5}
\newcommand{\colfive}{8}% redundant
\newcommand{\colsix}{9}
\newcommand{\colseven}{12}
\newcommand{\coleight}{15}
\newcommand{\colnine}{19}
\newcommand{\colten}{21}
\newcommand{\coleleven}{23}
\newcommand{\coltwelve}{25}
\newcommand{\colthirteen}{27}
\newcommand{\wire}[4]{\nccurve[nodesep=-.2pt,angleA=#3,angleB=#4]{#1}{#2}}
\newcommand{\down}[2]{\ncline[nodesepA=-.2pt,angleA=-90,angleB=90]{#1}{#2}}
\newcommand{\up}[2]{\ncline[nodesepA=-.2pt,angleA=90,angleB=-90]{#1}{#2}}
\newcommand{\leftdown}[2]{\wire{#1}{#2}{180}{90}}
\newcommand{\rightdown}[2]{\wire{#1}{#2}{0}{90}}
\newcommand{\leftup}[2]{\wire{#1}{#2}{180}{-90}}
\newcommand{\rightup}[2]{\wire{#1}{#2}{0}{-90}}
\newcommand{\leftdiagdown}[2]{\wire{#1}{#2}{-135}{90}}
\newcommand{\rightdiagdown}[2]{\wire{#1}{#2}{-45}{90}}
\newcommand{\leftdiagup}[2]{\wire{#1}{#2}{135}{-90}}
\newcommand{\rightdiagup}[2]{\wire{#1}{#2}{45}{-90}}
\begin{pspicture}(2,3)(4,15)
\psset{xunit=.3cm,yunit=1}
\blanklink{\colone}{\conc}{conc}
\neglink{\coltwo}{\rowfive}{neg}
\tlink{\colfour}{\rowfour}{tensor}
\boxlink{\colthree}{\rowthree}{x}{x}
\neglink{\colsix}{\rowthree}{neg1}
\boxlink{\colseven}{\rowfour}{y}{y}
\neglink{\coleight}{\rowfive}{neg2}
\tlink{\colseven}{\rowone}{tensor2}
\neglink{\colnine}{\rowzero}{neg3}
\tlink{\coltwelve}{\rowtwo}{tensor3}
\unitlink{\colthirteen}{\rowone}{I}
\tlink{\coleleven}{\rowfour}{tensor4}
\unitlink{\colten}{\rowthree}{J}
\blanklink{\coleleven}{\hypoth}{hypoth}
\leftdown{neg}{conc}
\nbput[npos=.95]{$(a\tensor d\perp)\perp$}
\rightdown{neg}{tensor}
\naput[npos=.5]{$a\tensor d\perp$}
\leftdiagdown{tensor}{x}
\nbput[npos=.5]{$a$}
\wire{tensor}{neg1}{-45}{180}
\naput[npos=.5]{$d\perp$}
\rightup{neg1}{y}
\nbput[npos=.5]{$d$}
\leftdown{neg2}{y}
\nbput[npos=.5]{$c$}
\wire{neg2}{tensor2}{0}{45}
\naput[npos=.5]{$c\perp$}
\leftdiagup{tensor2}{x}
\naput[npos=.5]{$b$}
\leftup{neg3}{tensor2}
\naput[npos=.9]{$b\tensor c\perp$}
\ncput[npos=.44]{\ovalnode{j1}{\rule{0pt}{4pt}}}
\wire{neg3}{tensor3}{0}{-135}
\naput[npos=.5]{$\!\!(b\tensor c\perp)\perp\!\!\!$}
\rightdiagdown{tensor3}{I}
\naput[npos=.4]{$\unit$}
\rightdiagdown{tensor4}{tensor3}
\naput[npos=.5]{$(b\tensor c\perp)\perp\tensor\unit$}
\leftdiagdown{tensor4}{J}
\nbput[npos=.4]{$\unit$}
\up{tensor4}{hypoth}
\nbput[npos=.85]{$\unit\tensor\big((b\tensor c\perp)\perp\tensor\unit\big)$}
\psset{linestyle=dotted}
\nccurve[angleA=-100,angleB=80]{J}{j1}
\nccurve[angleA=-90,angleB=-63,ncurv=.6]{I}{j1}
\end{pspicture}
\end{center}
\begin{lemma}\label{lax-bijection}
Lax linkings $S\to T$ are in bijection with canonical circuits $S\to
T$.
\end{lemma}
\begin{proof}
Deleting the thinning links and generators from a canonical circuit
$S\to T$ leaves the parse tree structures of $S$ and $T$.  For
example, the circuit above leaves
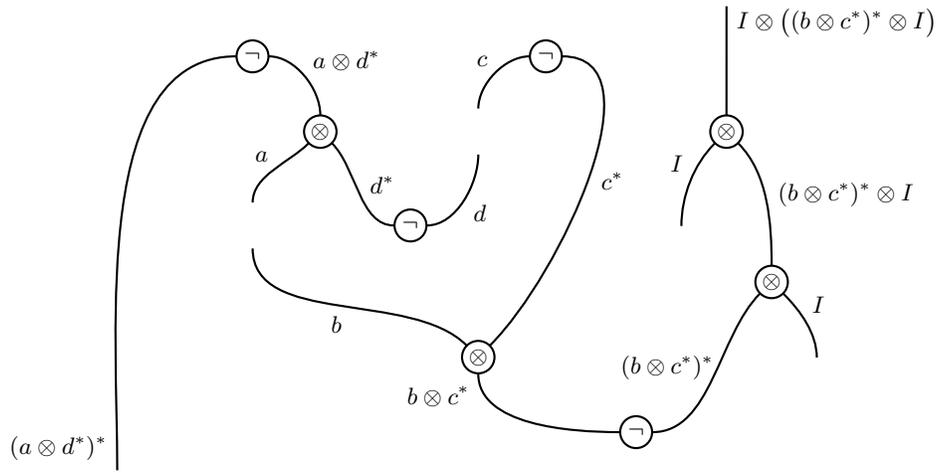
\begin{center}\footnotesize\psset{xunit=1.1cm,yunit=.5cm}
\newcommand{\blanklink}[3]{\pnode(#1,#2){#3}}
\newcommand{\unitlink}[3]{\cnodeput[framesep=1pt](#1,#2){#3}{$\unit$}}
\newcommand{\neglink}[3]{\cnodeput[framesep=2pt](#1,#2){#3}{$\neg$}}
\newcommand{\tlink}[3]{\cnodeput[framesep=1pt](#1,#2){#3}{$\tensor$}}
\newcommand{\boxlink}[4]{\cnodeput[linestyle=none,framesep=2ex](#1,#2){#3}{}}
\newcommand{\conc}{3}
\newcommand{\rowzero}{4}
\newcommand{\rowone}{6}
\newcommand{\rowtwo}{8}
\newcommand{\rowthree}{9.5}
\newcommand{\rowfour}{12}
\newcommand{\rowfive}{14}
\newcommand{\hypoth}{15.5}
\newcommand{\colone}{-4}
\newcommand{\coltwo}{2}
\newcommand{\colthree}{2}
\newcommand{\colfour}{5}
\newcommand{\colfive}{8}% redundant
\newcommand{\colsix}{9}
\newcommand{\colseven}{12}
\newcommand{\coleight}{15}
\newcommand{\colnine}{19}
\newcommand{\colten}{21}
\newcommand{\coleleven}{23}
\newcommand{\coltwelve}{25}
\newcommand{\colthirteen}{27}
\newcommand{\wire}[4]{\nccurve[nodesep=-.2pt,angleA=#3,angleB=#4]{#1}{#2}}
\newcommand{\down}[2]{\ncline[nodesepA=-.2pt,angleA=-90,angleB=90]{#1}{#2}}
\newcommand{\up}[2]{\ncline[nodesepA=-.2pt,angleA=90,angleB=-90]{#1}{#2}}
\newcommand{\leftdown}[2]{\wire{#1}{#2}{180}{90}}
\newcommand{\rightdown}[2]{\wire{#1}{#2}{0}{90}}
\newcommand{\leftup}[2]{\wire{#1}{#2}{180}{-90}}
\newcommand{\rightup}[2]{\wire{#1}{#2}{0}{-90}}
\newcommand{\leftdiagdown}[2]{\wire{#1}{#2}{-135}{90}}
\newcommand{\rightdiagdown}[2]{\wire{#1}{#2}{-45}{90}}
\newcommand{\leftdiagup}[2]{\wire{#1}{#2}{135}{-90}}
\newcommand{\rightdiagup}[2]{\wire{#1}{#2}{45}{-90}}
\begin{pspicture}(2,3)(4,15)
\psset{xunit=.3cm,yunit=1}
\blanklink{\colone}{\conc}{conc}
\neglink{\coltwo}{\rowfive}{neg}
\tlink{\colfour}{\rowfour}{tensor}
\boxlink{\colthree}{\rowthree}{x}{x}
\neglink{\colsix}{\rowthree}{neg1}
\boxlink{\colseven}{\rowfour}{y}{y}
\neglink{\coleight}{\rowfive}{neg2}
\tlink{\colseven}{\rowone}{tensor2}
\neglink{\colnine}{\rowzero}{neg3}
\tlink{\coltwelve}{\rowtwo}{tensor3}
\blanklink{\colthirteen}{\rowone}{I}
\tlink{\coleleven}{\rowfour}{tensor4}
\blanklink{\colten}{\rowthree}{J}
\blanklink{\coleleven}{\hypoth}{hypoth}
\leftdown{neg}{conc}
\nbput[npos=.95]{$(a\tensor d\perp)\perp$}
\rightdown{neg}{tensor}
\naput[npos=.5]{$a\tensor d\perp$}
\leftdiagdown{tensor}{x}
\nbput[npos=.5]{$a$}
\wire{tensor}{neg1}{-45}{180}
\naput[npos=.5]{$d\perp$}
\rightup{neg1}{y}
\nbput[npos=.5]{$d$}
\leftdown{neg2}{y}
\nbput[npos=.5]{$c$}
\wire{neg2}{tensor2}{0}{45}
\naput[npos=.5]{$c\perp$}
\leftdiagup{tensor2}{x}
\naput[npos=.5]{$b$}
\leftup{neg3}{tensor2}
\naput[npos=.9]{$b\tensor c\perp$}
\wire{neg3}{tensor3}{0}{-135}
\naput[npos=.5]{$\!\!(b\tensor c\perp)\perp\!\!\!$}
\rightdiagdown{tensor3}{I}
\naput[npos=.4]{$\unit$}
\rightdiagdown{tensor4}{tensor3}
\naput[npos=.5]{$(b\tensor c\perp)\perp\tensor\unit$}
\leftdiagdown{tensor4}{J}
\nbput[npos=.4]{$\unit$}
\up{tensor4}{hypoth}
\nbput[npos=.85]{$\unit\tensor\big((b\tensor c\perp)\perp\tensor\unit\big)$}
\end{pspicture}
\end{center}
The parse tree relationship is clearer if we (1) abstract away from
the dependency of the distinction between input and output
wires/ports\footnote{Lambek's \emph{covariables} and \emph{variables}: \cite[\S2.1]{BCST96}.} on the up/down direction in the page, by
explicitly orienting the wires from output to input, and (2) change
the $\neg$ label to $\ast$, to match the shape syntax.  For example,
the canonical circuit above (prior to deleting the generators and
thinning links) becomes the circuit below-left,
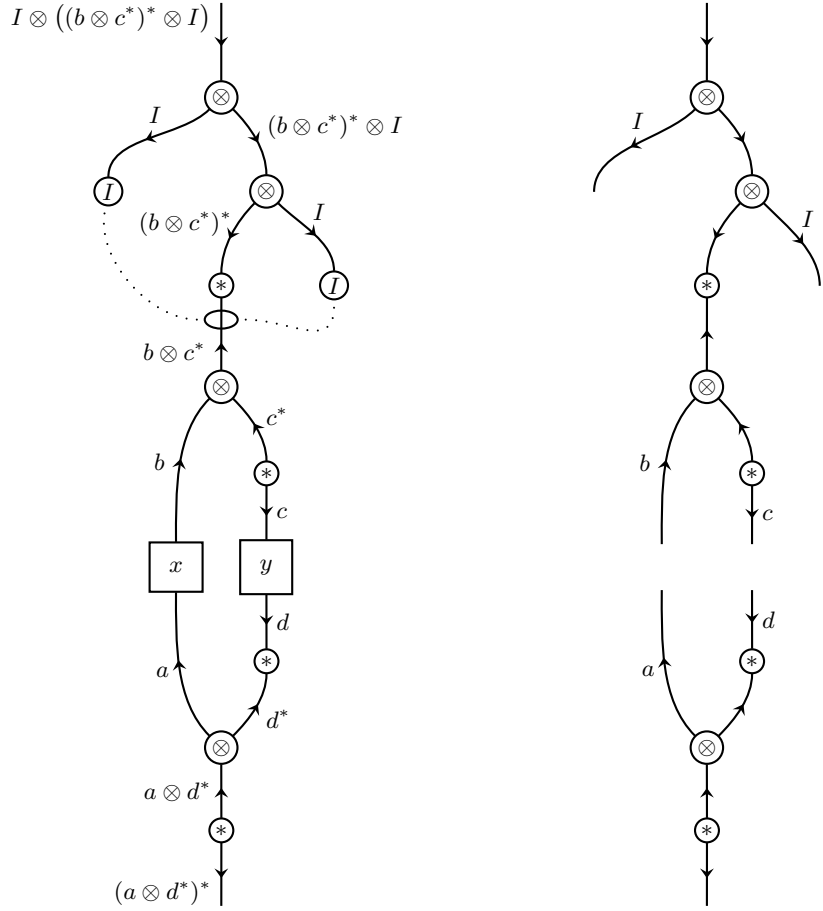
\begin{center}\label{vertical-circuit}\footnotesize\psset{xunit=1.1cm,yunit=.5cm}
\newcommand{\blanklink}[3]{\pnode(#1,#2){#3}}
\newcommand{\unitlink}[3]{\cnodeput[framesep=1pt](#1,#2){#3}{$\unit$}}
\newcommand{\neglink}[3]{\cnodeput[framesep=1pt](#1,#2){#3}{$\ast$}}
\newcommand{\tlink}[3]{\cnodeput[framesep=1pt](#1,#2){#3}{$\tensor$}}
\newcommand{\boxlink}[4]{\rput(#1,#2){\rnode{#3}{\psframebox[framesep=7pt]{$#4$}}}}
\newcommand{\conc}{3}
\newcommand{\rowone}{5}
\newcommand{\rowtwo}{7.2}
\newcommand{\rowthree}{9.5}
\newcommand{\rowfour}{12}
\newcommand{\rowfive}{14.5}
\newcommand{\rowsix}{16.8}
\newcommand{\rowseven}{19.5}
\newcommand{\roweight}{22}
\newcommand{\rownine}{24.5}
\newcommand{\hypoth}{27}
\newcommand{\colone}{-1}
\newcommand{\coltwo}{2}
\newcommand{\colthree}{4}
\newcommand{\colfour}{6}
\newcommand{\colfive}{9}
\newcommand{\wirevecat}[5]{\nccurve[nodesep=-.2pt,ArrowInside=->,ArrowInsidePos=#5,angleA=#3,angleB=#4]{#1}{#2}}
\newcommand{\wire}[4]{\wirevecat{#1}{#2}{#3}{#4}{.55}}
\newcommand{\down}[2]{\wire{#1}{#2}{-90}{90}}
\newcommand{\up}[2]{\wire{#1}{#2}{90}{-90}}
\newcommand{\upvecat}[3]{\wirevecat{#1}{#2}{90}{-90}{#3}}
\newcommand{\leftdown}[2]{\wire{#1}{#2}{180}{90}}
\newcommand{\rightdown}[2]{\wire{#1}{#2}{0}{90}}
\newcommand{\leftup}[2]{\wire{#1}{#2}{180}{-90}}
\newcommand{\rightup}[2]{\wire{#1}{#2}{0}{-90}}
\newcommand{\leftdiagdown}[2]{\wire{#1}{#2}{-135}{90}}
\newcommand{\rightdiagdown}[2]{\wire{#1}{#2}{-45}{90}}
\newcommand{\leftdiagup}[2]{\wire{#1}{#2}{135}{-90}}
\newcommand{\rightdiagup}[2]{\wire{#1}{#2}{45}{-90}}
\newcommand{\updiagright}[2]{\wire{#1}{#2}{90}{-135}}
\newcommand{\updiagleft}[2]{\wire{#1}{#2}{90}{-45}}
\newcommand{\diagleftdown}[2]{\wire{#1}{#2}{-135}{90}}
\newcommand{\diagrightdown}[2]{\wire{#1}{#2}{-45}{90}}
\vspace{-4ex}%
\begin{pspicture}(2,3)(.3,28)
\psset{xunit=.3cm,yunit=1}
\blanklink{\colthree}{\conc}{conc}
\neglink{\colthree}{\rowone}{neg}
\tlink{\colthree}{\rowtwo}{tensor}
\boxlink{\coltwo}{\rowfour}{x}{x}
\neglink{\colfour}{\rowthree}{neg1}
\boxlink{\colfour}{\rowfour}{y}{y}
\neglink{\colfour}{\rowfive}{neg2}
\tlink{\colthree}{\rowsix}{tensor2}
\neglink{\colthree}{\rowseven}{neg3}
\unitlink{\colfive}{\rowseven}{I}
\tlink{\colfour}{\roweight}{tensor3}
\unitlink{\colone}{\roweight}{J}
\tlink{\colthree}{\rownine}{tensor4}
\blanklink{\colthree}{\hypoth}{hypoth}
\down{neg}{conc}
\nbput[npos=.8]{$(a\tensor d\perp)\perp$}
\up{neg}{tensor}
\naput[npos=.5]{$a\tensor d\perp$}
\leftdiagup{tensor}{x}
\naput[npos=.5]{$a$}
\rightdiagup{tensor}{neg1}
\nbput[npos=.5]{$d\perp$}
\down{y}{neg1}
\naput[npos=.5]{$d$}
\down{neg2}{y}
\naput[npos=.5]{$c$}
\updiagleft{neg2}{tensor2}
\nbput[npos=.5]{$c\perp$}
\updiagright{x}{tensor2}
\naput[npos=.5]{$b\,$}
\upvecat{tensor2}{neg3}{.38}
\naput[npos=.25]{$b\tensor c\perp\:$}
\ncput[npos=.7]{\ovalnode{j1}{\rule{4pt}{0pt}}}
\diagleftdown{tensor3}{neg3}
\nbput[npos=.5]{$\!\!(b\tensor c\perp)\perp\!\!\!$}
\diagrightdown{tensor3}{I}
\naput[npos=.4]{$\unit$}
\diagrightdown{tensor4}{tensor3}
\naput[npos=.5]{$(b\tensor c\perp)\perp\tensor\unit$}
\diagleftdown{tensor4}{J}
\nbput[npos=.4]{$\unit$}
\down{hypoth}{tensor4}
\nbput[npos=.2]{$\unit\tensor\big((b\tensor c\perp)\perp\tensor\unit\big)$}
\psset{linestyle=dotted}
\nccurve[angleA=-100,angleB=180]{J}{j1}
\nccurve[angleA=-90,angleB=0,ncurv=1]{I}{j1}
\end{pspicture}
\hspace{30ex}
\footnotesize\psset{xunit=1.1cm,yunit=.5cm}
\renewcommand{\boxlink}[4]{\cnodeput[linestyle=none,framesep=2ex](#1,#2){#3}{}}
\begin{pspicture}(2,3)(.3,28)
\psset{xunit=.3cm,yunit=1}
\blanklink{\colthree}{\conc}{conc}
\neglink{\colthree}{\rowone}{neg}
\tlink{\colthree}{\rowtwo}{tensor}
\boxlink{\coltwo}{\rowfour}{x}{x}
\neglink{\colfour}{\rowthree}{neg1}
\boxlink{\colfour}{\rowfour}{y}{y}
\neglink{\colfour}{\rowfive}{neg2}
\tlink{\colthree}{\rowsix}{tensor2}
\neglink{\colthree}{\rowseven}{neg3}
\blanklink{\colfive}{\rowseven}{I}
\tlink{\colfour}{\roweight}{tensor3}
\blanklink{\colone}{\roweight}{J}
\tlink{\colthree}{\rownine}{tensor4}
\blanklink{\colthree}{\hypoth}{hypoth}
\down{neg}{conc}
\up{neg}{tensor}
\leftdiagup{tensor}{x}
\naput[npos=.5]{$a$}
\rightdiagup{tensor}{neg1}
\down{y}{neg1}
\naput[npos=.5]{$d$}
\down{neg2}{y}
\naput[npos=.5]{$c$}
\updiagleft{neg2}{tensor2}
\updiagright{x}{tensor2}
\naput[npos=.5]{$b\,$}
\up{tensor2}{neg3}
\diagleftdown{tensor3}{neg3}
\diagrightdown{tensor3}{I}
\naput[npos=.4]{$\unit$}
\diagrightdown{tensor4}{tensor3}
\diagleftdown{tensor4}{J}
\nbput[npos=.4]{$\unit$}
\down{hypoth}{tensor4}
\end{pspicture}
\end{center}
and deleting the generators $x$ and $y$, the thinning links, and
non-atomic labels more obviously leaves the parse trees of the shapes
$S=\unit\tensor\big((b\tensor c\perp)\perp\tensor\unit\big)$ and
$T=(a\tensor d\perp)\perp$, as shown above-right.
Generators and thinning links can then be viewed as the labelled edges
and unlabelled edges of a lax linking $S\to T$.  For example, the
canonical circuit above-left becomes the lax linking
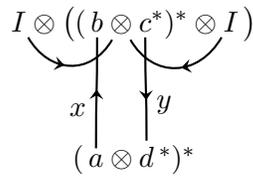
\begin{center}\begin{math}
\begin{array}{c}
\begin{psmatrix}[rowsep=3.3\baselineskip]
\rnode{I1}{\unit}\tensor\big((\,\rnode{b}{b}\:\rnode{t}{\tensor}\:\rnode{c}{c}\perp)\perp\tensor\rnode{I2}{\unit}\,\big) \\[8ex]
(\,\rnode{a}{a}\tensor\rnode{d}{d}\,{}\perp){}\perp
\vecanglespos{I1}{t}{-60}{-120}{.45}
\vecanglespos{I2}{t}{-120}{-60}{.45}
\uvecpos{a}{b}{.53}\naput[npos=.3]{x}
\dvecpos{c}{d}{.6}\naput[npos=.7]{y}
\end{psmatrix}\end{array}
\end{math}\end{center}\vspace{-4ex}
\end{proof}
Equivalence between canonical circuits is well-defined since
re-ordering of thinning links along a wire is a sub-relation of
circuit equivalence.
% Canonical circuits $C$ and $D$ are
%\defn{equivalent}\label{canonical-equivalent} if normal circuits $C'$
%and $D'$ representing $C$ and $D$, respectively, are equivalent.
%(This is independent of the choice of $C'$ and $D'$ since all normal
%circuits representing a canonical circuit are equivalent.)
\begin{lemma}\label{equiv-corresponds}
Lax linkings are lax equivalent iff the corresponding canonical
circuits are equivalent.
\end{lemma}
\begin{proof}
By the Empire Rewiring Proposition \cite[Prop.~3.3]{BCST96}, a
thinning link can be moved to any wire in its empire.  Since (by
definition) we are only dealing with circuits satisfying the
correctness criterion, such moves correspond to arbitrary retargeting
of edges from negative $\unit$s, between (correct) circuits.
\end{proof}
\begin{lemma}\label{equiv-to-standard}
Every lax linking $S\to T$ is lax equivalent to a standard linking
$S\to T$.
\end{lemma}
\begin{proof}
Using Corollary~\ref{shift-to-leaf} we can re-target all the edges from
negative $\unit$s to target leaves.  Then, suppose an edge from a leaf
$l$ targets a negative leaf $l'$, and suppose the edge from $l'$
targets $l''$.  Shift the edge from $l$ to target $l''$ instead.
Iterating this procedure leads to all edges targeting positive leaves,
yielding a standard linking.
\end{proof}
For any lax linking $f$, write $\LaxlinkingToCanonicalCircuit{f}$ for
the corresponding canonical circuit (via the bijection of
Lemma~\ref{lax-bijection}).  This induces a bijection
between nets and circuit-nets:
%
%Write $\LaxlinkingToCanonicalCircuit{(\strut\;)}$ for the bijection of
%  Lemma~\ref{lax-bijection}, thus $\LaxlinkingToCanonicalCircuit{l}$
%  denotes the canonical circuit corresponding to the lax linking $l$.
%
\begin{proposition}\label{net-bijection}
There is a bijection
$\netToCircuitNet{(\strut\;)}\;:\;\netsa(S,T)\to\netstara(S,T)$, for all shapes $S,T$.
% taking a net $n:S\to T$ in $\netsa$ to a net $\netToCircuitNet{n}:S\to T$ in $\netstara$.
\end{proposition}
\begin{proof}
Define $\netToCircuitNet{[f]}=[\LaxlinkingToCanonicalCircuit{f}]$, the
equivalence class of the canonical circuit
$\LaxlinkingToCanonicalCircuit{f}$.
%  Suppose $n\in\netsa(S,T)$.  Let $f$ be a
%standard linking representing $n$.
%%, with canonical circuit $\LaxlinkingToCanonicalCircuit{f}$ via
%%Lemma~\ref{lax-bijection}.
%
%Define the circuit-net $\netToCircuitNet{n}\in\netstara(S,T)$ as the
%equivalence class of the canonical circuit $\LaxlinkingToCanonicalCircuit{f}$.
%
This is independent of the choice of $f$ by
Lemma~\ref{equiv-corresponds}, injective by Lemma~\ref{lax-rewiring},
and surjective by Lemma~\ref{equiv-to-standard}.
\end{proof}
%% %
%% of any normal net $m$ representing
%% % Choose a normal net $m$ representing
%% $\LaxlinkingToCanonicalCircuit{l}$.
%% % and define $\netToCircuitNet{n}\in\netstara(S,T)$ as the equivalence class of $m$.
%% %
%% This is well-defined with respect to choice of $l$ and $m$: any
%% alternative $l'$ to $l$ is standard equivalent to $l$, thus
%% $\LaxlinkingToCanonicalCircuit{l'}$ is equivalent to
%% $\LaxlinkingToCanonicalCircuit{l}$ (by Lemma~\ref{equiv-corresponds});
%% every representative normal net of $\LaxlinkingToCanonicalCircuit{k'}$
%% is equivalent to every representative normal net of
%% $\LaxlinkingToCanonicalCircuit{l}$ (by definition of equivalence on
%% canonical circuits, preceding Lemma~\ref{equiv-corresponds}).
%
%
%% By the Lax Rewiring Lemma (Lemma~\ref{lax-rewiring}), this yields an
%% injection from nets $S\to T$ in $\netsa$ to nets $S\to T$ in
%% $\netstara$.  The injection is surjective since every lax linking is
%% equivalent to a standard linking, by Lemma~\ref{equiv-to-standard}.
%
To obtain an isomorphism of categories $\netsa\iso\netstara$, completing
the proof of Proposition~\ref{iso-prop}, we must prove that
$\netToCircuitNet{(\strut\;)}$ is functorial.
\begin{lemma}
Suppose $f:S\to T$ and $g:T\to U$ are standard $\A$-linkings.  Let
$\widehat{f}$ and $\widehat{g}$ be normal circuits representing the
canonical circuits $\LaxlinkingToCanonicalCircuit{f}$ and
$\LaxlinkingToCanonicalCircuit{g}$.  Let
$\hspace{.3ex}\widehat{f}\cup\widehat{g}\hspace{.1ex}$ be the circuit
obtained by pasting $\widehat{f}$ and $\widehat{g}$ at the $T\!$ wire.
There is a strategy of reduction and rewiring of thinning links for
$\widehat{f}\cup\widehat{g}$ leading to a normal circuit
$\;\widehat{f}\diamond\widehat{g}\,:\,S\to U$ whose canonical circuit is
$\LaxlinkingToCanonicalCircuit{(\seqcomp f g)}$ (that of the composite
of $f$ and $g$ in $\linka$).
\end{lemma}
\begin{proof}
Let $f'$ be the one-sided linking on the two-shape sequent $S\perp,T$
corresponding to $f$, and let $g'$ be the one-sided linking on
$T\perp,U$ corresponding to $g$.  Let $f'\cup g'$ be the disjoint
union of $f'$ and $g'$ on the cut sequent $S\perp,\cutpair{T},U$.  Let
$h'$ be the one-sided linking on $S\perp,U$ corresponding to the
composite two-sided linking $h=\seqcomp f g$.  Thus $h'$ is the normal form
resulting from cut elimination on $f'\cup g'$.  We shall mimic the cut
elimination steps on $\widehat{f}\cup\widehat{g}$ as reductions mixed
with rewiring of thinning links.

First, perform all $\tensor$ and $(-)\perp$ eliminations on $f'\cup
g'$.  This leaves the same leaf function $f'\cup g'$ on the cut
sequent $S\perp,\cutpair{a_1},\ldots,\cutpair{a_n},U$ where
$a_1,\ldots,a_n$ are the labels of the leaves of $T$.  Since these
eliminations affect only the parse trees of the shapes (and not the
edges of the leaf function), they can be mimicked directly on
$\widehat{f}\cup\widehat{g}$, as the tensor and tensor-unit reductions
\cite[\S3.1.1]{BCST96}.  Let $\widehat{f_0}\cup\widehat{g_0}$ denote
the end result of these reductions.

The normalisation of $f'\cup g'$ on the cut sequent
$S\perp,\cutpair{a_1},\ldots,\cutpair{a_n},U$ finishes with atomic
eliminations.  (See the definition of atomic elimination on
page~\pageref{atom-cut} for discrete $\A$, and its generalisation to
an arbitrary category $\A$ towards the end of Section~\ref{arb}.)
These atomic eliminations have one of two forms: reduction of (a) a
cut pair $\cutpair{\unit}$ or (b) a cut pair $\cutpair{a}$ for $a$ an
object of $\A$.

Consider (a).  Let $l_1,\ldots, l_n$ be the leaves having an edge to
the left $\unit$ of $\cutpair{\unit}$ and let $l$ be the target of the
edge from the right $\unit$ of $\cutpair{\unit}$.  Thus elimination
deletes $\cutpair{\unit}$ and moves the edges from $l_i$ to target $l$.
In the circuit we have the corresponding redex:
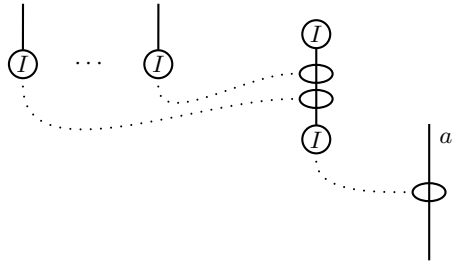
\newcommand{\localsettings}{\newcommand{\blanktop}{17}
\newcommand{\cdotsy}{13}
\newcommand{\lefti}{-3}
\newcommand{\righti}{3}
\newcommand{\topi}{15}
\newcommand{\boti}{8}
\newcommand{\barbellx}{10}
\newcommand{\wirex}{15}
\newcommand{\wiretop}{9}
\newcommand{\wirebot}{0}}
\begin{center}\footnotesize\psset{xunit=1.1cm,yunit=.5cm}
\newcommand{\blanklink}[3]{\pnode(#1,#2){#3}}
\newcommand{\unitlink}[3]{\cnodeput[framesep=1pt](#1,#2){#3}{$\unit$}}
\newcommand{\wire}[4]{\nccurve[nodesep=-.2pt,angleA=#3,angleB=#4]{#1}{#2}}
\newcommand{\down}[2]{\wire{#1}{#2}{-90}{90}}
\localsettings
\begin{pspicture}(3,0)(0,6.5)
\psset{xunit=.3cm,yunit=.4}
\blanklink{\lefti}{\blanktop}{blank1}
\blanklink{\righti}{\blanktop}{blank2}
\unitlink{\lefti}{\cdotsy}{i1}
\unitlink{\righti}{\cdotsy}{i2}
\down{blank1}{i1}
\down{blank2}{i2}
\ncline[linestyle=none]{i1}{i2}
\ncput{$\cdots$}
\unitlink{\barbellx}{\topi}{j1}
\unitlink{\barbellx}{\boti}{j2}
\down{j1}{j2}
\ncput[npos=.33]{\ovalnode{loop1}{\rule{4pt}{0pt}}}
\ncput[npos=.66]{\ovalnode{loop2}{\rule{4pt}{0pt}}}
\blanklink{\wirex}{\wiretop}{w1}
\blanklink{\wirex}{\wirebot}{w2}
\down{w1}{w2}
\ncput[npos=.5]{\ovalnode{loop3}{\rule{4pt}{0pt}}}
\naput[npos=.1]{$a$}
\psset{linestyle=dotted}
\nccurve[angleA=-90,angleB=180]{i1}{loop2}
\nccurve[angleA=-90,angleB=180]{i2}{loop1}
\nccurve[angleA=-90,angleB=180]{j2}{loop3}
\end{pspicture}
\end{center}
where the two left $\unit$-nodes tied by the ellipsis represent $n$
$\unit$-nodes corresponding to the leaves $l_1,\dots,l_n$, the highest
$\unit$-node corresponds to the left $\unit$ of $\cutpair{\unit}$, the
lowest $\unit$-node corresponds to the right $\unit$ of
$\cutpair{\unit}$, and $a$ is the label of $l$ (either $a=\unit$ or
$a$ is an object of $\A$).  Rewire the $n$ thinning links on the left
to target the $a$-wire:
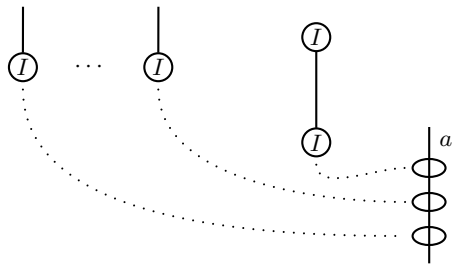
\begin{center}\footnotesize\psset{xunit=1.1cm,yunit=.5cm}
\newcommand{\blanklink}[3]{\pnode(#1,#2){#3}}
\newcommand{\unitlink}[3]{\cnodeput[framesep=1pt](#1,#2){#3}{$\unit$}}
\newcommand{\wire}[4]{\nccurve[nodesep=-.2pt,angleA=#3,angleB=#4]{#1}{#2}}
\newcommand{\down}[2]{\wire{#1}{#2}{-90}{90}}
\localsettings
\begin{pspicture}(3,0)(0,6.5)
\psset{xunit=.3cm,yunit=.4}
\blanklink{\lefti}{\blanktop}{blank1}
\blanklink{\righti}{\blanktop}{blank2}
\unitlink{\lefti}{\cdotsy}{i1}
\unitlink{\righti}{\cdotsy}{i2}
\down{blank1}{i1}
\down{blank2}{i2}
\ncline[linestyle=none]{i1}{i2}
\ncput{$\cdots$}
\unitlink{\barbellx}{\topi}{j1}
\unitlink{\barbellx}{\boti}{j2}
\down{j1}{j2}
\blanklink{\wirex}{\wiretop}{w1}
\blanklink{\wirex}{\wirebot}{w2}
\down{w1}{w2}
\ncput[npos=.3]{\ovalnode{loop3}{\rule{4pt}{0pt}}}
\ncput[npos=.55]{\ovalnode{loop1}{\rule{4pt}{0pt}}}
\ncput[npos=.8]{\ovalnode{loop2}{\rule{4pt}{0pt}}}
\naput[npos=.1]{$a$}
\psset{linestyle=dotted}
\nccurve[angleA=-90,angleB=180]{i1}{loop2}
\nccurve[angleA=-90,angleB=180]{i2}{loop1}
\nccurve[angleA=-90,angleB=180]{j2}{loop3}
\end{pspicture}
\end{center}
then reduce the $\unit$-node redex:
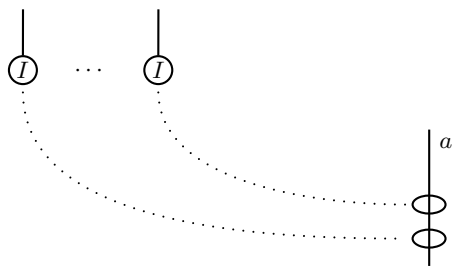
\begin{center}\footnotesize\psset{xunit=1.1cm,yunit=.5cm}
\newcommand{\blanklink}[3]{\pnode(#1,#2){#3}}
\newcommand{\unitlink}[3]{\cnodeput[framesep=1pt](#1,#2){#3}{$\unit$}}
\newcommand{\wire}[4]{\nccurve[nodesep=-.2pt,angleA=#3,angleB=#4]{#1}{#2}}
\newcommand{\down}[2]{\wire{#1}{#2}{-90}{90}}
\localsettings
\begin{pspicture}(3,0)(0,6.5)
\psset{xunit=.3cm,yunit=.4}
\blanklink{\lefti}{\blanktop}{blank1}
\blanklink{\righti}{\blanktop}{blank2}
\unitlink{\lefti}{\cdotsy}{i1}
\unitlink{\righti}{\cdotsy}{i2}
\down{blank1}{i1}
\down{blank2}{i2}
\ncline[linestyle=none]{i1}{i2}
\ncput{$\cdots$}
\blanklink{\wirex}{\wiretop}{w1}
\blanklink{\wirex}{\wirebot}{w2}
\down{w1}{w2}
\ncput[npos=.55]{\ovalnode{loop1}{\rule{4pt}{0pt}}}
\ncput[npos=.8]{\ovalnode{loop2}{\rule{4pt}{0pt}}}
\naput[npos=.1]{$a$}
\psset{linestyle=dotted}
\nccurve[angleA=-90,angleB=180]{i1}{loop2}
\nccurve[angleA=-90,angleB=180]{i2}{loop1}
\end{pspicture}
\end{center}
Case (b) is similar.  The circuit has a corresponding redex
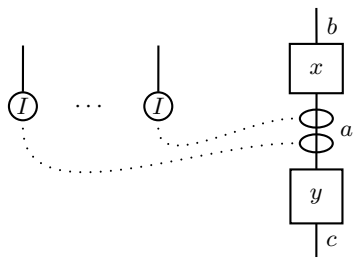
\begin{center}\footnotesize\psset{xunit=1.1cm,yunit=.5cm}
\newcommand{\blanklink}[3]{\pnode(#1,#2){#3}}
\newcommand{\unitlink}[3]{\cnodeput[framesep=1pt](#1,#2){#3}{$\unit$}}
\newcommand{\boxlink}[4]{\rput(#1,#2){\rnode{#3}{\psframebox[framesep=7pt]{$#4$}}}}
\newcommand{\wire}[4]{\nccurve[nodesep=-.2pt,angleA=#3,angleB=#4]{#1}{#2}}
\newcommand{\down}[2]{\wire{#1}{#2}{-90}{90}}
\localsettings
\renewcommand{\topi}{15.5}
\renewcommand{\boti}{7}
\newcommand{\boxtop}{19.5}
\newcommand{\boxbot}{3}
\begin{pspicture}(3,0)(0,7.5)
\psset{xunit=.3cm,yunit=.4}
\blanklink{\lefti}{\blanktop}{blank1}
\blanklink{\righti}{\blanktop}{blank2}
\unitlink{\lefti}{\cdotsy}{i1}
\unitlink{\righti}{\cdotsy}{i2}
\down{blank1}{i1}
\down{blank2}{i2}
\ncline[linestyle=none]{i1}{i2}
\ncput{$\cdots$}
\blanklink{\barbellx}{\boxtop}{boxtop}
\boxlink{\barbellx}{\topi}{j1}{x}
\boxlink{\barbellx}{\boti}{j2}{y}
\blanklink{\barbellx}{\boxbot}{boxbot}
\down{boxtop}{j1}
\naput{$b$}
\down{j1}{j2}
\naput{$\;\;a$}
\ncput[npos=.33]{\ovalnode{loop1}{\rule{4pt}{0pt}}}
\ncput[npos=.66]{\ovalnode{loop2}{\rule{4pt}{0pt}}}
\down{j2}{boxbot}
\naput{$c$}
%% \blanklink{\wirex}{\wiretop}{w1}
%% \blanklink{\wirex}{\wirebot}{w2}
%% \down{w1}{w2}
%% \naput[npos=.1]{$a$}
\psset{linestyle=dotted}
\nccurve[angleA=-90,angleB=180]{i1}{loop2}
\nccurve[angleA=-90,angleB=180]{i2}{loop1}
\end{pspicture}
\end{center}
where $x:b\to a$ and $y:a\to c$ are morphisms in $\A$.  Shift the
thinning links,
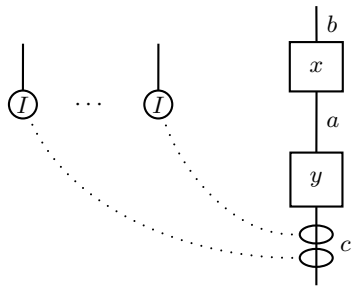
\begin{center}\footnotesize\psset{xunit=1.1cm,yunit=.5cm}
\newcommand{\blanklink}[3]{\pnode(#1,#2){#3}}
\newcommand{\unitlink}[3]{\cnodeput[framesep=1pt](#1,#2){#3}{$\unit$}}
\newcommand{\boxlink}[4]{\rput(#1,#2){\rnode{#3}{\psframebox[framesep=7pt]{$#4$}}}}
\newcommand{\wire}[4]{\nccurve[nodesep=-.2pt,angleA=#3,angleB=#4]{#1}{#2}}
\newcommand{\down}[2]{\wire{#1}{#2}{-90}{90}}
\localsettings
\renewcommand{\topi}{15.5}
\renewcommand{\boti}{8}
\newcommand{\boxtop}{19.5}
\newcommand{\boxbot}{1}
\begin{pspicture}(3,0)(0,7.5)
\psset{xunit=.3cm,yunit=.4}
\blanklink{\lefti}{\blanktop}{blank1}
\blanklink{\righti}{\blanktop}{blank2}
\unitlink{\lefti}{\cdotsy}{i1}
\unitlink{\righti}{\cdotsy}{i2}
\down{blank1}{i1}
\down{blank2}{i2}
\ncline[linestyle=none]{i1}{i2}
\ncput{$\cdots$}
\blanklink{\barbellx}{\boxtop}{boxtop}
\boxlink{\barbellx}{\topi}{j1}{x}
\boxlink{\barbellx}{\boti}{j2}{y}
\blanklink{\barbellx}{\boxbot}{boxbot}
\down{boxtop}{j1}
\naput{$b$}
\down{j1}{j2}
\naput{$a$}
\down{j2}{boxbot}
\naput{$\;\;c$}
%\blanklink{\wirex}{\wiretop}{w1}
%\blanklink{\wirex}{\wirebot}{w2}
%\down{w1}{w2}
%\naput[npos=.1]{$a$}
\ncput[npos=.35]{\ovalnode{loop1}{\rule{4pt}{0pt}}}
\ncput[npos=.65]{\ovalnode{loop2}{\rule{4pt}{0pt}}}
\psset{linestyle=dotted}
\nccurve[angleA=-65,angleB=-180]{i1}{loop2}
\nccurve[angleA=-70,angleB=-180]{i2}{loop1}
\end{pspicture}
\end{center}
then reduce:
\begin{center}\footnotesize\psset{xunit=1.1cm,yunit=.5cm}
\newcommand{\blanklink}[3]{\pnode(#1,#2){#3}}
\newcommand{\unitlink}[3]{\cnodeput[framesep=1pt](#1,#2){#3}{$\unit$}}
\newcommand{\boxlink}[4]{\rput(#1,#2){\rnode{#3}{\psframebox[framesep=7pt]{$#4$}}}}
\newcommand{\wire}[4]{\nccurve[nodesep=-.2pt,angleA=#3,angleB=#4]{#1}{#2}}
\newcommand{\down}[2]{\wire{#1}{#2}{-90}{90}}
\localsettings
\renewcommand{\topi}{15.5}
\renewcommand{\boti}{11}
\newcommand{\midi}{11.25}
\newcommand{\boxtop}{19.5}
\newcommand{\boxbot}{3}
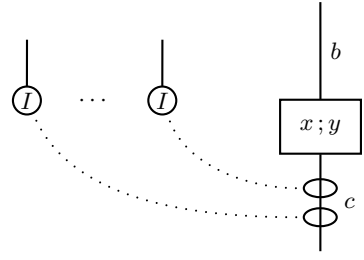
\begin{pspicture}(3,0)(0,7.5)
\psset{xunit=.3cm,yunit=.4}
\blanklink{\lefti}{\blanktop}{blank1}
\blanklink{\righti}{\blanktop}{blank2}
\unitlink{\lefti}{\cdotsy}{i1}
\unitlink{\righti}{\cdotsy}{i2}
\down{blank1}{i1}
\down{blank2}{i2}
\ncline[linestyle=none]{i1}{i2}
\ncput{$\cdots$}
\blanklink{\barbellx}{\boxtop}{boxtop}
\boxlink{\barbellx}{\midi}{yx}{\seqcomp{x\,} y}
\blanklink{\barbellx}{\boxbot}{boxbot}
\down{boxtop}{yx}
\naput{$b$}
\down{yx}{boxbot}
\naput{$\;\;c$}
%% \blanklink{\wirex}{\wiretop}{w1}
%% \blanklink{\wirex}{\wirebot}{w2}
%% \down{w1}{w2}
%% \naput[npos=.1]{$a$}
\ncput[npos=.35]{\ovalnode{loop1}{\rule{4pt}{0pt}}}
\ncput[npos=.65]{\ovalnode{loop2}{\rule{4pt}{0pt}}}
\psset{linestyle=dotted}
\nccurve[angleA=-65,angleB=-180]{i1}{loop2}
\nccurve[angleA=-70,angleB=-180]{i2}{loop1}
\end{pspicture}
\end{center}
where $\seqcomp{x\,} y\,:\,b\to c$ is the composite of the morphisms $x:b\to a$ and
$y:a\to c$ in $\A$.

Since the elimination steps are mimicked precisely, the resulting
normal circuit $\widehat{f}\diamond\widehat{g}$, modulo the order of
attachments of thinning links along the same wire (\ie, the canonical
circuit represented by $\widehat{f}\diamond\widehat{g}\,$), corresponds to the
normal one-sided linking $h'$, hence the composite two-sided linking
$h=\seqcomp f g$.
\end{proof}
\begin{corollary}
Let $f:S\to T$ and $g:T\to U$ be nets in $\netsa$.  Then
$\seqcomp{\netToCircuitNet{f}}{\netToCircuitNet{g}}\:=\:\netToCircuitNet{(\seqcomp f g)}\::\:S\to
U$ in $\netstara$.
%Let $f$ and $g$ be nets in $\netsa$ corresponding to the nets
%$\widehat{f}$ and $\widehat{g}$ in $\netstara$, by the bijection in
%Proposition~\ref{iso-prop}.  Then the composite $f;g$ in $\netsa$
%corresponds to the composite $\widehat{f};\widehat{g}$ in $\netstara$.
\end{corollary}
Thus the bijection $\netToCircuitNet{(\strut\;)}$ of
Proposition~\ref{net-bijection} preserves composition.  It preserves
identities because the identity linking and the identity circuit $S\to
S$ each amount to a dual pair of copies of the parse tree of $S$.
Therefore $\netToCircuitNet{(\strut\;)}\::\:\netsa\:\to\:\netstara$ is
functorial, providing an isomorphism of categories
$\netsa\iso\netstara$.
%% In other words, we have one half of the desired functoriality: the
%% bijection $\netToCircuitNet{(\strut\;)}$ of Proposition~\ref{iso-prop}
%% preserves composition.  Finally, $F(\id)=\id$ since both the identity
%% linking and the identity circuit $S\to S$ amount to dual copies of the
%% parse tree $S$.  
This completes the proof of Proposition~\ref{iso-prop}, whence the
Freeness Theorem: $\netsa$ is the free star-autonomous category
generated by $\A$.

\small
\bibliographystyle{alpha}

\end{document}